\documentclass[12pt,english,refpage,intoc,bibliography=totoc,index=totoc,BCOR7.5mm,captions=tableheading]{amsart}
\usepackage{lmodern}

\usepackage[T1]{fontenc}
\usepackage[latin9]{inputenc}
\usepackage{geometry}
\usepackage{color}
\usepackage{babel}
\usepackage{enumitem}
\usepackage{amstext}
\usepackage{amsthm}
\usepackage{amssymb}
\usepackage{cancel}
\PassOptionsToPackage{normalem}{ulem}
\usepackage{ulem}
\usepackage[unicode=true,
 bookmarks=true,bookmarksnumbered=true,bookmarksopen=false,
 breaklinks=false,pdfborder={0 0 1},backref=false,colorlinks=true]
 {hyperref}
\hypersetup{pdftitle={The LyX User's Guide},
 pdfauthor={LyX Team},
 pdfsubject={LyX},
 pdfkeywords={LyX},
 linkcolor=black, citecolor=black, urlcolor=blue, filecolor=blue, pdfpagelayout=OneColumn, pdfnewwindow=true, pdfstartview=XYZ, plainpages=false}

\makeatletter

\providecommand{\tabularnewline}{\\}

\numberwithin{equation}{section}
\numberwithin{figure}{section}
\theoremstyle{remark}
    \ifx\thechapter\undefined
      \newtheorem{note}{\protect\notename}
    \else
      \newtheorem{note}{\protect\notename}[chapter]
    \fi
\theoremstyle{plain}
    \ifx\thechapter\undefined
      \newtheorem{prop}{\protect\propositionname}
    \else
      \newtheorem{prop}{\protect\propositionname}[chapter]
    \fi
\theoremstyle{definition}
    \ifx\thechapter\undefined
      \newtheorem{defn}{\protect\definitionname}
    \else
      \newtheorem{defn}{\protect\definitionname}[chapter]
    \fi
\theoremstyle{plain}
    \ifx\thechapter\undefined
      \newtheorem{conjecture}{\protect\conjecturename}
    \else
      \newtheorem{conjecture}{\protect\conjecturename}[chapter]
    \fi
\theoremstyle{remark}
    \ifx\thechapter\undefined
      \newtheorem{rem}{\protect\remarkname}
    \else
      \newtheorem{rem}{\protect\remarkname}[chapter]
    \fi
\theoremstyle{definition}
    \ifx\thechapter\undefined
      \newtheorem{example}{\protect\examplename}
    \else
      \newtheorem{example}{\protect\examplename}[chapter]
    \fi
\theoremstyle{plain}
    \ifx\thechapter\undefined
	    \newtheorem{thm}{\protect\theoremname}
	  \else
      \newtheorem{thm}{\protect\theoremname}[chapter]
    \fi
\theoremstyle{plain}
    \ifx\thechapter\undefined
      \newtheorem{lem}{\protect\lemmaname}
    \else
      \newtheorem{lem}{\protect\lemmaname}[chapter]
    \fi
\theoremstyle{remark}
    \ifx\thechapter\undefined
      \newtheorem{claim}{\protect\claimname}
    \else
      \newtheorem{claim}{\protect\claimname}[chapter]
    \fi
\theoremstyle{plain}
    \ifx\thechapter\undefined
  \newtheorem{cor}{\protect\corollaryname}
\else
      \newtheorem{cor}{\protect\corollaryname}[chapter]
    \fi

\usepackage[figure]{hypcap}

\let\myTOC\tableofcontents
\renewcommand\tableofcontents{%
  \frontmatter
  \pdfbookmark[1]{\contentsname}{}
  \myTOC
  \mainmatter }

\usepackage{fancyhdr}

\@ifundefined{extrarowheight}
 {\usepackage{array}}{}
\makeatother

\providecommand{\claimname}{Claim}
\providecommand{\conjecturename}{Conjecture}
\providecommand{\corollaryname}{Corollary}
\providecommand{\definitionname}{Definition}
\providecommand{\examplename}{Example}
\providecommand{\lemmaname}{Lemma}
\providecommand{\notename}{Note}
\providecommand{\propositionname}{Proposition}
\providecommand{\remarkname}{Remark}
\providecommand{\theoremname}{Theorem}

\begin{document}
\title{The Collatz Conjecture \& Non-Archimedean Spectral Theory - Part I
- Arithmetic Dynamical Systems and Non-Archimedean Value Distribution
Theory}
\author{by Maxwell C. Siegel (University of Southern California)}
\address{3620 S. Vermont Ave., KAP 104, Los Angeles, CA 90089-2532}
\curraddr{1626 Thayer Ave., Los Angeles, CA, 90024}
\email{siegelmaxwellc@ucla.edu}
\date{27 March, 2023}
\begin{abstract}
Let $q$ be an odd prime, and let $T_{q}:\mathbb{Z}\rightarrow\mathbb{Z}$
be the Shortened $qx+1$ map, defined by $T_{q}\left(n\right)=n/2$
if $n$ is even and $T_{q}\left(n\right)=\left(qn+1\right)/2$ if
$n$ is odd. The study of the dynamics of these maps is infamous for
its difficulty, with the characterization of the dynamics of $T_{3}$
being an alternative formulation of the famous \textbf{Collatz Conjecture}.
This series of papers presents a new paradigm for studying such arithmetic
dynamical systems by way of a neglected area of ultrametric analysis
which we have termed \textbf{$\left(p,q\right)$-adic analysis}, the
study of functions from the $p$-adics to the $q$-adics, where $p$
and $q$ are distinct primes. In this, the first paper, working with
the $T_{q}$ maps as a toy model for the more general theory, for
each odd prime $q$, we construct a function $\chi_{q}:\mathbb{Z}_{2}\rightarrow\mathbb{Z}_{q}$
(the \textbf{Numen }of $T_{q}$) and prove the \textbf{Correspondence
Principle} (CP): $x\in\mathbb{Z}\backslash\left\{ 0\right\} $ is
a periodic point of $T_{q}$ if and only there is a $\mathfrak{z}\in\left(\mathbb{Z}_{2}\cap\mathbb{Q}\right)\backslash\left\{ 0\right\} $
so that $\chi_{q}\left(\mathfrak{z}\right)=x$.
\end{abstract}

\maketitle
\begin{quote}
\begin{flushright}
$\ldots$many people feel that functions $\mathbb{Z}_{p}\rightarrow\mathbb{Q}_{p}$
are more interesting than functions $\mathbb{Z}_{p}\rightarrow\mathbb{Q}_{q}$,
which is understandable. 
\par\end{flushright}
\begin{flushright}
\textemdash W. M. Schikhof\footnote{From \emph{Ultrametric Calculus} (\cite{Ultrametric Calculus}), on
the bottom of page 97.}
\par\end{flushright}
\end{quote}
\begin{note}
(September 2026) \textbf{Theorem \ref{thm:Divergent trajectories come from irrational z}}
from the original version of this paper has subsequently been demonstrated
to be false. The periodic-point Correspondence Principle remains valid.
An explicit counterexample and discussion of the correction are given
in \textbf{Section \ref{sec:Corrections}}.
\end{note}

\section*{Notation \& Preliminaries}

\subsubsection*{Notations for Sets}

For any prime $p$, we write $\mathbb{Z}_{p}$ and $\mathbb{Q}_{p}$,
to denote the ring of $p$-adic integers and fields of $p$-adic rational
numbers, respectively, equipped with the standard $p$-adic absolute
value $\left|\cdot\right|_{p}$. $\mathbb{Z}_{p}^{\times}$ is the
group of multiplicatively invertible $p$-adic integers; this is the
set of all $p$-adic integers which are not congruent to $0$ mod
$p$. We write $\mathbb{C}_{p}$ to denote the field of $p$-adic
complex numbers, i.e., the metric completion of the algebraic closure
of $\mathbb{Q}_{p}$. Note that $\mathbb{C}_{p}$ is \emph{not }spherically
complete.

We write $v_{p}\left(\cdot\right)$ to denote the $p$-adic valuation,
with $\left|\cdot\right|_{p}=p^{-v_{p}\left(\cdot\right)}$ and $v_{p}\left(0\right)\overset{\textrm{def}}{=}\infty$,
where $\overset{\textrm{def}}{=}$ means ``by definition''.

For any integer $N\geq1$, we write $\mathbb{Z}/N\mathbb{Z}$ to denote
the set $\left\{ 0,\ldots,N-1\right\} $. We write $\left(\mathbb{Z}/N\mathbb{Z}\right)^{\times}$
to denote the subset of $\left\{ 0,\ldots,N-1\right\} $ consisting
of all integers co-prime to $N$.

We write $\hat{\mathbb{Z}}_{p}$ to denote $\mathbb{Z}\left[1/p\right]/\mathbb{Z}$,
the Pontryagin dual of $\mathbb{Z}_{p}$. We identify $\hat{\mathbb{Z}}_{p}$
with the group of rational numbers in $\left[0,1\right)$ whose denominators
are powers of $p$. Viewing $\hat{\mathbb{Z}}_{p}$ as a subset of
$\mathbb{Q}_{p}$, for any $t\in\hat{\mathbb{Z}}_{p}\backslash\left\{ 0\right\} $,
observe that $t$ can be written in irreducible form as $k/p^{n}$
for some integers $n\geq1$ and $k\in\left(\mathbb{Z}/p^{n}\mathbb{Z}\right)^{\times}$.
Consequently, the $p$-adic absolute value of $t$ is $\left|t\right|_{p}=p^{n}$,
and $t\left|t\right|_{p}=k$, with both values being $0$ when $t=0$.

For any real number $x$, we write $\mathbb{N}_{x}$ to denote the
set of all integers $\geq x$ (thus, $\mathbb{N}_{0}=\left\{ 0,1,2,\ldots\right\} $;
$\mathbb{N}_{1}=\left\{ 1,2,\ldots\right\} $). We write $\mathbb{Z}_{p}^{\prime}$
to denote $\mathbb{Z}_{p}\backslash\mathbb{N}_{0}$; note that this
is the set of all $p$-adic integers with infinitely many non-zero
$p$-adic digits.

For aesthetic reasons, we write $p$-adic variables using lower-case
$\mathfrak{fraktur}$ font.

\subsubsection*{Notation for Congruences}

For any $\mathfrak{z}\in\mathbb{Z}_{p}$ and any integer $n\geq0$,
we write $\left[\mathfrak{z}\right]_{p^{n}}$ to denote the unique
integer in $\mathbb{Z}/p^{n}\mathbb{Z}$ which is congruent to $\mathfrak{z}$
modulo $p^{n}$. For terminological purposes, we follow Yvette Amice
in referring to the series representation $\sum_{n=-n_{0}}^{\infty}c_{n}p^{n}$
of a $p$-adic number $\mathfrak{y}\in\mathbb{Q}_{p}$ (be it rational
or integral) as $\mathfrak{y}$'s \textbf{Hensel series}.

Note that we adopt the convention of viewing $p$-adic integers as
power series in the variable $p$ (and, likewise, we view $p$-adic
rational numbers as Laurent series in the variable $p$), and, importantly
we write $p$-adic digits from \emph{left to right}, in order of increasing
powers of $p$. Thus, in $2$-adic notation, we would write:
\begin{align*}
1 & =\centerdot_{2}1\\
2 & =\centerdot_{2}01\\
3 & =\centerdot_{2}11\\
4 & =\centerdot_{2}001\\
 & \vdots
\end{align*}
and so on. Here, we use $\centerdot_{2}$ as the ``$2$-adic decimal
point''. Throughout this series, all $2$-adic numbers written in
base $2$ will have a $\centerdot_{2}$ in them. For example:
\begin{equation}
11\centerdot_{2}011=\frac{1}{2^{2}}+\frac{1}{2}+2^{1}+2^{2}=6+\frac{3}{4}
\end{equation}

We employ a non-standard notation for congruences:
\begin{align}
x & \overset{a}{\equiv}y\nonumber \\
 & \Updownarrow\\
x & =y\mod a\nonumber 
\end{align}
To give some examples, for $\mathfrak{z}\in\mathbb{Z}_{p}$, $\mathfrak{z}\overset{p^{n}}{\equiv}k$
means ``$\mathfrak{z}$ is congruent to $k$ mod $p^{n}$''; i.e.,
$\mathfrak{z}\in k+p^{n}\mathbb{Z}_{p}$. Given $\mathfrak{x},\mathfrak{y}\in\mathbb{Q}_{p}$,
we write $\mathfrak{x}\overset{1}{\equiv}\mathfrak{y}$ to denote
that $\mathfrak{x}-\mathfrak{y}=0\mod\mathbb{Z}_{p}$. In particular,
note that $\mathfrak{y}\overset{1}{\equiv}\left\{ \mathfrak{y}\right\} _{p}$
for all $\mathfrak{y}\in\mathbb{Q}_{p}$, that $s,t\in\hat{\mathbb{Z}}_{p}$
denote the same element of $\hat{\mathbb{Z}}_{p}$ if and only if
$s\overset{1}{\equiv}t$, and that $\mathfrak{z}\overset{1}{\equiv}k$
is true for all $\mathfrak{z}\in\mathbb{Z}_{p}$ and all $k\in\mathbb{Z}$.
Also, for a rational integer $n$, $n\overset{2}{\equiv}0$ means
$n$ is even, while $n\overset{2}{\equiv}1$ means $n$ is odd. Naturally,
we write $x\overset{a}{\cancel{\equiv}}y$ to mean ``$x$ is \emph{not
}congruent to $y$ mod $a$''.

\subsubsection*{Notation for Convergence}

Given an equation with a limit or sum, we write $\overset{\mathbb{Q}_{p}}{=}$
to indicate that the convergence is with respect to the topology of
$\mathbb{Q}_{p}$. Similarly, $\overset{\mathbb{R}}{=}$ means convergence
with respect to the topology of $\mathbb{R}$. We write $\overset{\mathbb{Q}}{=}$
and $\overset{\overline{\mathbb{Q}}}{=}$ to denote convergence with
respect to the discrete topologies on $\mathbb{Q}$ and $\overline{\mathbb{Q}}$,
respectively; note then that a sequence of rational numbers $\left\{ a_{n}\right\} _{n\geq0}$
satisfies $\lim_{n\rightarrow\infty}a_{n}\overset{\mathbb{Q}}{=}0$
if and only if $a_{n}=0$ for all sufficiently large $n$.

\subsubsection*{Preliminaries from Dynamical Systems}

Given a map $T$, we write $T^{\circ k}$ to denote the composition
of $k$ copies of $T$. We also define $T^{\circ0}$ to be the identity
map.

For the remaining definitions, $T$ is a map $T:\mathbb{Z}\rightarrow\mathbb{Z}$.
Recall that a \textbf{periodic point }of $T$ is an $x\in\mathbb{Z}$
so that $T^{\circ k}\left(x\right)=x$ for some $k\geq1$. We say
$x$ is \textbf{pre-periodic }if there exists a $k\geq1$ so that
$T^{\circ k}\left(x\right)$ is a periodic point of $T$. We say $x$
is a \textbf{strictly pre-periodic }point if $x$ is a pre-periodic
point of $T$ which is \emph{not }a periodic point of $T$. The \textbf{forward
orbit }of $x\in\mathbb{Z}$ under $T$ is the sequence:
\begin{equation}
\left\{ T^{\circ k}\left(x\right)\right\} _{k\geq0}
\end{equation}
A \textbf{cycle }$\Omega$ of $T$ is the forward orbit of a periodic
point of $T$. On the other hand, a \textbf{divergent point }of $T$,
meanwhile, is an $x\in\mathbb{Z}$ so that: 
\begin{equation}
\lim_{k\rightarrow\infty}\left|T^{\circ k}\left(x\right)\right|\overset{\mathbb{R}}{=}\infty
\end{equation}
A \textbf{divergent trajectory }of $T$ is the forward orbit of a
divergent point.

An \textbf{orbit class }of $T$ is a non-empty set $V\subseteq\mathbb{Z}$
so that the pre-image of $V$ under $T$ is equal to $V$: $T^{-1}\left(V\right)=V$.
An orbit class $V$ of $T$ is said to \textbf{irreducible }if it
cannot be written as $V=V_{1}\cup V_{2}$, where $V_{1}$ and $V_{2}$
are orbit classes of $T$ with $V_{1}\cap V_{2}=\varnothing$. A \textbf{backwards
orbit }of some $x_{0}\in\mathbb{Z}$ is a sequence $\left\{ x_{n}\right\} _{n\geq1}$
in $\mathbb{Z}$ so that $x_{n+1}\in T^{-1}\left(\left\{ x_{n}\right\} \right)$
for all $n\geq0$. Note that if $T$ is not injective, a given $x_{0}\in\mathbb{Z}$
can possess more than one backward orbit, potentially infinitely many.

We say an orbit class $V$ is \textbf{periodic }if it is irreducible
and contains a cycle. We say an orbit class $V$ is \textbf{divergent
}if it is irreducible and contains a divergent trajectory.

We also recall the following two elementary results from the theory
of dynamical systems:
\begin{prop}
Let $T:\mathbb{Z}\rightarrow\mathbb{Z}$ be any map. Then, the collection
of irreducible orbit classes of $T$ form a partition of $\mathbb{Z}$
into at most countably infinitely many pair-wise disjoint sets. These
sets are the equivalence classes of $\mathbb{Z}$ under the relation
$\sim$ defined by:
\begin{equation}
x\sim y\Leftrightarrow\exists m,n\geq0:T^{\circ m}\left(x\right)=T^{\circ n}\left(y\right)
\end{equation}
\end{prop}
\begin{prop}
\label{prop:Classification of T_3's dynamics}Let $T:\mathbb{Z}\rightarrow\mathbb{Z}$
be any map. Every irreducible orbit class of $T$ is either periodic
or divergent. Consequently, every $x\in\mathbb{Z}$ is either a divergent
point of $T$ or\emph{ }a pre-periodic\textbf{ }point of $T$.
\end{prop}

\subsubsection*{Preliminaries from $p$-Adic Analysis}

For a refresher on the topic, the reader can refer to Robert's book
\cite{Robert's Book}, Schikhof's book \cite{Ultrametric Calculus},
or the author's doctoral dissertation \cite{my dissertation}. For
the purpose of this paper, we need only one result from the theory
of ultrametric analysis:
\begin{prop}
Let $p$ be a prime, and let $\left\{ \mathfrak{a}_{n}\right\} _{n\geq0}$
be a sequence in $\mathbb{Q}_{p}$. Then, the series $\sum_{n=0}^{\infty}\mathfrak{a}_{n}$
converges in $\mathbb{Q}_{p}$ if and only if $\lim_{n\rightarrow\infty}\left|\mathfrak{a}_{n}\right|_{p}\overset{\mathbb{R}}{=}0$.

\newpage{}
\end{prop}

\section*{\label{sec:Introduction}Introduction}
\begin{quote}
\begin{flushright}
\emph{This book of mine has little need of preface, for indeed it
is 'all preface' from beginning to end.}
\par\end{flushright}
\end{quote}
\begin{quotation}
\textemdash D'Arcy Thompson, \emph{On Growth and Form} \cite{D'Arcy}
\end{quotation}
\vphantom{}

We begin with a trouble-maker.
\begin{defn}
The\textbf{ Collatz map} $C:\mathbb{Z}\rightarrow\mathbb{Z}$ is the
function defined by:
\begin{equation}
C\left(n\right)\overset{\textrm{def}}{=}\begin{cases}
\frac{n}{2} & \textrm{if }n\overset{2}{\equiv}0\\
3n+1 & \textrm{if }n\overset{2}{\equiv}1
\end{cases}\label{eq:Collatz Map}
\end{equation}
\end{defn}
The infamous \textbf{Collatz Conjecture} is the assertion:
\begin{conjecture}[\textbf{Collatz Conjecture, Ver. 1 }\cite{Lagarias' Survey}]
For every integer $n\geq1$, there exists a $k\geq0$ so that $C^{\circ k}\left(n\right)=1$.
\end{conjecture}
Rather than work with $C\left(n\right)$ directly, we shall consider
a well-known acceleration of $C\left(n\right)$.
\begin{defn}
The \textbf{Shortened Collatz map }is the function $T_{3}:\mathbb{Z}\rightarrow\mathbb{Z}$
defined by:
\begin{equation}
T_{3}\left(n\right)\overset{\textrm{def}}{=}\begin{cases}
\frac{n}{2} & \textrm{if }n\overset{2}{\equiv}0\\
\frac{3n+1}{2} & \textrm{if }n\overset{2}{\equiv}1
\end{cases}\label{eq:Shortened Collatz Map-1}
\end{equation}
which is ``shortened'' in the sense that:

\begin{equation}
T_{3}\left(n\right)=\begin{cases}
C\left(n\right) & \textrm{if }n\overset{2}{\equiv}0\\
C\left(C\left(n\right)\right) & \textrm{if }n\overset{2}{\equiv}1
\end{cases}
\end{equation}
\end{defn}
With this, the Collatz Conjecture can be reformulated as:
\begin{conjecture}[\textbf{Collatz Conjecture, Ver. 2}]
For every integer $n\geq1$, there exists a $k\geq0$ so that $T_{3}^{\circ k}\left(n\right)=1$.
\end{conjecture}
A moment's thought shows that this conjecture actually consists of
two sub-conjectures. First, there is the so-called \textbf{Weak Collatz
Conjecture}:
\begin{conjecture}[\textbf{Weak Collatz Conjecture \cite{Tao Blog}}]
\label{conj:Weak Collatz Conjecture}The only periodic points of
$T_{3}$ in the positive integers are $1$ and $2$.
\end{conjecture}
Second is the aptly named \textbf{No Divergent Points Conjecture}:
\begin{conjecture}[\textbf{No Divergent Points Conjecture \cite{Tao Blog}}]
\label{conj:No divergent points}No positive integer is a divergent
point of $T_{3}$.
\end{conjecture}
Technically speaking, these conjectures (and others like them) fall
into the category of arithmetic dynamics. With its beginnings in Joe
Silverman's classic text \cite{Silverman}, the subject of arithmetic
studies rational maps on algebraic varieties from the perspective
of function iteration dynamics. This dovetails classical concepts
from the theory of discrete dynamical systems with the subjects and
vocabulary of algebraic number theory and arithmetic geometry, forming
a powerful, fertile approach which synthesizes insights from multiple
mathematical disciplines. Unfortunately, this kind of arithmetic dynamics
has very little to say about the Collatz Conjecture; $T_{3}$ is simply
too far removed from the foundational algebraic structures studied
in arithmetic geometry (algebraic varieties, algebraic groups, abelian
varieties, schemes, etc.) for these methods to apply.

In this regard, Collatz poses a much \emph{meta}mathematical challenge
as much as a mathematical one. It is noteworthy that, unlike many
other infamous, difficult number-theoretic problems (\textbf{s}uch
as \textbf{Fermat's Last Theorem}, the \textbf{Riemann Hypothesis},
or the \textbf{Birch and Swinnerton-Dyer Conjecture}) the Collatz
Conjecture does not seem to readily connect to other areas of mathematics\textemdash the
sole exception being Collatz's well-established connection to transcendental
number theory, itself a notoriously difficult field \cite{Tao Blog}.
Collatz stands alone, inscrutable and unapproachable, seemingly lacking
a broader context in which it might be better understood\textemdash until
now. The work done in the author's doctoral dissertation (\cite{my dissertation})
shows that this unapproachability is longer the case, and to that
end, this series of papers is intended to present a refined distillation
of those findings for a broader mathematical audience.

To be clear, \emph{we }\emph{\uline{do not}}\emph{ have a proof
of Collatz, nor even a partial one}. Rather, what we offer what is
an undiscovered country, one rich with unexpected structure. As in
\cite{my dissertation}, this series of papers focus on three main
innovations. These are:
\begin{enumerate}
\item A novel \emph{unifying }analytic formalism for studying Collatz-type
maps: a function, the \textbf{numen}, which we associate to a given
Collatz-type map. With this formalism, the author's original research
program attempted to show that the questions, ``Is $x$ a periodic
point of the map?'' and, ``Is $x$ a divergent point of the map''
could be restated in terms of the value-distribution of the numen;
that is, ``Is $x$ a value attained by the numen?'' However, the
example constructed in Section \ref{sec:Corrections} renders this
program purely conjectural at the moment. Nevertheless, as will be
discussed in future work, there is a potential avenue to recover the
original research program involving a quasi-adèlic treatment of $\chi_{q}$. 
\item The discovery of an unexpected level of depth and subtlety in an almost
totally neglected area of non-archimedean analysis:\textbf{ $\left(p,q\right)$-Adic
Analysis}, the study of functions from the $p$-adics to the $q$-adics,
where $p$ and $q$ are distinct primes.
\item The use of (1) and (2) to reformulate Collatz-type conjectures into
problems of \textbf{Non-Archimedean Spectral Theory}\footnote{Let $K$ be a metrically complete non-archimedean valued field, let
$\mathcal{X}$ be a unital Banach algebra over $K$ (with $1_{\mathcal{X}}$
as its multiplicative identity element), and let $\chi\in\mathcal{X}$.
By \textbf{Non-Archimedean Spectral Theory}, we mean the study of
those scalars $\lambda\in K$ for which $\chi-\lambda1_{\mathcal{X}}$
is a multiplicative invertible element of $\mathcal{X}$ (a \textbf{unit}
of $\mathcal{X}$).}. Combined, (2) and (3) show that we can reformulate Collatz-type
conjectures using a totally unexplored set of tools.
\end{enumerate}
The reformulation mentioned in (3) has three main ingredients. The
first is the novel analytic formalism described in (1): a function
$\chi_{3}:\mathbb{Z}_{2}\rightarrow\mathbb{Z}_{3}$; the \textbf{numen}\footnote{The name is from the Latin, meaning ``the spirit or power presiding
over a thing or place''. Originally, the author called $\chi_{3}$
the ``characteristic function'' of $T_{3}$, however, calling it
this would cause a thorny conflict with the use of a bonafide ``characteristic
function'' (in the sense of ``characteristic function of a random
variable'') used by Tao in his 2019 paper on the Collatz Conjecture
\cite{Tao Probability paper}. Although the thrust of our work here
and its progenitors in the author's PhD dissertation are all but orthogonal
to the thoroughly probabilistic approach used by Tao to analyze the
Collatz map by way of \textbf{Syracuse Random Variables}, both Tao's
approach and our own rest squarely on $\chi_{3}$ and its properties.
However, this commonality might not be apparent at first glance due
to the wildly different perspective (and notation) used by Tao in
his paper. The interested reader can consult \textbf{Section \ref{subsec:The-Tao-of}}
below for a comprehensive account of the relationship between our
work and Tao's.}\textbf{ }(plural: numina) of $T_{3}$ \cite{my dissertation}. $\chi_{3}$
can be characterized as\textemdash and, if so desired, \emph{defined}
\emph{as}\textemdash as the unique function $f:\mathbb{Z}_{2}\rightarrow\mathbb{Z}_{3}$
satisfying the functional equations:
\begin{align}
f\left(2\mathfrak{z}\right) & =\frac{f\left(\mathfrak{z}\right)}{2}\\
f\left(2\mathfrak{z}+1\right) & =\frac{3f\left(\mathfrak{z}\right)+1}{2}
\end{align}
for all $\mathfrak{z}\in\mathbb{Z}_{2}$, subject to the following
limit condition:
\begin{equation}
\lim_{n\rightarrow\infty}f\left(\mathfrak{z}\right)\overset{\mathbb{Q}_{3}}{=}f\left(\left[\mathfrak{z}\right]_{2^{n}}\right),\textrm{ }\forall\mathfrak{z}\in\mathbb{Z}_{2}\label{eq:Rising continuity}
\end{equation}
(See \textbf{Lemma} \textbf{\ref{lem:Chi_H functional equation on N_0 and uniqueness}}
on page \pageref{lem:Chi_H functional equation on N_0 and uniqueness}
for the proof.) In keeping with the terminology of the author's dissertation,
limit conditions in the manner of (\ref{eq:Rising continuity}) will
be called \textbf{rising-continuity}. Specifically:
\begin{defn}[Rising Continuity \cite{my dissertation}]
Let $p$ be a prime, and let $K$ be a metrically complete non-archimedean
valued field. We say a function $f:\mathbb{Z}_{p}\rightarrow K$ is
\textbf{rising-continuous }if:
\begin{equation}
\lim_{n\rightarrow\infty}f\left(\mathfrak{z}\right)\overset{K}{=}f\left(\left[\mathfrak{z}\right]_{p^{n}}\right),\textrm{ }\forall\mathfrak{z}\in\mathbb{Z}_{p}\label{eq:definition of rising continuity}
\end{equation}
where, as indicated, the convergence of the limit occurs in the topology
of $K$. Note that we only require the convergence to occur point-wise.
\end{defn}
\begin{rem}
It can be shown that the set of rising-continuous functions $f:\mathbb{Z}_{p}\rightarrow K$
forms a non-archimedean Banach algebra under the usual operations
of point-wise addition, point-wise multiplication, and scalar multiplication.
This algebra contains $C\left(\mathbb{Z}_{p},K\right)$\textemdash the
Banach algebra of continuous functions $f:\mathbb{Z}_{p}\rightarrow K$\textemdash as
a dense subalgebra. Importantly, every rising-continuous function
is uniquely represented by a \textbf{van der Put series} \cite{my dissertation}.
\end{rem}
We use the term \textbf{$\left(p,q\right)$-adic function }to refer
to a function from the $p$-adic numbers to the $q$-adic numbers,
and let $p$ (resp. $q$) be $\infty$ when the domain (resp. co-domain)
of our functions is a Euclidean space over the real or complex numbers
($\mathbb{R}$, $\mathbb{C}$, $\mathbb{R}^{2}$, etc.); in this terminology,
we would describe $\chi_{3}$ as a $\left(2,3\right)$-adic function.
The study of such functions dates back to the work of W.M. Schikhof,
A.C.M van Rooij, and other mathematicians of the Dutch school of Non-Archimedean
(Functional) Analysis in the mid-1960s \cite{Schikhof's Thesis,Ultrametric Calculus,van Rooij - Non-Archmedean Functional Analysis}.
For better and for worse, the pioneers of non-archimedean functional
analysis believed $\left(p,q\right)$-adic analysis to be ``uninteresting''.
As such, they left it by the wayside in favor of more fruitful (and,
generally, more abstract) areas of non-archimedean analysis.

As we shall see in later installments of this series of papers, though
this judgment was entirely understandable, it was, nevertheless, premature.
At the risk of spoiling the surprise, the current state of affairs
is merely a failure of imagination. Outside of hard analysis\textemdash and,
especially, in subjects like algebraic number theory and algebraic
geometry, where the $p$-adics find their most consequential applications\textemdash mathematicians
tend to be drawn to \emph{continuous }functions. This is a perfectly
reasonable preference to have, and, as long as one insists on working
with continuous functions, Schikhof's quip about the lesser interest
of functions \emph{$\mathbb{Z}_{p}\rightarrow\mathbb{Q}_{q}$} is
absolutely correct; $C\left(\mathbb{Z}_{p},\mathbb{Q}_{q}\right)$
is rigid and almost disappointingly simplistic. But this does not
mean that $\left(p,q\right)$-adic analysis has no interesting results.
Rather, it merely tells us we cannot expect to get interesting results
in $\left(p,q\right)$-adic analysis as long as we continue to mandate
that our functions be continuous.

As it turns out, we can get much richer\textemdash and in several
ways, entirely \emph{novel}\textemdash results if we allow ourselves
to work with functions such as the rising-continuous functions: not
fully continuous, yes, but nearly so. As far as we can tell, the rising-continuous
functions seem to be a more ``natural'' choice of functions to work
with in the $\left(p,q\right)$-adic context than their more rigid,
fully continuous cousins. In this aspect, our work in this series
of papers is as much a development of $\left(p,q\right)$-adic analysis
as it is a development of Collatz studies. The two seem to be inextricably
linked; by examining one, we find new questions to ask of the other.

This, on its own, should be enough reason to care about both $\chi_{3}$
and $\left(p,q\right)$-adic analysis as a whole, but, in fact, there
is more, and this stems from the second ingredient of our reformulation:
the result we call the \textbf{Correspondence Principle} (\textbf{CP}).
The CP comes in two flavors, one for periodic points, and another
for divergent points:
\begin{quotation}
\emph{I. (}\textbf{\emph{Corollary \ref{cor:CP v4}}}\emph{, page
\pageref{cor:CP v4}) $x\in\mathbb{Z}\backslash\left\{ 0\right\} $
is a periodic point of $T_{3}$ }if and only if\emph{ there exists
a $\mathfrak{z}\in\mathbb{Z}_{2}^{\prime}\cap\mathbb{Q}$ so that
$\chi_{3}\left(\mathfrak{z}\right)=x$. (}\textbf{\emph{CPPP}}\emph{)}

\emph{\vphantom{}}

\emph{II. (}\textbf{\emph{Theorem \ref{thm:Divergent trajectories come from irrational z}}}\emph{,
page \pageref{thm:Divergent trajectories come from irrational z})
If there is a $\mathfrak{z}\in\mathbb{Z}_{2}\backslash\mathbb{Q}$
so that $\chi_{3}\left(\mathfrak{z}\right)\in\mathbb{Z}$, then $\chi_{3}\left(\mathfrak{z}\right)$
is a divergent point of $T_{3}$. (}\textbf{\emph{CPDP}}\emph{)}
\end{quotation}
\begin{rem}
Unlike (I), the proof of (II) is \emph{completely non-constructive}.
It is a near-immediate consequence of (I) and \textbf{Proposition}
\textbf{\ref{prop:Classification of T_3's dynamics}}. We strongly
suspect that the converse of (II) is true, though we have no proof
of this.
\end{rem}
\begin{rem}
Alas, \textbf{Theorem \ref{thm:Divergent trajectories come from irrational z}}\textbf{\emph{
}}is false. See \textbf{Section \ref{sec:Corrections} }for an explicit
counterexample and further discussion.
\end{rem}
The third and final ingredient of our spectral-theoretic reformulation
of Collatz emerges from a new $\left(p,q\right)$-adic generalization\footnote{In light of \cite{my paper}, this could be phrased as a \emph{frame-theoretic}
generalization.} of the classic \textbf{Wiener Tauberian Theorem }(WTT) of harmonic
analysis, first established in \cite{my dissertation}. A proof of
this result will be given in a later paper in this series, along with
an in-depth presentation of the pertinent methods of Fourier analysis.
For now, it suffices to say that the methods outlined in the series
of papers touch upon a variety of unusual and completely novel analytic
phenomena; in the author's estimation, these are in dire need of a
firm theoretical foundation.

From a birds'-eye-view, the purpose of this series of papers is to
give a comprehensive exposition of a special case of the broader programme
of study the author initiated in his PhD dissertation. These deals
with a general family of Collatz-type maps defined on $\mathbb{Z}^{d}$,
where $d$ is an integer $\geq1$ \cite{my dissertation}. We call
these \textbf{Hydra maps}\footnote{The name is due to the Hydra-like appearances of some visualizations
of orbit classes of the Collatz map, as well as the difficulty defeating
the Hydra of Greek myth. }, and they are defined as follows:
\begin{defn}[\textbf{A One-Dimensional Hydra Map}]
Let $p$ be a prime number\footnote{The theory presented here can likely be made to work when $p$ is
composite, though it will entails certain technical complications.} $\geq2$. A \textbf{(one-dimensional) $p$-Hydra map} is a map $H:\mathbb{Z}\rightarrow\mathbb{Z}$
of the form:
\begin{equation}
H\left(n\right)=\begin{cases}
\frac{a_{0}n+b_{0}}{d_{0}} & \textrm{if }n\overset{p}{\equiv}0\\
\frac{a_{1}n+b_{1}}{d_{1}} & \textrm{if }n\overset{p}{\equiv}1\\
\vdots & \vdots\\
\frac{a_{p-1}n+b_{p-1}}{d_{p-1}} & \textrm{if }n\overset{p}{\equiv}p-1
\end{cases},\textrm{ }\forall n\in\mathbb{Z}\label{eq:Def of a Hydra Map on Z}
\end{equation}
where $a_{j}$, $b_{j}$, and $d_{j}$ are integer constants (with
$a_{j},d_{j}\geq0$ for all $j$) so that the following two properties
hold:

I. (``co-primality'') $a_{j},d_{j}>0$ and $\gcd\left(a_{j},d_{j}\right)=1$
for all $j\in\left\{ 0,\ldots,p-1\right\} $.

\vphantom{}

II. (``integrality'') For each $j\in\left\{ 0,\ldots,p-1\right\} $,
$\frac{a_{j}n+b_{j}}{d_{j}}$ is an integer if and only if $n\overset{p}{\equiv}j$.

Also, for each $j\in\mathbb{Z}/p\mathbb{Z}$, we write $H_{j}:\mathbb{Q}\rightarrow\mathbb{Q}$
to denote the map: 
\begin{equation}
H_{j}\left(x\right)\overset{\textrm{def}}{=}\frac{a_{j}x+b_{j}}{d_{j}}
\end{equation}
and call this the \textbf{$j$th branch }of $H$.
\end{defn}
A guiding principle for \cite{my dissertation} was the author's conviction
that the Collatz map ought to be studied \emph{in relation} \emph{to
other maps of its type}. Doing so is the only way we can hope to determine
whether a given observation of the Collatz map and its dynamics is
truly related to map's conjectured behavior, or merely a mathematical
coincidence. Moreover, it is highly desirable to have a unified formalism
for working with maps of thiis type. The techniques and terminologies
used in Collatz-studies are often idiosyncratic and highly variable;
having a single ``language'' that unites and subsumes all these
approaches would help bring order to that chaos. Even more significantly,
a unifying formalism would be more amenable to generalization.

Case in point, the method of the \textbf{numen }admits an immediate
and straight-forward generalization to Collatz-type maps defined on
rings of algebraic integers, such as $\mathbb{Z}\left[\sqrt{3}\right]$.
(The first example of such a map is due to Leigh; an account of it
and other maps of that type can be found in K.R. Matthews' excellent
power-point slides \cite{Matthews' slides}.) Our formalism treats
these as ``multi-dimensional'' Hydra maps, which are defined as
follows. (It should be noted that, for technical reasons, we will
state the definition in terms of ideals of a ring of algebraic integers,
but in practice, one works with affine linear maps on $\mathbb{Z}^{d}$
for an integer $d\geq2$; the details are given below.)
\begin{defn}[\textbf{A Multi-Dimensional Hydra Map}]
Let $\mathbb{F}$ be a number field with $\left[\mathbb{F}:\mathbb{Q}\right]=d$.
When $d=1$, we obtain a Hydra map on $\mathbb{Z}\rightarrow\mathbb{Z}$
by defining what $H$ does to a given $n$ based on the value of $n$
modulo some prime $p$. For $\mathbb{F}$, the analogue of this is
a map on $\mathcal{O}_{\mathbb{F}}\rightarrow\mathcal{O}_{\mathbb{F}}$
(the ring of $\mathbb{F}$-integers) whose action on a given $z\in\mathcal{O}_{\mathbb{F}}$
is determined by the value of $z$ modulo $\mathfrak{I}$, where $\mathfrak{I}$
is some proper ideal of $\mathcal{O}_{\mathbb{F}}$ (the $d=1$ case
has $\mathfrak{I}$ as the ideal of $\mathbb{Z}$ generated by the
prime $p$).

Fixing a non-zero proper ideal $\mathfrak{I}$ of $\mathcal{O}_{\mathbb{F}}$
of index $\iota$ (the index being defined as $\left|\mathcal{O}_{\mathbb{F}}/\mathfrak{I}\right|$,
the cardinality of the quotient group $\mathcal{O}_{\mathbb{F}}/\mathfrak{I}$),
the \textbf{Structure Theorem for Finitely-Generated Modules over
a Principal Ideal Domain} tells us that there is an isomorphism of
additive groups:
\begin{equation}
\mathcal{O}_{\mathbb{F}}/\mathfrak{I}\cong\prod_{n=1}^{r}\mathbb{Z}/p_{n}\mathbb{Z}\label{eq:Direct Product Representation of O_F / I}
\end{equation}
where $r\in\left\{ 1,\ldots,\left|\mathcal{O}_{\mathbb{F}}/\mathfrak{I}\right|\right\} $,
and where the $p_{n}$s are integers $\geq2$ so that $p_{n}\mid p_{n+1}$
for all $n\in\left\{ 1,\ldots,r-1\right\} $. We call $r$ the \textbf{depth
}of $\mathfrak{I}$; note that $r\leq d$.

Because our goal is to use the structure of $\mathbb{F}$ as a $d$-dimensional
vector space over $\mathbb{Q}$ to represent our multi-dimensional
Hydra map using matrices, we will need to choose a basis for this
vector space. By elementary linear algebra, we know that there is
a set $\mathcal{B}\subset\mathcal{O}_{\mathbb{F}}$ which forms a
$\mathbb{Z}$-basis of the group $\left(\mathcal{O}_{\mathbb{F}},+\right)$
(and hence, forms a basis of $\mathbb{F}$ as a vector space over
$\mathbb{Q}$), so that, for any element element $z\in\mathcal{O}_{\mathbb{F}}$,
the image of $z$ under the projection from $\mathcal{O}_{\mathbb{F}}/\mathfrak{I}$
to $\mathbb{Z}/p_{n}\mathbb{Z}$ is the value of $x_{n}$ mod $p_{n}$,
where $\mathbf{x}=\left(x_{1},\ldots,x_{d}\right)\in\mathbb{Z}^{d}$
is the $d$-tuple representing $\mathbf{z}$ in $\mathcal{B}$-coordinates.
\textbf{We choose once and for all such a basis $\mathcal{B}$}.

Using $\mathcal{B}$, we identify $\mathcal{O}_{\mathbb{F}}$ with
the lattice $\mathbb{Z}^{d}$. For a Hydra map $H$ on $\mathbb{Z}$,
we distinguish between the restrictions of $H$ to $\mathbb{N}_{0}$
and $-\mathbb{N}_{0}$; in other words, we work on the half-lattices
of non-negative and non-positive integers, respectively. We realize
this feature of Hydra maps on $\mathcal{O}_{\mathbb{F}}$ by writing
$\mathcal{O}_{\mathbb{F},\mathcal{B}}^{+}$ (resp. $\mathcal{O}_{\mathbb{F},\mathcal{B}}^{-}$)
to denote the set of all elements of $\mathcal{O}_{\mathbb{F}}$ whose
$\mathcal{B}$-coordinate vectors have non-negative (resp. non-positive)
integer entries. Finally, we write $\mathfrak{I}_{0},\ldots,\mathfrak{I}_{\iota}$
to denote the cosets of $\mathfrak{I}$ in $\mathcal{O}_{\mathbb{F}}$,
with $\mathfrak{I}_{0}$ denoting $\mathfrak{I}$ itself.

With all these conventions in place, an \textbf{$\left(\mathbb{F},\mathfrak{I},\mathcal{B}\right)$-Hydra
map} \textbf{on $\mathcal{O}_{\mathbb{F}}$} (more succinctly, an
$\mathfrak{I}$-Hydra map) is a surjective map $\tilde{H}:\mathcal{O}_{\mathbb{F}}\rightarrow\mathcal{O}_{\mathbb{F}}$
of the form:
\begin{equation}
\tilde{H}\left(z\right)=\begin{cases}
\frac{a_{0}z+b_{0}}{d_{0}} & \textrm{if }z\in\mathfrak{I}_{0}\\
\vdots & \vdots\\
\frac{a_{\iota-1}z+b_{\iota-1}}{d_{\iota-1}} & \textrm{if }z\in\mathfrak{I}_{\iota}
\end{cases}\label{eq:Definition of I-hydra map-1}
\end{equation}
Here, the $a_{j}$s, $b_{j}$s, and $d_{j}$s are elements of $\mathcal{O}_{\mathbb{F}}$
so that:

I. $a_{j},d_{j}\neq0$ for all $j\in\left\{ 0,\ldots,\iota-1\right\} $.

\vphantom{}

II. $\gcd\left(a_{j},d_{j}\right)=1$ for all $j\in\left\{ 0,\ldots,\iota-1\right\} $.

\vphantom{}

III. For all $j\in\left\{ 0,\ldots,\iota-1\right\} $, the quantity
$\frac{a_{j}z+b_{j}}{d_{j}}$ is an element of $\mathcal{O}_{\mathbb{F},\mathcal{B}}^{+}$
if and only if $z\in\mathcal{O}_{\mathbb{F},\mathcal{B}}^{+}\cap\mathfrak{I}_{j}$.

\vphantom{}

IV. For all $j\in\left\{ 0,\ldots,\iota-1\right\} $, the quantity
$\frac{a_{j}z+b_{j}}{d_{j}}$ is an element of $\mathcal{O}_{\mathbb{F},\mathcal{B}}^{-}$
if and only if $z\in\mathcal{O}_{\mathbb{F},\mathcal{B}}^{-}\cap\mathfrak{I}_{j}$.

\vphantom{}

V. For all $j\in\left\{ 0,\ldots,\iota-1\right\} $, the ideal $\left\langle d_{j}\right\rangle _{\mathcal{O}_{\mathbb{F}}}$
in $\mathcal{O}_{\mathbb{F}}$ generated by $d_{j}$ is contained
in $\mathfrak{I}$.

\vphantom{}

VI. For all $j\in\left\{ 0,\ldots,\iota-1\right\} $, the the matrix
representation in $\mathcal{B}$-coordinates on $\mathbb{Q}^{d}$
of the ``multiplication by $\frac{a_{j}}{d_{j}}$'' map on $\mathbb{F}$
is of the form:
\begin{equation}
\frac{\mathbf{A}}{\mathbf{D}}\overset{\textrm{def}}{=}\mathbf{D}^{-1}\mathbf{A}\label{eq:A / D notation-1}
\end{equation}
where $\mathbf{A},\mathbf{D}$ are invertible $d\times d$ matrices
so that:

\vphantom{}

VI-i. $\mathbf{A}=\tilde{\mathbf{A}}\mathbf{P}$, where $\mathbf{P}$
is a permutation matrix\footnote{That is, a matrix of $0$s and $1$s which is a representation of
the action on $\mathbb{Z}^{d}$ of an element of the symmetric group
on $d$ objects by way of a permutation of the coordinate entries
of the $d$-tuples in $\mathbb{Z}^{d}$.} and where $\mathbf{\tilde{\mathbf{A}}}$ is a diagonal matrix whose
non-zero entries are positive integers.

\vphantom{}

VI-ii. $\mathbf{D}$ is a diagonal matrix whose non-zero entries are
positive integers such that every entry on the diagonal of $\mathbf{R}\mathbf{D}^{-1}$
is a positive integer. Note that this then forces the $\left(r+1\right)$th
through $d$th diagonal entries of $\mathbf{D}$ to be equal to $1$.

Lastly, we define the \textbf{depth }of $\tilde{H}$ as the depth
of the ideal $\mathfrak{I}$ (i.e., $r$).

While the above definition is conceptually straight-forward, it is
far too messy to use in compuations, so, in practice, all work with
these maps is done over $\mathbb{Z}^{d}$. Having defined $\tilde{H}$,
we can forget about $\mathbb{F}$, $\mathcal{O}_{\mathbb{F}}$, $\mathfrak{I}$,
and $\mathcal{B}$ and pass to isomorphic copies thereof. Instead
of $\mathcal{O}_{\mathbb{F}}$, we work in the lattice $\mathbb{Z}^{d}$.
It then follows that $\tilde{H}$ can be realized as a map $H:\mathbb{Z}^{d}\rightarrow\mathbb{Z}^{d}$
with branches of the form:
\begin{equation}
\mathbf{x}\in\mathbb{Z}^{d}\mapsto\mathbf{D}_{\mathbf{j}}^{-1}\left(\mathbf{A}_{\mathbf{j}}\mathbf{x}+\mathbf{b}_{\mathbf{j}}\right)\in\mathbb{Z}^{d}
\end{equation}
for invertible matrices $\mathbf{D}_{\mathbf{j}}$ and $\mathbf{A}_{\mathbf{j}}$
and $d\times1$ column vectors $\mathbf{b}_{\mathbf{j}}$. Here, the
$\mathbf{j}$s are no longer strings, but $r$-tuples $\mathbf{j}=\left(j_{1},\ldots,j_{r}\right)$
in an additive group of the form:
\begin{equation}
\mathbb{Z}^{r}/P\mathbb{Z}^{r}\overset{\textrm{def}}{=}\prod_{m=1}^{r}\left(\mathbb{Z}/p_{m}\mathbb{Z}\right)
\end{equation}
where $P=\left(p_{1},\ldots,p_{r}\right)$. That being done, we then
define a \textbf{$d$-dimensional $P$-Hydra map }$H:\mathbb{Z}^{d}\rightarrow\mathbb{Z}^{d}$
to be a map of the form:
\begin{equation}
H\left(\mathbf{x}\right)\overset{\textrm{def}}{=}\sum_{\mathbf{j}\in\mathbb{Z}^{r}/P\mathbb{Z}^{r}}\mathbf{D}_{\mathbf{j}}^{-1}\left(\mathbf{A}_{\mathbf{j}}\mathbf{x}+\mathbf{b}_{\mathbf{j}}\right)\left[\mathbf{x}\overset{P}{\equiv}\mathbf{j}\right],\textrm{ }\forall\mathbf{x}\in\mathbb{Z}^{d}\label{eq:Form a d-dimensional P-Hydra map}
\end{equation}
where $\mathbf{x}\overset{P}{\equiv}\mathbf{j}$ is our shorthand
for the system of congruences:
\begin{align*}
x_{1} & \overset{p_{1}}{\equiv}j_{1}\\
 & \vdots\\
x_{r} & \overset{p_{r}}{\equiv}j_{r}
\end{align*}
and where $\left[\mathbf{x}\overset{P}{\equiv}\mathbf{j}\right]$
is Iverson Bracket notation, being $1$ when $\mathbf{x}\overset{P}{\equiv}\mathbf{j}$
is true and $0$ otherwise.
\end{defn}
All of the techniques to be discussed in this series of papers apply
almost verbatim to the above multi-dimensional Hydra maps. In e-mail
correspondence with the author in 2017, K.R. Matthews explained that
multi-dimensional Hydra maps are almost entirely unstudied.

Though our methods are new, they are not entirely unanticipated. As
we show in \textbf{Section \ref{subsec:The-Tao-of}}, Terence Tao's
breakthrough 2019 result on Collatz (\cite{Tao Probability paper})
can be succinctly stated in terms of a kind of Fourier transform of
$\chi_{3}$; likewise \cite{Gonzales}'s generalization of \cite{Tao Probability paper}
can be stated in these same terms using the numen of the associated
Hydra maps. In this regard, $\chi_{3}$ is an independent and alternative
realization of Tao's \textbf{Syracuse Random Variables}.\textbf{ }These
same methods ought to be applicable to multi-dimensional Hydra maps
(and, again, would have the same Fourier-analytic correspondence to
those maps' associated numina\footnote{The plural form of ``numen''.}).
Such a result would, obviously, be a major leap in our understanding
of Hydra maps on $\mathbb{Z}^{d}$. This also supports our contention
that the numen formalism appears to be the ``natural'' approach
to studying Collatz-type dynamical systems.

As yet another example of the naturality of our approach, the numen
formalism also extends to handle Collatz-type maps in spaces of positive
characteristic. Matthews' slides give examples of Collatz-type maps
on rings of polynomials with coefficients in finite fields; the arguments
in this paper apply to that situation almost verbatim, and the author
intends to publish the details at a later date. For example, in his
slides (\cite{Matthews' slides}), Matthews defines a map $T:\mathbb{F}_{2}\left[X\right]\rightarrow\mathbb{F}_{2}\left[X\right]$
by:
\begin{equation}
T\left\{ f\right\} \overset{\textrm{def}}{=}\begin{cases}
\frac{f}{X} & \textrm{if }f\overset{X}{\equiv}0\\
\frac{\left(X^{2}+1\right)f+1}{X} & \textrm{if }f\overset{X}{\equiv}1
\end{cases},\textrm{ }\forall f\in\mathbb{F}_{2}\left[X\right]
\end{equation}
The numen of this map, $\chi_{T}$, accepts an element of $\mathbb{Z}_{2}$
as input and outputs an element of the ring of integers $R$ of the
completion of the field $\mathbb{F}_{2}\left(X\right)$ with respect
to the absolute value induced by the valuation associated to the prime
ideal $\left\langle X-1\right\rangle \subset\mathbb{F}_{2}\left[X\right]$.
A bit of algebra shows that this ring is isomorphic to $\mathbb{F}_{2}\left[\left[Y\right]\right]$,
integer ring of the local field $\mathbb{F}_{2}\left(\left(Y\right)\right)$
of formal Laurent series in $Y$ with coefficients in $\mathbb{F}_{2}$.
$\chi_{T}$ is characterized as the unique function $\mathbb{Z}_{2}\rightarrow R$
satisfying the functional equations:
\begin{align}
\chi_{T}\left(2\mathfrak{z}\right) & =\frac{\chi_{T}\left(\mathfrak{z}\right)}{X}\\
\chi_{T}\left(2\mathfrak{z}+1\right) & =\frac{\left(X^{2}+1\right)\chi_{T}\left(\mathfrak{z}\right)+1}{X}
\end{align}
subject to the rising-continuity condition:
\begin{equation}
\chi_{T}\left(\mathfrak{z}\right)\overset{R}{=}\lim_{n\rightarrow\infty}\chi_{T}\left(\left[\mathfrak{z}\right]_{2^{n}}\right),\textrm{ }\forall\mathfrak{z}\in\mathbb{Z}_{2}
\end{equation}
where the convergence is in the topology of $R$. The isomorphic copy
of $\chi_{T}$ taking values in $\mathbb{F}_{2}\left(\left(Y\right)\right)$
can be written as:
\begin{align}
\chi_{T}\left(2\mathfrak{z}\right) & =\frac{\chi_{T}\left(\mathfrak{z}\right)}{Y-1}\\
\chi_{T}\left(2\mathfrak{z}+1\right) & =\frac{Y\chi_{T}\left(\mathfrak{z}\right)+1}{Y-1}
\end{align}
and is the unique function from $\mathbb{Z}_{2}$ to $\mathbb{F}_{2}\left(\left(Y\right)\right)$
satisfying these functional equations subject to the rising continuity
condition:
\begin{equation}
\chi_{T}\left(\mathfrak{z}\right)\overset{\mathbb{F}_{2}\left(\left(Y\right)\right)}{=}\lim_{n\rightarrow\infty}\chi_{T}\left(\left[\mathfrak{z}\right]_{2^{n}}\right),\textrm{ }\forall\mathfrak{z}\in\mathbb{Z}_{2}
\end{equation}
where the convergence is in the topology of $\mathbb{F}_{2}\left(\left(Y\right)\right)$.
Replacing $Y$ with $X+1$ yields the previous version. Once the author
works out the details of the Fourier analysis of such functions, it
would be interesting to see if the number-theoretic rule of thumb
that function fields are easier to work in than number fields also
extends to include Collatz-type Conjectures. Moreover, the methods
of $\left(p,q\right)$-adic analysis described in this series of papers
ought to have analogues when the co-domain of the functions being
studied is a local field of positive characteristic, rather than a
local field of characteristic zero.

Because working at this level of generality would entail a cumbersome
amount of technical detail, for this series of papers, we will restrict
our analysis to a simple family of Hydra maps, the \textbf{Shortened
$qx+1$ maps}.
\begin{defn}
Let $q$ be an odd prime. The \textbf{Shortened $qx+1$ map }is the
function $T_{q}:\mathbb{Z}\rightarrow\mathbb{Z}$ defined by:
\begin{equation}
T_{q}\left(n\right)\overset{\textrm{def}}{=}\begin{cases}
\frac{n}{2} & \textrm{if }n\overset{2}{\equiv}0\\
\frac{qn+1}{2} & \textrm{if }n\overset{2}{\equiv}1
\end{cases}\label{eq:Shortened qx+1 Map}
\end{equation}
\end{defn}
The $T_{q}$s are an ideal ``toy'' case of the more general theory
described above. They allow for variable behavior while keeping those
same variations under the control of a single parameter ($q$). In
this regard, we have a direct parallel with the notion of a \textbf{bifurcation
}of a dynamical system, where one studies the system's long-term behavior
as a function of one or more of its scalar parameters \cite{Strogatz}.
Indeed, it will be very useful to view the $T_{q}$s as a family of
dynamical systems parameterized by $q$. In certain cases, we will
be able to treat $q$ as a real-valued or even \emph{complex}-valued
parameter.

$T_{5}$ is arguably the second-most well-studied Hydra map after
$T_{3}$ and $C$, and is somewhat infamous for how starkly its behavior
contrasts with $T_{3}$'s. In the 1970s, Riho Terras developed groundbreaking
probabilistic tools (\textbf{parity vectors})\textbf{ }for proving
that the set of all divergent points of $T_{3}$ in $\mathbb{N}_{0}$
has density $0$ \cite{Terras 76,Terras 79,Lagarias-Kontorovich Paper}.
These same tools can be used to show that the set of all divergent
points of $T_{5}$ in $\mathbb{N}_{0}$ has density $1$. Despite
this, not a \emph{single} integer has been proven to be a divergent
point of $T_{5}$ \cite{my dissertation,Lagarias-Kontorovich Paper}!
($7$ is believed to be the smallest positive integer which $T_{5}$
iterates to $\infty$, thoug, of course, no proof of this is currently
known). The same preponderance of divergent points also applies to
$T_{q}$ for $q\geq7$.

Just like with $\chi_{3}$, we associate a numen $\chi_{q}:\mathbb{Z}_{2}\rightarrow\mathbb{Z}_{q}$
to each $T_{q}$. As one would hope for a generalization of $\chi_{3}$,
it can be shown (\textbf{Lemma} \textbf{\ref{lem:Chi_H functional equation on N_0 and uniqueness}}
on page \pageref{lem:Chi_H functional equation on N_0 and uniqueness})
that $\chi_{q}$ is the unique function $f:\mathbb{Z}_{2}\rightarrow\mathbb{Z}_{q}$
satisfying the functional equations:
\begin{align}
f\left(2\mathfrak{z}\right) & =\frac{f\left(\mathfrak{z}\right)}{2}\\
f\left(2\mathfrak{z}+1\right) & =\frac{qf\left(\mathfrak{z}\right)+1}{2}
\end{align}
for all $\mathfrak{z}\in\mathbb{Z}_{2}$, subject to the rising-continuity
condition:
\begin{equation}
\lim_{n\rightarrow\infty}f\left(\mathfrak{z}\right)\overset{\mathbb{Z}_{q}}{=}f\left(\left[\mathfrak{z}\right]_{2^{n}}\right),\textrm{ }\forall\mathfrak{z}\in\mathbb{Z}_{2}
\end{equation}
Even better, $\chi_{q}$ satisfies the Correspondence Principle:
\begin{quotation}
\vphantom{}

\emph{I. (}\textbf{\emph{Corollary \ref{cor:CP v4}}}\emph{, page
\pageref{cor:CP v4}) $x\in\mathbb{Z}\backslash\left\{ 0\right\} $
is a periodic point of $T_{q}$ }if and only if\emph{ there exists
a $\mathfrak{z}\in\mathbb{Z}_{2}^{\prime}\cap\mathbb{Q}$ so that
$\chi_{q}\left(\mathfrak{z}\right)=x$. (}\textbf{\emph{CPPP}}\emph{)}

\emph{\vphantom{}}

\emph{II. (}\textbf{\emph{Theorem \ref{thm:Divergent trajectories come from irrational z}}}\emph{,
page \pageref{thm:Divergent trajectories come from irrational z})
If there is a $\mathfrak{z}\in\mathbb{Z}_{2}\backslash\mathbb{Q}$
so that $\chi_{q}\left(\mathfrak{z}\right)\in\mathbb{Z}$, then $\chi_{q}\left(\mathfrak{z}\right)$
is a divergent point of $T_{q}$. (}\textbf{\emph{CPDP}}\emph{)}
\end{quotation}
\begin{rem}
Unfortunately, \textbf{Theorem \ref{thm:Divergent trajectories come from irrational z}}
is incorrect. See \textbf{Section \ref{sec:Corrections} }for an explicit
counterexample and further discussion.
\end{rem}
Not only is this version of the CP\textbf{ }the main result of the
present paper, it also motivates all of the subsequent papers of this
series. The interest of this result lies in its novelty: it reframes
the study of $T_{q}$'s dynamics in terms of the \textbf{value distribution}
of $\chi_{q}$\textemdash which is to say, the study of the values
attained by $\chi_{q}\left(\mathfrak{z}\right)$ as $\mathfrak{z}$
varies in $\mathbb{Z}_{2}$. Investigating $\chi_{q}$'s value distribution
is far from elementary; moreover, this series of papers makes no claims
at solving Collatz or any related problems. If the methods of this
series end up playing a role in an eventual resolution of the Collatz
Conjecture, it will only be because of significant advances which
will have occurred in the interrim. But that is not our purpose here.

As we will see, the $\left(p,q\right)$-adic analytic setting utilized
by the later papers in this series abounds with tractable secondary
and tertiary questions, worthy of exploration in their own right,
independent of their potential applicability to $\chi_{q}$ and Collatz.
While it would be ideal if a better understanding of $\left(p,q\right)$
would lead to a proof of Collatz, we cannot know for certain if this
ideal can be made reality until a proof (or disproof) of Collatz is
finally produced. Still, this ought to have no bearing on whether
or not $\left(p,q\right)$ itself is worthy of study: it is a new
subject with mysteries and tantalizing connections all its own. That,
in it of itself, makes it deserving of our attention.
\begin{quotation}
\newpage{}
\end{quotation}

\section*{An Outline of this Paper}

After opening with a motivating example, \textbf{Section 1.1} introduces
the \textbf{string formalism }we will use in our study of $T_{q}$.
\textbf{Section 1.2} then constructs $\chi_{q}$, first as a $\mathbb{Q}$-valued
function of strings, and then as a function $\mathbb{N}_{0}\rightarrow\mathbb{Q}$.
\textbf{Section 1.2} also establishes the system of functional equations
that characterize $\chi_{q}:\mathbb{N}_{0}\rightarrow\mathbb{Q}$,
which is then used to show that $\chi_{q}$ admits a unique continuation
to a rising-continuous $\left(2,q\right)$-adic function $\chi_{q}:\mathbb{Z}_{2}\rightarrow\mathbb{Z}_{q}$.
We prove that this continuation satisfies the natural extension of
the aforementioned functional equations, and that those equations,
coupled with the rising-continuity condition, provide a complete characterization
of $\chi_{q}:\mathbb{Z}_{2}\rightarrow\mathbb{Z}_{q}$.

\textbf{Section 2} is dedicated to a proof the \textbf{Correspondence
Principle} (CP). This section is broken into two subsections: \textbf{Section
2.1}, in which our arguments' heuristics are presented and the necessary
preparatory work is done, and \textbf{Section 2.2}, in which the actual
proof(s) of the CP is/are given.

Much in the manner of this introduction, \textbf{Section 3} continues
in the vein suggested by our epigrammatic quote from D'Arcy Thompson.
\textbf{Section 3} provides a brief tour of how our numen formalism
connects with pre-existing literature on Collatz-type maps. We will
also point out several possible avenues for future exploration, though
more are to follow in the later papers of this series.

\textbf{Section 3.1} gives an account of how our methods dovetail
with Tao's work \cite{Tao Probability paper}. In particular, we elucidate
the relationship between $\chi_{3}$ and Tao's \textbf{Syracuse Random
Variables}. This is done by restating the key proposition of \cite{Tao Probability paper}\textemdash a
decay estimate on the characteristic function of the Syracuse Random
Variables\textemdash in terms of a simple variant of the Fourier transform
of $\chi_{3}$. \textbf{Section 3.2} gives a minor results about reverse
trajectories of points under $T_{q}$ and their attainability as elements
of the image of $\chi_{q}$. The paper concludes with \textbf{Section
3.3}, which gives a somewhat informal discussion of the implications
the Correspondence Principle has for the previously mentioned connection
between the Collatz Conjecture (specifically, the \textbf{Weak Collatz
Conjecture}) and improvements to lower bounds for $\left|q^{m}-2^{n}\right|$.
\begin{rem}
As it regards the relevant background material, \cite{my dissertation}
contains a thorough account of the history of work on the Collatz
Conjecture, as well as a discussion of some of the broader motivations
and goals of this research paradigm. Interested readers should consult
it. Lagarias' survey \cite{Lagarias' Survey} is also an excellent
introduction, although it lacks the broad scope of \cite{my dissertation}.
In that respect, K.R. Matthews' slides \cite{Matthews' slides} on
Markov-chain approaches to Collatz and related arithmetical dynamical
systems is much closer to \cite{my dissertation} in the larger variety
of Collatz-type maps it considers.
\end{rem}
\vphantom{}

This series of paper is a distillation of the highlights of the author's
PhD (Mathematics) dissertation \cite{my dissertation} done at the
University of Southern California under the obliging supervision of
Professors Sheldon Kamienny and Nicolai Haydn. Thanks must also be
given to Jeffery Lagarias, Steven J. Miller, Alex Kontorovich, Andrei
Khrennikov, K.R. Matthews, Susan Montgomery, Amy Young and all the
helpful staff of the USC Mathematics Department, and all the kindly
strangers became and acquainted with along the way.

\vphantom{}

This research did not receive any specific grant from funding agencies
in the public, commercial, or not-for-profit sectors.

\vphantom{}

Declarations of interest: none.

\section{\label{sec:1 - The Construction of Chi_q}The Construction of $\chi_{q}$}

IN THIS SECTION, WE WRITE $H$ TO DENOTE $T_{q}$, WHERE $q$ IS AN
ARBITRARY PRIME NUMBER. WE WILL ALSO USE $T_{q}$ TO DENOTE $T_{q}$.

\subsection{Strings}

This subsection details the formalism of \textbf{strings}, a foundational
notion in formal language theory and theoretical computer science
(particularly the study of automata); see, for instance, the first
chapter of \cite{Automatic sequences}. Our notation deviates somewhat
from the standard conventions in those subjects: for example, we write
$\textrm{String}\left(p\right)$ to denote what \cite{Automatic sequences}
writes as $\Sigma_{p}^{*}$, and $\textrm{String}_{\infty}\left(p\right)$
to denote what \cite{Automatic sequences} writes as $\Sigma_{p}^{\infty}$.
We also write strings as bold letters ($\mathbf{i}$, $\mathbf{j}$,
etc.), rather than italics ($x$, $y$, etc.), and denote string concatenation
by a binary operator ($\wedge$), rather than by mere juxtaposition
(thus, $\mathbf{i}\wedge\mathbf{j}$, rather than $xy$). We choose
these conventions primarily for their psychological benefits, benefits
which should hopefully become apparent within a couple of pages. Before
we introduce our notation, however, let us see why it is needed.
\begin{defn}
For each $j\in\left\{ 0,1\right\} $, we write $H_{j}:\mathbb{Q}\rightarrow\mathbb{Q}$
to denote the \textbf{$j$th branch} of $H$, defined as the function:
\begin{equation}
H_{j}\left(x\right)\overset{\textrm{def}}{=}\frac{a_{j}x+b_{j}}{2}\label{eq:Definition of H_j}
\end{equation}
where $a_{0}=1$, $b_{0}=0$ and $a_{1}=q$, $b_{1}=1$; thus:
\begin{align}
H_{0}\left(x\right) & =\frac{x}{2}\\
H_{1}\left(x\right) & =\frac{qx+1}{2}
\end{align}
We also call $H_{0}$ and $H_{1}$ the \textbf{even }and \textbf{odd
}branches of $H$, respectively.
\end{defn}
Given an integer $n$, observe that the sequence $n,H\left(n\right),H\left(H\left(n\right)\right),\ldots$
of iterates of $n$ under $H$ can be expressed in terms of compositions
of $H_{0}$ and $H_{1}$. Obviously, the particular composition sequence
of $H_{0}$s and $H_{1}$s we get depends on the $n$ we started with,
and the challenge of studying the dynamics of $H$ comes from the
difficulty of untangling the relation between $n$ and what iteration
of $H$ does to $n$. To that end, our work here begins with a question:
\begin{quote}
What happens if we considered \emph{arbitrary }composition sequences
of $H_{0}$ and $H_{1}$? Can we get something useful out of that?
\end{quote}
Let's find out.
\begin{example}
\label{exa:1}Classically, one of the key statistics of the Collatz
map (or, for that matter, $T_{q}$) is Terras' notion of a \textbf{parity
vector}. Given an integer $x$, the parity vector of $x$ under $H$
is the sequence generated by the values, mod $2$, of the iterates
of $x$ under $H$: $\left\{ \left[H^{\circ k}\left(x\right)\right]_{2}\right\} _{k\geq0}$.
The parity vector encodes the ``motions'' of $x$ under $H$, because
the value of $H^{\circ k}\left(x\right)$ mod $2$ determines which
branch of $H$ we will apply to $H^{\circ k}\left(x\right)$ in order
to obtain $H^{\circ k+1}\left(x\right)$.

For this example, let us fix $q=3$. Then, setting $x=3$ we have:
\begin{align}
H\left(3\right) & =H_{1}\left(3\right)=\frac{3\times3+1}{2}=5\\
H^{\circ2}\left(3\right) & =H\left(5\right)=H_{1}\left(5\right)=\frac{3\times5+1}{2}=8\\
H^{\circ3}\left(3\right) & =H\left(8\right)=H_{0}\left(8\right)=\frac{8}{2}=4
\end{align}
and so:
\begin{equation}
H^{\circ3}\left(3\right)=H_{0}\left(H_{1}\left(H_{1}\left(3\right)\right)\right)
\end{equation}
In this way, given any $x\in\mathbb{Z}$ and any $k\geq1$, there
is a \emph{unique} sequence $j_{1},\ldots,j_{k}$ of $0$s and $1$s
so that:
\begin{equation}
H^{\circ k}\left(n\right)=\left(H_{j_{1}}\circ H_{j_{2}}\circ\cdots\circ H_{j_{k}}\right)\left(n\right)
\end{equation}
Since an integer $x$ is a periodic point of $H$ if and only if there
is a $k\geq1$ so that $H^{\circ k}\left(x\right)=x$, by replacing
$H^{\circ k}\left(x\right)$ with the \textbf{composition sequence
}$\left(H_{j_{1}}\circ H_{j_{2}}\circ\cdots\circ H_{j_{k}}\right)\left(x\right)$,
we get:
\begin{align}
x & =H^{\circ k}\left(x\right)\nonumber \\
 & \Updownarrow\nonumber \\
x & =\left(H_{j_{1}}\circ H_{j_{2}}\circ\cdots\circ H_{j_{k}}\right)\left(x\right)\label{eq:Example for strings}
\end{align}
\emph{Classically}, the sequence $j_{1},\ldots,j_{k}$ (the first
$k$ elements of $x$'s parity vector, albeit in reverse order) would,
like the parity vector from which it is derived, be viewed \emph{a
posteriori }to the value of $x$ that generated it. First we find
the periodic point $x$, and then, by applying $H$ to $x$, we discover
the finite number of $0$s and $1$s $j_{1},\ldots,j_{k}\in\left\{ 0,1\right\} $
(where $k\geq1$) so that $x=\left(H_{j_{1}}\circ H_{j_{2}}\circ\cdots\circ H_{j_{k}}\right)\left(x\right)$.

But, what if we reversed this relationship? Can we obtain

Well, given any finite-length sequence of $0$s and $1$s $j_{1},\ldots,j_{k}\in\left\{ 0,1\right\} $,
since $H_{j_{1}},\ldots,H_{j_{k}}$ are maps of the form $r\mapsto ar+b$
(where $r\in\mathbb{Q}$), there are constants $a$ and $b$ (depending
on the $j$s we were given) so that:
\begin{equation}
\left(H_{j_{1}}\circ H_{j_{2}}\circ\cdots\circ H_{j_{k}}\right)\left(r\right)=ar+b,\textrm{ }\forall r\in\mathbb{Q}
\end{equation}
(This also holds for $r\in\mathbb{R}$, but never mind that.) Consequently,
we can solve $r=\left(H_{j_{1}}\circ H_{j_{2}}\circ\cdots\circ H_{j_{k}}\right)\left(r\right)$
and thereby obtain an expression for the composition sequence's fixed
point in terms of the $j$s:
\begin{align}
x & =\left(H_{j_{1}}\circ H_{j_{2}}\circ\cdots\circ H_{j_{k}}\right)\left(x\right)\\
 & \Updownarrow\nonumber \\
x & =ax+b\nonumber \\
 & \Updownarrow\nonumber \\
x & =\frac{b}{1-a}
\end{align}
Since $a\neq1$, this method always gives an answer, and, moroever,
the answer it gives is the \emph{unique }rational number fixed by
the composition sequence specified by our choice of $j$s.

To give a specific example, letting $q$ be arbitrary, and picking
the composition sequence $H_{1}\circ H_{1}\circ H_{0}\circ H_{1}$,
the unique rational number $x$ fixed is the solution of the equation:
\begin{equation}
x=\left(H_{1}\circ H_{1}\circ H_{0}\circ H_{1}\right)\left(x\right)=\frac{q^{3}}{16}x+\frac{q^{2}}{16}+\frac{q}{4}+\frac{1}{2}
\end{equation}
and so:
\begin{equation}
x=\frac{\frac{q^{2}}{16}+\frac{q}{4}+\frac{1}{2}}{1-\frac{q^{3}}{16}}=\frac{q^{2}+4q+8}{16-q^{3}}
\end{equation}
When $q=3$, this is $-29/11$. Since $-29/11$ is not an integer,
we can conclude that $T_{3}$ has no periodic point in $\mathbb{Z}$
which is fixed by the composition sequence $H_{1}\circ H_{1}\circ H_{0}\circ H_{1}$.

In this way, we see that (\ref{eq:Example for strings}) is an explicit
expression for periodic points of $T_{q}$. Moreover, we can study
this expression by viewing it as a function of our choice of $j$s.
\end{example}
In it of themselves, these insights are not new. For the specific
case of $H=T_{3}$, writing out the right-hand side of:
\begin{equation}
x=\left(H_{j_{1}}\circ\cdots\circ H_{j_{k}}\right)\left(x\right)
\end{equation}
in terms of powers of $2$ and $3$ produces a Diophantine equation
criterion for periodic points of the Collatz map. This was first discovered
by Böhm and Sontacchi in 1978 \cite{B=0000F6hm and Sontacchi} (though,
see also \cite{Wirsching's book on 3n+1}):
\begin{thm}[\textbf{Böhm-Sontacchi Criterion}\footnote{This terminology is the author's.}
\textbf{(1978)}]
\label{thm:Bohm-Sontacchi}Let $x\in\mathbb{Z}$ be a periodic point
of the Collatz map $C$ (\ref{eq:Collatz Map}). Then, there are integers
$m,n\geq1$ and a strictly increasing sequence of positive integers
$b_{1}<b_{2}<\cdots<b_{n}$ so that: 
\begin{equation}
x=\frac{\sum_{k=1}^{n}2^{m-b_{k}-1}3^{k-1}}{2^{m}-3^{n}}\label{eq:The Bohm-Sontacchi Criterion}
\end{equation}
\end{thm}
\begin{rem}
More general versions of this formula can be found in Matthews' slides
\cite{Matthews' slides}.
\end{rem}
While the existence of a Diophantine equation governing periodic points
of $C$ and Collatz-type maps in general is not new, our approach
differs from past work in a key way: instead of viewing the exponents
of $2$ and $3$ in the Böhm-Sontacchi Criterion as unknowns to be
determined, we can, in the manner of \textbf{Example \ref{exa:1}},
\emph{parameterize }those numbers in terms of a single variable: our
choice of $j$s. Specifically, the act of solving the equation (\ref{eq:Example for strings})
for $x$ from a given sequence of $j$s defines a scalar-valued function
of a single, string-valued variable.
\begin{defn}[\textbf{Strings}]
Let $p$ be an integer $\geq2$. For each $n\in\mathbb{N}_{1}$,
we write $\textrm{String}\left(p\right)$ to denote the set of all
finite sequences whose entries belong to the set $\left\{ 0,\ldots,p-1\right\} $
(``\textbf{strings}''). We also include the empty set ($\varnothing$)
as an element of $\textrm{String}\left(p\right)$, and call it the
\textbf{empty string}. We denote an arbitrary string as $\mathbf{j}=\left(j_{1},\ldots,j_{n}\right)$,
where $n$ is the \textbf{length }of the string; we denote the length
of $\mathbf{j}$ by $\left|\mathbf{j}\right|$, and define the empty
string as the unique string of length $0$. In this manner, any $\mathbf{j}$
can be written as $\mathbf{j}=\left(j_{1},\ldots,j_{\left|\mathbf{j}\right|}\right)$.
Also, we say a string $\mathbf{j}$ is \textbf{non-zero }when it contains
at least one non-zero entry.

We write $\textrm{String}_{\infty}\left(p\right)$ to denote the set
of all all sequences (finite \emph{or }infinite) whose entries belong
to the set $\left\{ 0,\ldots,p-1\right\} $. As with $\textrm{String}\left(p\right)$,
we also include the empty string ($\varnothing$) as an element of
$\textrm{String}_{\infty}\left(p\right)$. We call the elements of
$\textrm{String}_{\infty}\left(p\right)$ \textbf{strings} as well.
\end{defn}
With strings at our disposal, we can succinctly state composition
sequences of the branches of $H$.
\begin{defn}[\textbf{Composition sequences}]
Given any $\mathbf{j}\in\textrm{String}\left(2\right)$, we define
the \textbf{composition sequence }$H_{\mathbf{j}}:\mathbb{Q}\rightarrow\mathbb{Q}$
as the affine linear map:
\begin{equation}
H_{\mathbf{j}}\left(x\right)=\left(H_{j_{1}}\circ\cdots\circ H_{j_{\left|\mathbf{j}\right|}}\right)\left(x\right),\textrm{ }\forall\mathbf{j}\in\textrm{String}\left(2\right),\forall x\in\mathbb{Q}\label{eq:Def of composition sequence}
\end{equation}
In particular, we identify $H_{\varnothing}$ with the identity map;
$H_{\varnothing}\left(x\right)=x$ for all $x\in\mathbb{Q}$.
\end{defn}
While it might initially seem unusual to have our \emph{first }index
correspond to the branch we apply \emph{last}, this convention turns
out to be as vital psychologically as it is mathematically. Our motivation
here is the natural connection to $p$-adic integers, viewed here\textemdash as
Hensel himself viewed them\textemdash as formal power series in $p$
with coefficients in $\left\{ 0,\ldots,p-1\right\} $:
\begin{equation}
\mathfrak{z}=\sum_{n=0}^{\infty}d_{n}p^{n}=d_{0}+d_{1}p+d_{2}p^{2}+\cdots
\end{equation}
Here the $d_{n}$s are the \textbf{$p$-adic digits} of $\mathfrak{z}$,
written in left to right order, as is conventional for power series.
We adopt the left to right order in representing $\mathfrak{z}$ in
terms of their digits, with the above $\mathfrak{z}$ being written
as:
\begin{equation}
\mathfrak{z}=\centerdot_{p}d_{0}d_{1}d_{2}d_{3}\ldots
\end{equation}
where $\centerdot_{p}$ is the $p$-adic decimal point. The first
big trick of this paper is to transform the set $\textrm{String}\left(p\right)$
into the much more analytically manageable set $\mathbb{Z}_{p}$ by
identifying a given string with the sequence of digits of a $p$-adic
integer:
\begin{equation}
\mathfrak{z}\sim\underbrace{\left(d_{0},d_{1},d_{2},\ldots\right)}_{\mathbf{j}}
\end{equation}
The function $\textrm{Dig}\textrm{Sum}_{p}$ introduced below lets
us make this identification rigorous.
\begin{defn}
We write $\textrm{Dig}\textrm{Sum}_{p}:\textrm{String}_{\infty}\left(p\right)\rightarrow\mathbb{Z}_{p}$
to denote the function:
\begin{equation}
\textrm{DigSum}_{p}\left(\mathbf{j}\right)\overset{\textrm{def}}{=}\sum_{k=1}^{\left|\mathbf{j}\right|}j_{k}p^{k-1}\label{eq:Definition of DigSum_p of bold j}
\end{equation}
where $\left|\mathbf{j}\right|$ is defined to be $\infty$ if $\mathbf{j}$
is of infinite length.
\end{defn}
We go about identifying strings with numbers like so:
\begin{defn}[\textbf{Identification with $\mathbb{N}_{0}$}]
Given $\mathbf{j}\in\textrm{String}_{\infty}\left(p\right)$ and
$\mathfrak{z}\in\mathbb{Z}_{p}$, we say $\mathbf{j}$ \textbf{represents
}/ \textbf{is} \textbf{associated to }$\mathfrak{z}$ (and vice-versa),
written $\mathbf{j}\sim\mathfrak{z}$ or $\mathfrak{z}\sim\mathbf{j}$
whenever $\mathbf{j}$ is the sequence of the $p$-adic digits of
$\mathfrak{z}$; that is:
\begin{equation}
\mathfrak{z}\sim\mathbf{j}\Leftrightarrow\mathbf{j}\sim\mathfrak{z}\Leftrightarrow\mathfrak{z}=j_{1}+j_{2}p+j_{3}p^{2}+\cdots\label{eq:Definition of n-bold-j correspondence.}
\end{equation}
or, equivalently:
\begin{equation}
\mathbf{j}\sim\mathfrak{z}\Leftrightarrow\textrm{DigSum}_{p}\left(\mathbf{j}\right)=\mathfrak{z}
\end{equation}
Additionaly, we observe that $\sim$ induces an equivalence relation
on both $\textrm{String}\left(p\right)$ and $\textrm{String}_{\infty}\left(p\right)$
by the rule:
\begin{equation}
\mathbf{i}\sim\mathbf{j}\Leftrightarrow D_{p}\left(\mathbf{i}\right)=D_{p}\left(\mathbf{j}\right)\label{eq:Definition of string rho equivalence relation, rational integer version}
\end{equation}
which is to say, $\mathbf{i}\sim\mathbf{j}$ if and only if $\mathbf{i}$
and $\mathbf{j}$ represent the same $p$-adic integer.

Note that in both $\textrm{String}\left(p\right)$ and $\textrm{String}_{\infty}\left(p\right)$,
the shortest string representing the number $0$ is then the empty
string.

Lastly, we write $\textrm{String}\left(p\right)/\sim$ and $\textrm{String}_{\infty}\left(p\right)/\sim$
to denote the set of equivalence classes of $\textrm{String}\left(p\right)$
and $\textrm{String}_{\infty}\left(p\right)$, respectively, under
the equivalence relation (\ref{eq:Definition of string rho equivalence relation, rational integer version}).
\end{defn}
\begin{prop}
\label{prop:string number equivalence}\ 

I. $\textrm{String}\left(p\right)/\sim$ and $\textrm{String}_{\infty}\left(p\right)/\sim$
are in bijective correspondences with $\mathbb{N}_{0}$ and $\mathbb{Z}_{p}$,
respectively, by way of the map $\textrm{Dig}\textrm{Sum}_{p}$.

\vphantom{}

II. Two finite strings $\mathbf{i}$ and $\mathbf{j}$ satisfy $\mathbf{i}\sim\mathbf{j}$
if and only if $\mathbf{i}$ can be obtained from $\mathbf{j}$ (and
vice-versa) by adding or removing $0$s to the right of $\mathbf{j}$.

\vphantom{}

III. Given a finite string $\mathbf{i}$ and an infinite string $\mathbf{j}$,
$\mathbf{i}\sim\mathbf{j}$ can only occur if $\mathbf{j}$ contains
finitely many non-zero entries.
\end{prop}
Proof: Immediate from the definitions.

Q.E.D.

\vphantom{}

From an algebraic point of view, strings are useful because we can
concatenate them.
\begin{defn}[\textbf{Concatenation}]
We introduce the \textbf{concatenation operation} $\wedge:\textrm{String}_{\infty}\left(p\right)\times\textrm{String}_{\infty}\left(p\right)\rightarrow\textrm{String}_{\infty}\left(p\right)$,
defined by:
\begin{equation}
\mathbf{i}\wedge\mathbf{j}=\left(i_{1},\ldots,i_{\left|\mathbf{i}\right|}\right)\wedge\left(j_{1},\ldots,j_{\left|\mathbf{j}\right|}\right)\overset{\textrm{def}}{=}\begin{cases}
\left(i_{1},\ldots,i_{\left|\mathbf{i}\right|},j_{1},\ldots,j_{\left|\mathbf{j}\right|}\right) & \textrm{if }\left|\mathbf{i}\right|<\infty\\
\mathbf{i} & \textrm{if }\left|\mathbf{i}\right|=\infty
\end{cases}\label{eq:Definition of Concatenation}
\end{equation}
Additionally, for any integer $m\geq1$ and any finite string $\mathbf{j}$,
we write $\mathbf{j}^{\wedge m}$ to denote the concatenation of $m$
copies of $\mathbf{j}$:
\begin{equation}
\mathbf{j}^{\wedge m}\overset{\textrm{def}}{=}\left(\underbrace{j_{1},\ldots,j_{\left|\mathbf{j}\right|},j_{1},\ldots,j_{\left|\mathbf{j}\right|},\ldots,j_{1},\ldots,j_{\left|\mathbf{j}\right|}}_{m\textrm{ times}}\right)\label{eq:Definition of concatenation exponentiation}
\end{equation}
\end{defn}
Equipping $\textrm{String}\left(p\right)$ (or $\textrm{String}_{\infty}\left(p\right)$)
with the operation $\wedge$ makes it a monoid; the empty string is
the identity element. For us, strings are useful precisely because
our notation for composition sequence is, itself, a homomorphism of
monoids\footnote{This suggests it may be interesting to explore how the methods of
this section apply to the study of arbitrary sequences of compositions
of maps $\left\{ \phi_{0},\ldots,\phi_{p-1}\right\} $ on $\mathbb{Q}$.}:
\begin{prop}
\label{prop:H string formalism}\ 
\begin{equation}
H_{\mathbf{i}\wedge\mathbf{j}}\left(x\right)=H_{\mathbf{i}}\left(H_{\mathbf{j}}\left(x\right)\right),\textrm{ }\forall x\in\mathbb{R},\textrm{ }\forall\mathbf{i},\mathbf{j}\in\textrm{String}\left(2\right)\label{eq:H string formalism}
\end{equation}
\end{prop}
Proof: 
\begin{equation}
H_{\mathbf{i}\wedge\mathbf{j}}\left(x\right)=H_{i_{1},\ldots,i_{\left|\mathbf{i}\right|},j_{1},\ldots,j_{\left|\mathbf{j}\right|}}\left(x\right)=\left(H_{i_{1},\ldots,i_{\left|\mathbf{i}\right|}}\circ H_{j_{1},\ldots,j_{\left|\mathbf{j}\right|}}\right)\left(x\right)=H_{\mathbf{i}}\left(H_{\mathbf{j}}\left(x\right)\right)
\end{equation}

Q.E.D.

\vphantom{}

Before moving on to the next section, we need to introduce two functions
vital for working with elements of $\textrm{String}_{\infty}\left(2\right)$
and the $2$-adic integers they represent.
\begin{defn}[\textbf{$\lambda_{2}$ and $\#_{1}$}]
We write $\lambda_{2}:\mathbb{N}_{0}\rightarrow\mathbb{N}_{0}$ to
denote the function:
\begin{equation}
\lambda_{2}\left(n\right)\overset{\textrm{def}}{=}\left\lceil \log_{2}\left(n+1\right)\right\rceil ,\textrm{ }\forall n\in\mathbb{N}_{0}\label{eq:definition of lambda rho}
\end{equation}
which gives the number of $2$-adic (or, equivalently, binary) digits
of $n$, otherwise known as the number of \textbf{bits }in $n$ Note
that every integer $n\geq0$ can be uniquely written as:
\begin{equation}
n=c_{0}+c_{1}2+\cdots+c_{\lambda_{2}\left(n\right)-1}2^{\lambda_{2}\left(n\right)-1}
\end{equation}
where the $c_{j}$s are integers in $\left\{ 0,1\right\} $. Additionally,
note that:
\begin{equation}
\left\lceil \log_{2}\left(n+1\right)\right\rceil =\left\lfloor \log_{2}n\right\rfloor +1,\textrm{ }\forall n\in\mathbb{N}_{1}\label{eq:Floor and Ceiling expressions for lambda_rho}
\end{equation}
and that $\lambda_{2}\left(n\right)\leq\left|\mathbf{j}\right|$ is
satisfied for any string $\mathbf{j}$ representing $n$.

Next, for each $n\geq0$, we write $\#_{1}\left(n\right)$ to denote
the number of $1$s present in the $2$-adic expansion of $n$. In
an abuse of notation, we also use $\#_{1}$ to denote the number of
$1$s entries in a given string: 
\begin{equation}
\#_{1}\left(\mathbf{j}\right)\overset{\textrm{def}}{=}\textrm{number of }1\textrm{s in }\mathbf{j},\textrm{ }\forall\mathbf{j}\in\textrm{String}\left(2\right)\label{eq:Definition of number of ks in rho adic digits of bold j}
\end{equation}

Lastly, we note the functional equations satisfied by $\#_{1}$ and
$\lambda_{2}$.
\end{defn}
\begin{prop}[\textbf{Functional Equations for $\lambda_{2}$ and $\#_{1}$}]
\ 

I.
\begin{align}
\#_{1}\left(2^{n}a+b\right) & =\#_{1}\left(a\right)+\#_{1}\left(b\right),\textrm{ }\forall a\in\mathbb{N}_{0},\textrm{ }\forall n\in\mathbb{N}_{1},\textrm{ }\forall b\in\left\{ 0,1\right\} \label{eq:number-symbol functional equations}
\end{align}

\vphantom{}

II.
\begin{align}
\lambda_{2}\left(2^{n}a+b\right) & =\lambda_{2}\left(a\right)+n,\textrm{ }\forall a\in\mathbb{N}_{1},\textrm{ }\forall n\in\mathbb{N}_{1},\textrm{ }\forall b\in\left\{ 0,1\right\} \label{eq:lambda functional equations}
\end{align}
\end{prop}
Proof: A straightforward computation.

Q.E.D.

\subsection{\label{subsec:2.2.2 Construction-of}Construction of $\chi_{q}$}

With strings now at our disposal, we can begin a proper analysis of
the equation:
\begin{equation}
x=\left(H_{j_{1}}\circ\cdots\circ H_{j_{k}}\right)\left(x\right)
\end{equation}
To begin, since the composition of any finite number of affine linear
maps (such as the $H_{j}$s) is an affine linear map, given any $\mathbf{j}\in\textrm{String}\left(2\right)$,
we can express $H_{\mathbf{j}}\left(r\right)$ as:

\begin{equation}
H_{\mathbf{j}}\left(r\right)=H_{\mathbf{j}}^{\prime}\left(0\right)r+H_{\mathbf{j}}\left(0\right),\textrm{ }\forall\mathbf{j}\in\textrm{String}\left(2\right),\textrm{ }\forall r\in\mathbb{Q}\label{eq:ax+b formula for h_bold_j}
\end{equation}
where $H_{\mathbf{j}}^{\prime}\left(0\right)$ is the derivative of
$H_{\mathbf{j}}\left(r\right)$ evaluated at $r=0$; note that, by
the chain rule of differential calculus, $H_{\mathbf{j}}^{\prime}\left(0\right)$
is the product of $H_{j}^{\prime}\left(0\right)$ for each $j$ in
$\mathbf{j}$; also, the choice of $r=0$ is arbitrary, seeing as,
for fixed $\mathbf{j}$, $H_{\mathbf{j}}^{\prime}\left(r\right)$
is constant as a function of $r$.

Next, viewing $\mathbf{j}$ as the ``instructions'' for how to apply
the branches of $H$, we observe that a \emph{necessary }condition
for an integer $x$ to be a periodic point of $H$ is that there exists
a $\mathbf{j}\in\textrm{String}\left(2\right)$ of length $\geq1$
so that $H_{\mathbf{j}}\left(x\right)=x$. Plugging the condition
$H_{\mathbf{j}}\left(x\right)=x$ into (\ref{eq:ax+b formula for h_bold_j})
gives us:
\begin{equation}
x=H_{\mathbf{j}}^{\prime}\left(0\right)x+H_{\mathbf{j}}\left(0\right)\label{eq:affine formula for little x}
\end{equation}
Since $H_{j}^{\prime}\left(0\right)\neq1$ for all $j$, the chain
rule tells us that, for $\mathbf{j}\in\textrm{String}\left(2\right)$,
$H_{\mathbf{j}}^{\prime}\left(0\right)=1$ occurs if and only if $\mathbf{j}=\varnothing$.
Thus, since $\mathbf{j}$ has positive length, $H_{\mathbf{j}}^{\prime}\left(0\right)\neq1$,
and so there is a unique $x\in\mathbb{Q}$ satisfying (\ref{eq:affine formula for little x}).
This is given by:
\begin{equation}
x=\frac{H_{\mathbf{j}}\left(0\right)}{1-H_{\mathbf{j}}^{\prime}\left(0\right)}\label{eq:Formula for little x in terms of bold j}
\end{equation}
So, given $x\in\mathbb{Z}\backslash\left\{ 0\right\} $, the existence
of a $\mathbf{j}\in\textrm{String}\left(2\right)$ satisfying (\ref{eq:Formula for little x in terms of bold j})
is a \emph{necessary} condition for $x$ to be a periodic point of
$H$; this is the content of the Böhm-Sontacchi criterion. However,
we will take things further.

The significance of (\ref{eq:Formula for little x in terms of bold j})
is that it gives us a formula for the unique fixed point of $H_{\mathbf{j}}:\mathbb{Q}\rightarrow\mathbb{Q}$
\emph{in terms of} the string variable $\mathbf{j}$. Indeed, the
parameterization of $x$ in terms of $\mathbf{j}$ is the map:
\begin{equation}
\mathbf{j}\mapsto\frac{H_{\mathbf{j}}\left(0\right)}{1-H_{\mathbf{j}}^{\prime}\left(0\right)}\label{eq:Formula for X}
\end{equation}
By using $\textrm{DigSum}_{2}$ to identify $\mathbf{j}$s with non-negative
(and later, $2$-adic) integers, we open the door to analytical investigations
of the right-hand side of (\ref{eq:Formula for X}), first as a function
of a non-negative integer variable, and later, as a function of a
$2$-adic integer.

Now, while we \emph{could} analyze the map (\ref{eq:Formula for X})
directly, this turns out to be unwise; the expression in the denominator
is a troublemaker. Indeed, a simple computation shows that the value
of $H_{\mathbf{j}}^{\prime}\left(0\right)$ in terms of $\mathbf{j}$
is:
\begin{equation}
H_{\mathbf{j}}^{\prime}\left(0\right)=\frac{q^{\#_{1}\left(\mathbf{j}\right)}}{2^{\left|\mathbf{j}\right|}}
\end{equation}
and so, (\ref{eq:Formula for X}) will blow up for values of $\mathbf{j}$
which make this fraction close to $1$, which is to say, values of
$\mathbf{j}$ which make real number:
\begin{equation}
\left|\#_{1}\left(\mathbf{j}\right)\ln q-\left|\mathbf{j}\right|\ln2\right|
\end{equation}
close to 0 (in the archimedean sense). Controlling the size of this
number would take us directly into transcendental number theory (specifically,
\textbf{Baker's Theorem }on linear forms in logarithms)\textemdash a
can of wyrms in its own right. So, rather than bite off more than
we can chew and attempt to study (\ref{eq:Formula for X}) in its
entirety, we will instead focus on the its numerator, the map $\mathbf{j}\mapsto H_{\mathbf{j}}\left(0\right)$.
We call this the \textbf{numen }of $H=T_{q}$, and denote it by $\chi_{q}$;
the definition is given below. Even better, this decision turns out
to subsume the map (\ref{eq:Formula for X}), which we can recover
from $\chi_{q}$ for free by way of a functional equation \textbf{(Lemma
\ref{lem:Chi_H o B_p functional equation}}, from \textbf{Section
\ref{sec:2 The-Correspondence-Principle}}.
\begin{defn}
\label{def:Chi_H on N_0 in strings}The \textbf{numen of the Shortened
$qx+1$ map }is the function $\chi_{q}:\textrm{String}\left(2\right)\rightarrow\mathbb{Q}$
defined by: 
\begin{equation}
\chi_{q}\left(\mathbf{j}\right)\overset{\textrm{def}}{=}H_{\mathbf{j}}\left(0\right),\textrm{ }\forall\mathbf{j}\in\textrm{String}\left(2\right)\label{eq:Definition of Chi_H of bold j}
\end{equation}
\end{defn}
\begin{rem}
Since $H_{\varnothing}$ is the identity map on $\mathbb{Q}$, note
that $\chi_{q}\left(\varnothing\right)=0$.
\end{rem}
Next, we show that $\chi_{q}$ remains well-defined when we use $\sim$
to identify strings with non-negative integers.
\begin{prop}
\label{prop:Chi_H is well defined mod twiddle}$\chi_{q}$ is well-defined
on $\textrm{String}\left(2\right)/\sim$. That is, $\chi_{q}\left(\mathbf{i}\right)=\chi_{q}\left(\mathbf{j}\right)$
for all $\mathbf{i},\mathbf{j}\in\textrm{String}\left(2\right)$ for
which $\mathbf{i}\sim\mathbf{j}$.
\end{prop}
Proof: Let $\mathbf{j}\in\textrm{String}\left(2\right)$ be any non-empty
string. Then, by \textbf{Proposition \ref{prop:H string formalism}}:
\begin{equation}
\chi_{q}\left(\mathbf{j}\wedge0\right)=H_{\mathbf{j}\wedge0}\left(0\right)=H_{\mathbf{j}}\left(H_{0}\left(0\right)\right)=H_{\mathbf{j}}\left(0\right)=\chi_{q}\left(\mathbf{j}\right)
\end{equation}
since $H_{0}\left(0\right)=0$. By \textbf{Proposition \ref{prop:string number equivalence}},
given any two equivalent finite strings $\mathbf{i}\sim\mathbf{j}$,
we can obtain one of the two by concatenating finitely many $0$s
to the right of the other. As such, without loss of generality, suppose
that $\mathbf{i}=\mathbf{j}\wedge0^{\wedge n}$, where $0^{\wedge n}$
is a string of $n$ consecutive $0$s. Then: 
\begin{equation}
\chi_{q}\left(\mathbf{i}\right)=\chi_{q}\left(\mathbf{j}\wedge0^{\wedge n}\right)=\chi_{q}\left(\mathbf{j}\wedge0^{\wedge\left(n-1\right)}\right)=\cdots=\chi_{q}\left(\mathbf{j}\right)
\end{equation}
Hence, $\chi_{q}$ is well defined on $\textrm{String}\left(2\right)/\sim$.

Q.E.D.

\vphantom{}
\begin{lem}[$\chi_{q}$ \textbf{on} $\mathbb{N}_{0}$]
\label{lem:Chi_H on N_0}We can realize $\chi_{q}$ as a function
$\mathbb{N}_{0}\rightarrow\mathbb{Q}$ by defining: 
\begin{equation}
\chi_{q}\left(n\right)\overset{\textrm{def}}{=}\chi_{q}\left(\mathbf{j}\right)=H_{\mathbf{j}}\left(0\right)\label{eq:Definition of Chi_H of n}
\end{equation}
where $\mathbf{j}\in\textrm{String}\left(2\right)$ is any string
representing $n$. As such, $\chi_{q}\left(0\right)=\chi_{q}\left(\varnothing\right)=H_{\varnothing}\left(0\right)=0$.
\end{lem}
Proof: By \textbf{Proposition \ref{prop:string number equivalence}},
$\textrm{String}\left(2\right)/\sim$ is in a bijective correspondence
with $\mathbb{N}_{0}$. \textbf{Proposition \ref{prop:Chi_H is well defined mod twiddle}
}tells us that $\chi_{q}$ is well-defined on $\textrm{String}\left(2\right)/\sim$,
so, by using the aforementioned bijection to identify $\mathbb{N}_{0}$
with $\textrm{String}\left(2\right)/\sim$, the rule $\chi_{q}\left(n\right)\overset{\textrm{def}}{=}\chi_{q}\left(\mathbf{j}\right)$
is well-defined, being independent of our the string $\mathbf{j}\in\textrm{String}\left(2\right)$
chosen to represent $n$.

Q.E.D.

\vphantom{}

Although we will eventually shelve the string formalism in favor of
treating $\chi_{q}$ as a function on $\mathbb{N}_{0}$ (and, later,
$\mathbb{Z}_{2}$), strings have the advantage of \emph{greatly }simplifying
the proofs of the results that get our analysis of $\chi_{q}$ off
the ground. As such, for the remainder of this section, we will maintain
a bilingual dictionary, expressing everything in terms of both strings
and the numbers those strings represent.

First, we deal with $H_{\mathbf{j}}^{\prime}\left(0\right)$, that
troublesome term from (\ref{eq:Formula for X})'s denominator.
\begin{defn}
\label{def:M_H on N_0 in strings}We define $M_{q}:\mathbb{N}_{0}\rightarrow\mathbb{Q}$
by:
\begin{equation}
M_{q}\left(n\right)\overset{\textrm{def}}{=}M_{q}\left(\mathbf{j}\right)\overset{\textrm{def}}{=}H_{\mathbf{j}}^{\prime}\left(0\right),\textrm{ }\forall n\geq1\label{eq:Definition of M_H}
\end{equation}
where $\mathbf{j}$ is the shortest element of $\textrm{String}\left(2\right)$
representing $n$. We define $M_{q}\left(0\right)=M_{q}\left(\varnothing\right)$
to be $1$.
\end{defn}
\begin{example}
\emph{Unlike} $\chi_{q}$, observe that $M_{q}$ is \emph{not }well-defined
over $\textrm{String}\left(2\right)/\sim$. For example, the strings
$\left(0,1\right)$ and $\left(0,1,0\right)$ both represent the integer
$2$, but:
\begin{align}
H_{0,1}\left(x\right) & =H_{0}\left(H_{1}\left(x\right)\right)=\frac{q}{4}x+\frac{1}{4}\\
H_{0,1,0}\left(x\right) & =H_{0}\left(H_{1}\left(H_{0}\left(x\right)\right)\right)=\frac{q}{8}x+\frac{1}{4}
\end{align}
Thus, $\left(0,1\right)\sim\left(0,1,0\right)$, but $M_{q}\left(0,1\right)=q/4$,
which is not equal to $M_{q}\left(0,1,0\right)=q/8$.
\end{example}
With these definitions in place, we can now use strings to establish
the key functional equation identities for $\chi_{q}$ and $M_{q}$.
These identities will dominate all of the analysis to come, particularly
the \textbf{concatenation identities }of $\chi_{q}$ and $M_{q}$;
we will also sometimes refer to these as the \textbf{functional equations
}of $\chi_{q}$ and $M_{q}$.
\begin{prop}[\textbf{Expression for $H_{\mathbf{j}}$}]
\label{prop:H_boldj in terms of M_H and Chi_H} 
\begin{equation}
H_{\mathbf{j}}\left(x\right)=M_{q}\left(\mathbf{j}\right)x+\chi_{q}\left(\mathbf{j}\right),\textrm{ }\forall\mathbf{j}\in\textrm{String}\left(2\right),\textrm{ }\forall x\in\mathbb{Q}\label{eq:Formula for composition sequences of H}
\end{equation}
\end{prop}
Proof: $H_{\mathbf{j}}\left(x\right)=H_{\mathbf{j}}^{\prime}\left(0\right)x+H_{\mathbf{j}}\left(0\right)$.

Q.E.D.

\vphantom{}
\begin{lem}[$\chi_{q}$ \textbf{Concatenation Identity}]
\label{lem:Chi_H concatenation identity}\ 

\begin{equation}
\chi_{q}\left(\mathbf{i}\wedge\mathbf{j}\right)=H_{\mathbf{i}}\left(\chi_{q}\left(\mathbf{j}\right)\right),\textrm{ }\forall\mathbf{i},\mathbf{j}\in\textrm{String}\left(2\right)\label{eq:Chi_H concatenation identity}
\end{equation}
That is to say, for $\mathbf{i}=\left(i_{1},i_{2},\ldots,i_{m}\right)\in\textrm{String}\left(2\right)$,
where $m\geq1$, we have: 
\begin{equation}
\chi_{q}\left(i_{1}+i_{2}2+\cdots+i_{m}2^{m-1}+2^{m}n\right)=H_{\mathbf{i}}\left(\chi_{q}\left(n\right)\right),\textrm{ }\forall n\in\mathbb{N}_{0}\label{eq:Chi_H concatenation identity, alternate version}
\end{equation}
\end{lem}
Proof: Letting $\mathbf{i}$, $\mathbf{j}$, and $x$ be arbitrary,
we have that: 
\begin{equation}
H_{\mathbf{i}\wedge\mathbf{j}}\left(x\right)=H_{\mathbf{i}\wedge\mathbf{j}}^{\prime}\left(0\right)x+H_{\mathbf{i}\wedge\mathbf{j}}\left(0\right)=H_{\mathbf{i}\wedge\mathbf{j}}^{\prime}\left(0\right)x+H_{\mathbf{i}}\left(H_{\mathbf{j}}\left(0\right)\right)
\end{equation}
Setting $x=0$ yields: 
\begin{equation}
\underbrace{H_{\mathbf{i}\wedge\mathbf{j}}\left(0\right)}_{\chi_{q}\left(\mathbf{i}\wedge\mathbf{j}\right)}=\underbrace{H_{\mathbf{i}}\left(H_{\mathbf{j}}\left(0\right)\right)}_{H_{\mathbf{j}}\left(\chi_{q}\left(\mathbf{j}\right)\right)}
\end{equation}

Q.E.D.

\vphantom{}
\begin{prop}[$M_{q}$ \textbf{Concatenation Identity}]
\label{prop:M_H concatenation identity}\ 
\begin{equation}
M_{q}\left(\mathbf{i}\wedge\mathbf{j}\right)=M_{q}\left(\mathbf{i}\right)M_{q}\left(\mathbf{j}\right),\textrm{ }\forall\mathbf{i},\mathbf{j}\in\textrm{String}\left(2\right)\label{eq:M_H concatenation identity}
\end{equation}
That is to say, for $\mathbf{i}=\left(i_{1},i_{2},\ldots,i_{m}\right)\in\textrm{String}\left(2\right)$,
where $m\geq1$, we have: 
\begin{equation}
M_{q}\left(i_{1}+i_{2}2+\cdots+i_{m}2^{m-1}+2^{m}n\right)=M_{q}\left(\mathbf{i}\right)M_{q}\left(n\right),\textrm{ }\forall n\in\mathbb{N}_{0}\label{eq:Inductive identity for M_H}
\end{equation}
\end{prop}
Proof: Since our notation for composition sequences is a monoid homomorphism
(\textbf{Proposition \ref{prop:H string formalism}}), given $\mathbf{i},\mathbf{j}$,
we can write:
\begin{equation}
H_{\mathbf{i}\wedge\mathbf{j}}\left(x\right)=H_{\mathbf{i}}\left(H_{\mathbf{j}}\left(x\right)\right),\textrm{ }\forall x\in\mathbb{Q}
\end{equation}
Formally differentiating both sides with respect to $x$ and applying
the chain rule gives: 
\begin{equation}
H_{\mathbf{i}\wedge\mathbf{j}}^{\prime}\left(x\right)=H_{\mathbf{i}}^{\prime}\left(H_{\mathbf{j}}\left(x\right)\right)H_{\mathbf{j}}^{\prime}\left(x\right)
\end{equation}
Now, set $x=0$: 
\begin{equation}
\underbrace{H_{\mathbf{i}\wedge\mathbf{j}}^{\prime}\left(0\right)}_{M_{q}\left(\mathbf{i}\wedge\mathbf{j}\right)}=H_{\mathbf{i}}^{\prime}\left(H_{\mathbf{j}}\left(0\right)\right)\underbrace{H_{\mathbf{j}}^{\prime}\left(0\right)}_{M_{q}\left(\mathbf{j}\right)}
\end{equation}
Since $H_{\mathbf{i}}\left(x\right)$ is an affine linear map of the
form $H_{\mathbf{i}}^{\prime}\left(0\right)x+H_{\mathbf{i}}\left(0\right)$,
its derivative is the constant function $H_{\mathbf{i}}^{\prime}\left(y\right)=H_{\mathbf{i}}^{\prime}\left(0\right)$.
Thus: 
\begin{equation}
M_{q}\left(\mathbf{i}\wedge\mathbf{j}\right)=H_{\mathbf{i}}^{\prime}\left(H_{\mathbf{j}}\left(0\right)\right)M_{q}\left(\mathbf{j}\right)=H_{\mathbf{i}}^{\prime}\left(0\right)M_{q}\left(\mathbf{j}\right)=M_{q}\left(\mathbf{i}\right)M_{q}\left(\mathbf{j}\right)
\end{equation}

Q.E.D.
\begin{rem}
This makes $M_{q}$ into a homomorphism of monoids, sending $\left(\textrm{String}\left(2\right),\wedge\right)$
to $\left(\mathbb{Q}^{\times},\times\right)$, the monoid obtained
by equipping the set of non-zero rational numbers with ordindary multiplication
as the monoid operation.
\end{rem}
\vphantom{}

Next, we restate the concatenation identities in terms of systems
of functional equations.
\begin{prop}[\textbf{Functional Equations for} $M_{q}$]
\label{prop:M_H functional equation}\ 
\end{prop}
\begin{equation}
M_{q}\left(2n+j\right)=\frac{a_{j}}{2}M_{q}\left(n\right),\textrm{ }\forall n\geq0\textrm{ \& }\forall j\in\left\{ 0,1\right\} :2n+j\neq0\label{eq:M_H functional equations}
\end{equation}
where $a_{0}=1$ and $a_{1}=q$.
\begin{rem}
To be clear, this functional equation \emph{does not hold} when $2n+j=0$. 
\end{rem}
Proof: Let $n\geq1$, and let $\mathbf{j}=\left(j_{1},\ldots,j_{m}\right)$
be the shortest string in $\textrm{String}\left(2\right)$ representing
$n$. Letting $j\in\left\{ 0,1\right\} $ be arbitrary, setting $\mathbf{k}=\left(j,j_{1},j_{2},\ldots,j_{m}\right)$,
we have that $\mathbf{k}$ then represents $j+2n$:
\begin{equation}
j\wedge\mathbf{j}=\mathbf{k}\sim\left(j+2n\right)
\end{equation}
\textbf{Proposition \ref{prop:M_H concatenation identity}} lets us
write: 
\begin{equation}
M_{q}\left(2n+j\right)=M_{q}\left(j\wedge\mathbf{j}\right)=M_{q}\left(j\right)M_{q}\left(\mathbf{j}\right)=\frac{a_{j}}{2}M_{q}\left(n\right)
\end{equation}
When $n=0$, these equalities hold for $j=1$, seeing as: 
\begin{equation}
\left(1+2\cdot0\right)\sim j\wedge\varnothing
\end{equation}
As for the exceptional case ($n=j=0$), there, we obtain: 
\begin{equation}
M_{q}\left(2\cdot0+0\right)=M_{q}\left(0\right)\overset{\textrm{def}}{=}1
\end{equation}
but: 
\begin{equation}
\frac{a_{0}}{2}M_{q}\left(0\right)\overset{\textrm{def}}{=}\frac{a_{0}}{2}\times1=\frac{1}{2}\neq1
\end{equation}

Q.E.D.

\vphantom{}

The reason $\chi_{q}$'s concatenation identity is so important is
because it actually completely characterization of $\chi_{q}$.
\begin{lem}[\textbf{Functional Equations for} $\chi_{q}$]
\label{lem:Chi_H functional equation on N_0 and uniqueness}$\chi_{q}$
is the unique rational-valued function on $\mathbb{N}_{0}$ satisfying
the system of functional equations: 
\begin{align}
\chi_{q}\left(2n+j\right) & =H_{j}\left(\chi_{q}\left(n\right)\right),\textrm{ }\forall n\in\mathbb{N}_{0},\textrm{ }\forall j\in\left\{ 0,1\right\} \label{eq:Chi_H functional equations}
\end{align}
\end{lem}
\begin{rem}
These functional equations can also be written as:
\begin{align}
\chi_{q}\left(2n\right) & =\frac{\chi_{q}\left(n\right)}{2}\\
\chi_{q}\left(2n+1\right) & =\frac{q\chi_{q}\left(n\right)+1}{2}
\end{align}
for all $n\in\mathbb{N}_{0}$.
\end{rem}
\begin{rem}
This lemma is equivalent to the statement that $\chi_{q}:\textrm{String}\left(2\right)\rightarrow\mathbb{Q}$
is the \emph{unique} function $\textrm{String}\left(2\right)\rightarrow\mathbb{Q}$
satisfying (\ref{eq:Chi_H concatenation identity}).
\end{rem}
Proof:

I. Let $\mathbf{i}\sim n$ and let $j\in\left\{ 0,1\right\} $ be
arbitrary. Then $2n+j\sim j\wedge\mathbf{i}$, and hence, by $\chi_{q}$'s
concatenation identity (\textbf{Lemma \ref{lem:Chi_H concatenation identity}}):
\begin{equation}
\chi_{q}\left(2n+j\right)=\chi_{q}\left(j\wedge\mathbf{i}\right)=H_{j}\left(\chi_{q}\left(\mathbf{i}\right)\right)=H_{j}\left(\chi_{q}\left(n\right)\right)
\end{equation}
Thus, $\chi_{q}$ is a solution of the given system of functional
equations. The reason why we need not exclude the case where $n=j=0$
(as we had to do with $M_{q}$) is because $H_{0}\left(\chi_{q}\left(0\right)\right)=H_{0}\left(0\right)=0$.

\vphantom{}

II. On the other hand, let $f:\mathbb{N}_{0}\rightarrow\mathbb{Q}$
be any function so that: 
\begin{equation}
f\left(2n+j\right)=H_{j}\left(f\left(n\right)\right),\textrm{ }\forall n\in\mathbb{N}_{0},\textrm{ }\forall j\in\left\{ 0,1\right\} 
\end{equation}
Setting $n=j=0$ gives: 
\begin{equation}
f\left(0\right)=H_{0}\left(f\left(0\right)\right)=\frac{f\left(0\right)}{2}
\end{equation}
This forces $f\left(0\right)=0$. Plugging in $n=0$ leaves us with:
\begin{equation}
f\left(j\right)=H_{j}\left(f\left(0\right)\right)=H_{j}\left(0\right),\textrm{ }\forall j\in\left\{ 0,1\right\} 
\end{equation}
Then, writing $n$ as a $2$-adic Hensel series: 
\begin{equation}
n=j_{1}+j_{2}2+\cdots+j_{L}2^{L-1}
\end{equation}
we see that the identity $f\left(2n+j\right)=H_{j}\left(f\left(n\right)\right)$
is equivalent to: 
\begin{align*}
f\left(j+j_{1}2+j_{2}2^{2}+\cdots+j_{L}2^{L}\right) & =H_{j}\left(f\left(j_{1}+j_{2}2+\cdots+j_{L}2^{L-1}\right)\right)\\
 & =H_{j}\left(H_{j_{1}}\left(f\left(j_{2}+j_{q}2+\cdots+j_{L}2^{L-2}\right)\right)\right)\\
 & \vdots\\
 & =\left(H_{j}\circ H_{j_{1}}\circ\cdots\circ H_{j_{L}}\right)\left(f\left(0\right)\right)\\
 & =\left(H_{j}\circ H_{j_{1}}\circ\cdots\circ H_{j_{L}}\right)\left(0\right)\\
 & =H_{j,j_{1},\ldots,j_{L}}\left(0\right)
\end{align*}
So, for any string $\mathbf{j}^{\prime}=\left(j,j_{1},\ldots,j_{L}\right)$,
we have: 
\begin{equation}
f\left(\mathbf{j}^{\prime}\right)=H_{\mathbf{j}^{\prime}}\left(0\right)\overset{\textrm{def}}{=}\chi_{q}\left(\mathbf{j}^{\prime}\right)
\end{equation}
where, note, $\mathbf{j}^{\prime}=j\wedge\mathbf{j}$ and $\mathbf{j}^{\prime}\sim n$.
In other words, if $f$ solves the given system of functional equations,
$f\left(\mathbf{j}\right)=\chi_{q}\left(\mathbf{j}\right)$ for all
$\mathbf{j}\in\textrm{String}\left(2\right)$. Thus, $f\left(n\right)=\chi_{q}\left(n\right)$,
which proves the uniqueness of the system's solution.

Q.E.D.

\vphantom{}

Next we compute explicit formulae for $M_{q}$ and $\chi_{q}$. Aside
from being satisfying in their own right, they let us make the $q$-adic
estimates needed to prove the convergence of $\chi_{q}$'s extension
from a function $\mathbb{N}_{0}\rightarrow\mathbb{Q}$ to a function
$\mathbb{Z}_{2}\rightarrow\mathbb{Z}_{q}$.
\begin{prop}
\label{prop:Explicit Formulas for M_H}$M_{q}\left(n\right)$ and
$M_{q}\left(\mathbf{j}\right)$ can be explicitly given in terms of
the constants associated to $H$ by the formulae: 
\begin{equation}
M_{q}\left(n\right)=\frac{q^{\#_{1}\left(n\right)}}{2^{\lambda_{2}\left(n\right)}},\textrm{ }\forall n\in\mathbb{N}_{0}\label{eq:Formula for M_H of n}
\end{equation}
and: 
\begin{equation}
M_{q}\left(\mathbf{j}\right)=\frac{q^{\#_{1}\left(\mathbf{j}\right)}}{2^{\left|\mathbf{j}\right|}},\textrm{ }\forall\mathbf{j}\in\textrm{String}\left(2\right)\label{eq:formula for M_H of bold-j}
\end{equation}
respectively.
\end{prop}
Proof: We omit the proof, because it is part of the proof of the next
proposition (\textbf{Proposition \ref{prop:Explicit formula for Chi_H of bold j}}).

Q.E.D.
\begin{rem}
The main thing to note is that (\ref{eq:formula for M_H of bold-j})
shows $M_{q}\left(\mathbf{j}\right)$ will decay $q$-adically to
$0$ as the number of $1$s in $\mathbf{j}$ tends to $\infty$:
\begin{equation}
\left|M_{q}\left(\mathbf{j}\right)\right|_{q}=q^{-\#_{1}\left(\mathbf{j}\right)}\label{eq:M_H decay estimate}
\end{equation}
\end{rem}
\begin{prop}
\label{prop:Explicit formula for Chi_H of bold j} 
\begin{equation}
\chi_{q}\left(\mathbf{j}\right)=\sum_{m=1}^{\left|\mathbf{j}\right|}\frac{b_{j_{m}}}{2}\prod_{k=1}^{m-1}\frac{a_{j_{k}}}{2},\textrm{ }\forall\mathbf{j}\in\textrm{String}\left(2\right)\label{eq:Formula for Chi_H in terms of bold-j}
\end{equation}
where the $k$-product is defined to be $1$ when $m=1$. (Recall
that $a_{0}=1$, $a_{1}=q$, $b_{0}=0$, and $b_{1}=1$.)
\end{prop}
Proof: Let $\mathbf{j}=\left(j_{1},\ldots,j_{\left|\mathbf{j}\right|}\right)\in\textrm{String}\left(2\right)$
be arbitrary. Since $\chi_{q}\left(\mathbf{j}\right)=\chi_{q}\left(\mathbf{i}\right)$
for any $\mathbf{i}\in\textrm{String}\left(2\right)$ for which $\mathbf{i}\sim\mathbf{j}$,
we can assume without loss of generality that $\mathbf{j}$'s right-most
entry is non-zero; that is: $j_{\left|\mathbf{j}\right|}=1$. The
proof follows by examining \textbf{Proposition \ref{prop:H_boldj in terms of M_H and Chi_H}}'s
formula: 
\begin{equation}
H_{\mathbf{j}}\left(x\right)=M_{q}\left(\mathbf{j}\right)x+\chi_{q}\left(\mathbf{j}\right),\textrm{ }\forall x\in\mathbb{R}
\end{equation}
carefully writing out the composition sequence in full, and multiplying
it all out: 
\begin{align*}
H_{\mathbf{j}}\left(x\right) & =H_{j_{1},\ldots,j_{\left|\mathbf{j}\right|}}\left(x\right)\\
 & =\frac{a_{j_{1}}\frac{a_{j_{2}}\left(\cdots\right)+b_{j_{2}}}{2}+b_{j_{1}}}{2}\\
 & =\overbrace{\left(\prod_{k=1}^{\left|\mathbf{j}\right|}\frac{a_{j_{k}}}{2}\right)}^{M_{q}\left(\mathbf{j}\right)}x+\overbrace{\frac{b_{j_{1}}}{2}+\frac{a_{j_{1}}}{2}\frac{b_{j_{2}}}{2}+\frac{a_{j_{2}}}{2}\frac{a_{j_{1}}}{2}\frac{b_{j_{3}}}{2}+\cdots+\left(\prod_{k=1}^{\left|\mathbf{j}\right|-1}\frac{a_{j_{k}}}{2}\right)\frac{b_{j_{\left|\mathbf{j}\right|}}}{2}}^{\chi_{q}\left(\mathbf{j}\right)}\\
 & =M_{q}\left(\mathbf{j}\right)x+\underbrace{\sum_{m=1}^{\left|\mathbf{j}\right|}\frac{b_{j_{m}}}{2}\prod_{k=1}^{m-1}\frac{a_{j_{k}}}{2}}_{\chi_{q}\left(\mathbf{j}\right)}
\end{align*}
where we have used the convention that $\prod_{k=1}^{m-1}\frac{a_{j_{k}}}{2}=1$
when $m=1$. This proves:
\begin{equation}
\chi_{q}\left(\mathbf{j}\right)=\sum_{m=1}^{\left|\mathbf{j}\right|}\frac{b_{j_{m}}}{2}\prod_{k=1}^{m-1}\frac{a_{j_{k}}}{2}
\end{equation}

Q.E.D.

\vphantom{}

Now, the main result of this subsection: the $\left(2,q\right)$-adic
continuation of $\chi_{q}$, and its characterization via functional
equations.
\begin{lem}[\textbf{$\left(2,q\right)$-adic Characterization of }$\chi_{q}$]
\label{lem:Unique rising continuation and p-adic functional equation of Chi_H}The
limit\footnote{Here, the $\mathbb{Z}_{q}$ over the equality indicates that the convergence
is occurring in the topology of $\mathbb{Z}_{q}$.}: 
\begin{equation}
\chi_{q}\left(\mathfrak{z}\right)\overset{\mathbb{Z}_{q}}{=}\lim_{n\rightarrow\infty}\chi_{q}\left(\left[\mathfrak{z}\right]_{2^{n}}\right)\label{eq:Rising Continuity Formula for Chi_H}
\end{equation}
exists for all $\mathfrak{z}\in\mathbb{Z}_{2}$, and thereby defines
an extension of $\chi_{q}$ to a function $\chi_{q}:\mathbb{Z}_{2}\rightarrow\mathbb{Z}_{q}$,
with $\chi_{q}$ being continuous on $\mathbb{Z}_{2}^{\prime}$. Moreover:

\vphantom{}

I. For all $\mathfrak{z}\in\mathbb{Z}_{2}$:

\begin{align}
\chi_{q}\left(2\mathfrak{z}\right) & =\frac{\chi_{q}\left(\mathfrak{z}\right)}{2}\label{eq:Functional Equations for Chi_H over the rho-adics}\\
\chi_{q}\left(2\mathfrak{z}+1\right) & =\frac{q\chi_{q}\left(\mathfrak{z}\right)+1}{2}
\end{align}

\vphantom{}

II. The extension $\chi_{q}:\mathbb{Z}_{2}\rightarrow\mathbb{Z}_{q}$
defined by the limit \emph{(\ref{eq:Rising Continuity Formula for Chi_H})}
is the \textbf{unique} function $f:\mathbb{Z}_{2}\rightarrow\mathbb{Z}_{q}$
satisfying the functional equations:

\emph{
\begin{align}
f\left(2\mathfrak{z}\right) & =\frac{f\left(\mathfrak{z}\right)}{2}\label{eq:unique p,q-adic functional equation of Chi_H}\\
f\left(2\mathfrak{z}+1\right) & =\frac{qf\left(\mathfrak{z}\right)+1}{2}
\end{align}
}subject to the \textbf{rising-continuity condition}:
\begin{equation}
f\left(\mathfrak{z}\right)\overset{\mathbb{Z}_{q}}{=}\lim_{n\rightarrow\infty}f\left(\left[\mathfrak{z}\right]_{2^{n}}\right),\textrm{ }\forall\mathfrak{z}\in\mathbb{Z}_{2}\label{eq:unique p,q-adic rising continuity of Chi_H}
\end{equation}
\end{lem}
Proof: First, we show the existence of the limit (\ref{eq:Rising Continuity Formula for Chi_H}).
If $\mathfrak{z}\in\mathbb{N}_{0}\cap\mathbb{Z}_{2}$, then $\left[\mathfrak{z}\right]_{2^{N}}=\mathfrak{z}$
for all sufficiently large $N$, and the limit (\ref{eq:Rising Continuity Formula for Chi_H})
exists, and we are done.

So, suppose $\mathfrak{z}\in\mathbb{Z}_{2}^{\prime}=\mathbb{Z}_{2}\backslash\mathbb{N}_{0}$.
Next, let $\mathbf{j}_{n}=\left(j_{1},\ldots,j_{n}\right)$ be the
shortest string representing $\left[\mathfrak{z}\right]_{2^{n}}$;
that is: 
\begin{equation}
\left[\mathfrak{z}\right]_{2^{n}}=\sum_{k=1}^{n}j_{k}2^{k-1}
\end{equation}
By \textbf{Lemma \ref{eq:Definition of Chi_H of n}}, the equation
(\ref{eq:Rising Continuity Formula for Chi_H}) can be written as:
\begin{equation}
\chi_{q}\left(\mathfrak{z}\right)=\lim_{n\rightarrow\infty}\chi_{q}\left(\mathbf{j}_{n}\right)
\end{equation}
Using \textbf{Proposition \ref{prop:Explicit formula for Chi_H of bold j}},
we get:
\begin{equation}
\chi_{q}\left(\mathbf{j}_{n}\right)=H_{\mathbf{j}_{n}}\left(0\right)=\sum_{m=1}^{\left|\mathbf{j}_{n}\right|}\frac{b_{j_{m}}}{2}\prod_{k=1}^{m-1}\frac{a_{j_{k}}}{2}
\end{equation}
where the product is defined to be $1$ when $m=1$. Then, with \textbf{Proposition
\ref{prop:Explicit Formulas for M_H}}, we can write:
\begin{equation}
\prod_{k=1}^{m-1}\frac{a_{j_{k}}}{2}=M_{q}\left(\mathbf{j}_{m-1}\right)
\end{equation}
where, again, the product is defined to be $1$ when $m=1$. With
this, the previous equation for $\chi_{q}\left(\mathbf{j}_{n}\right)$
can be written as: 
\begin{equation}
\chi_{q}\left(\mathbf{j}_{n}\right)=\sum_{m=1}^{\left|\mathbf{j}_{n}\right|}\frac{b_{j_{m}}}{2}M_{q}\left(\mathbf{j}_{m-1}\right)
\end{equation}
Taking limits gives us: 
\begin{equation}
\lim_{n\rightarrow\infty}\chi_{q}\left(\left[\mathfrak{z}\right]_{2^{n}}\right)=\lim_{n\rightarrow\infty}\chi_{q}\left(\mathbf{j}_{n}\right)=\lim_{n\rightarrow\infty}\sum_{m=1}^{\left|\mathbf{j}_{n}\right|}\frac{b_{j_{m}}}{2}M_{q}\left(\mathbf{j}_{m-1}\right)\label{eq:Formal rising limit of Chi_H}
\end{equation}

Now, since the $b_{j}$s are either $0$ or $1$, using (\ref{eq:M_H decay estimate}),
we have that the $q$-adic absolute value of the $m$th term of the
series in (\ref{eq:Formal rising limit of Chi_H}) is dominated by
$q^{-\#_{1}\left(\mathbf{j}_{m-1}\right)}$:
\begin{equation}
\left|\frac{b_{j_{m}}}{2}M_{q}\left(\mathbf{j}_{m-1}\right)\right|_{q}\leq\left|M_{q}\left(\mathbf{j}_{m-1}\right)\right|_{q}=q^{-\#_{1}\left(\mathbf{j}_{m-1}\right)}
\end{equation}
Thanks to the ultrametric topology of $\mathbb{Z}_{q}$, the series
in (\ref{eq:Formal rising limit of Chi_H}) converges in $\mathbb{Z}_{q}$\emph{
if and only if} its $m$th term tends to $0$ $q$-adically. Because
we assumed $\mathfrak{z}$ was in $\mathbb{Z}_{2}^{\prime}$, $\mathfrak{z}$
necessarily has infinitely many non-zero $2$-adic digits\textemdash in
particular, infinitely many $1$s\textemdash which forces $\#_{1}\left(\mathbf{j}_{m-1}\right)\rightarrow\infty$
as $m\rightarrow\infty$. This proves the point-wise $q$-adic convergence
of (\ref{eq:Formal rising limit of Chi_H}) for all $\mathfrak{z}\in\mathbb{Z}_{2}^{\prime}$.

\vphantom{}

To prove (I), by \textbf{Lemma \ref{lem:Chi_H functional equation on N_0 and uniqueness}},
we know that $\chi_{q}$ satisfies the functional equations (\ref{eq:Functional Equations for Chi_H over the rho-adics})
for $\mathfrak{z}\in\mathbb{N}_{0}$. Taking limits of the functional
equations using (\ref{eq:Rising Continuity Formula for Chi_H}) then
proves the equations hold for all $\mathfrak{z}\in\mathbb{Z}_{2}$.

\vphantom{}

For (II), we have just shown that $\chi_{q}$ satisfies the rising-continuity
condition and the functional equations. So, conversely, suppose $f:\mathbb{Z}_{2}\rightarrow\mathbb{Z}_{q}$
satisfies (\ref{eq:unique p,q-adic functional equation of Chi_H})
and (\ref{eq:unique p,q-adic rising continuity of Chi_H}). Then,
by \textbf{Lemma \ref{lem:Chi_H functional equation on N_0 and uniqueness}},
the restriction of $f$ to $\mathbb{N}_{0}$ must be equal to $\chi_{q}$.
Equation (\ref{eq:Rising Continuity Formula for Chi_H}) then shows
that the values $f\left(\mathfrak{z}\right)\overset{\mathbb{Z}_{q}}{=}\lim_{n\rightarrow\infty}f\left(\left[\mathfrak{z}\right]_{2^{n}}\right)$
obtained by taking limits necessarily force $f\left(\mathfrak{z}\right)=\chi_{q}\left(\mathfrak{z}\right)$
for all $\mathfrak{z}\in\mathbb{Z}_{2}$, thus proving the uniqueness
of the rising-continuous solutions of (\ref{eq:unique p,q-adic functional equation of Chi_H}).

Q.E.D.

\vphantom{}
\begin{rem}
In the limiting process $\lim_{n\rightarrow\infty}\chi_{q}\left(\left[\mathfrak{z}\right]_{2^{n}}\right)$
used in our proof of \textbf{Lemma \ref{lem:Unique rising continuation and p-adic functional equation of Chi_H}},
we get a dichotomy based on whether $\mathfrak{z}$ is an element
of $\mathbb{N}_{0}$ or an element of $\mathbb{Z}_{2}^{\prime}$.
For $\mathfrak{z}\in\mathbb{N}_{0}$, the limit $\lim_{n\rightarrow\infty}\chi_{q}\left(\left[\mathfrak{z}\right]_{2^{n}}\right)$
is, in a sense, degenerate, because the sequence $\left\{ \chi_{q}\left(\left[\mathfrak{z}\right]_{2^{n}}\right)\right\} _{n\geq0}$
converges after only finitely many terms. This is important, because
it means the convergence in this case \emph{does not depend on the
topology used to define} $\lim_{n\rightarrow\infty}\chi_{q}\left(\left[\mathfrak{z}\right]_{2^{n}}\right)$;
equivalently, it converges with respect to the \emph{discrete topology}.
On the other hand, for $\mathfrak{z}\in\mathbb{Z}_{2}^{\prime}$,
as our proof shows, we need to interpret $\lim_{n\rightarrow\infty}\chi_{q}\left(\left[\mathfrak{z}\right]_{2^{n}}\right)$
in the context of the $q$-adic topology in order for it to be well-defined
for all such $\mathfrak{z}$.
\end{rem}
\begin{rem}
By identifying of strings with $p$-adic integers, the concatenation
operation $\wedge$ can be realized as a binary operation on $\mathbb{Z}_{p}$.
For those unfamiliar with the $p$-adic numbers, it is a very instructive
exercise to use the topological definition of continuity to show that
the binary operation induced by $\wedge$ on $\mathbb{Z}_{p}$ is
not continuous with respect to the $p$-adic topology.

\newpage{}
\end{rem}

\section{\label{sec:2 The-Correspondence-Principle}The Correspondence Principle}

IN THIS SECTION, WE WRITE $H$ TO DENOTE $T_{q}$, WHERE $q$ IS AN
ARBITRARY PRIME NUMBER. WE WILL ALSO USE $T_{q}$ TO DENOTE $T_{q}$.

\vphantom{}

In this section, we establish a correspondence between the periodic
and divergent points of $H$ in $\mathbb{Z}\backslash\left\{ 0\right\} $
with the rational integer values that $\chi_{q}$ attains for $2$-adic
integer inputs. As mentioned previously, the CP for divergent points
(CPDP) will occur as a non-constructive consequence of the CP for
periodic points (CPPP), so, for conceptual understanding, we will
focus the exposition on the case of a periodic point. The main concept
to keep in mind is that for any $x\in\mathbb{Z}$ and any $L\geq1$,
there exists a unique $\mathbf{j}\in\textrm{String}\left(2\right)$
of length $L$ so that $H^{\circ L}\left(x\right)=H_{\mathbf{j}}\left(x\right)$.
This $\mathbf{j}$ is the string which encodes the movements of $x$
under $H$.

\subsection{Preparatory Work}

Obviously, the CPPP's if and only if characterization of periodic
points of $H$ gives us two directions to prove. The preparatory work
for the proof will tackle each of these one at a time. For pedagogical
purposes, in both cases, we will go over the main ideas of the proof
before beginning the associated preparatory work.

For the first direction, we start with an $x\in\mathbb{Z}\backslash\left\{ 0\right\} $
which is \emph{known} to be a periodic point of $H$, with $x$ generating
a cycle $\Omega$ of length $L$. Letting $\mathbf{j}\in\textrm{String}\left(2\right)$
be the unique length $L$ string so that $H^{\circ L}\left(x\right)=H_{\mathbf{j}}\left(x\right)$,
observe that, as a periodic point, it must be that $x=H^{\circ L}\left(x\right)=H_{\mathbf{j}}\left(x\right)$.
In particular, when we treat $H_{\mathbf{j}}$ as an affine linear
self-map of $\mathbb{Q}$, $x$ is then a fixed point of $H_{\mathbf{j}}$,
and we can apply $H_{\mathbf{j}}$ repeatedly to $x$ without changing
anything:
\begin{equation}
H_{\mathbf{j}}^{\circ m}\left(x\right)=H_{\mathbf{j}}^{\circ m-1}\left(x\right)=\cdots=H_{\mathbf{j}}\left(x\right)=x,\textrm{ }\forall m\geq1
\end{equation}
Using the formula:
\begin{equation}
H_{\mathbf{j}}\left(x\right)=M_{q}\left(\mathbf{j}\right)x+\chi_{q}\left(\mathbf{j}\right)
\end{equation}
we can express the affine-linear map $H_{\mathbf{j}}^{\circ m}$ in
terms of $M_{q}\left(\mathbf{j}\right)$, $\chi_{q}\left(\mathbf{j}\right)$,
and $m$:
\begin{equation}
x=H_{\mathbf{j}}^{\circ m}\left(x\right)=\left(M_{q}\left(\mathbf{j}\right)\right)^{m}x+\frac{1-\left(M_{q}\left(\mathbf{j}\right)\right)^{m}}{1-M_{q}\left(\mathbf{j}\right)}\chi_{q}\left(\mathbf{j}\right)\label{eq:Ready for limit}
\end{equation}
Since $x$ is a \emph{non-zero }periodic point of $H$, $\mathbf{j}$
must contain\footnote{Else, $\mathbf{j}$ is a string of $\left|\mathbf{j}\right|$ $0$s
and $x=H_{\mathbf{j}}\left(x\right)=x/2^{\left|\mathbf{j}\right|}$,
which forces $x=0$.} at least one $1$. Consequently, our decay estimate (\ref{eq:M_H decay estimate})
for the $q$-adic absolute value of $M_{q}\left(\mathbf{j}\right)$
tells us that $\left|M_{q}\left(\mathbf{j}\right)\right|_{q}^{m}\rightarrow0$
as $m\rightarrow\infty$. Letting $m\rightarrow\infty$, the limit
(\ref{eq:Ready for limit}) converges $q$-adically to:
\begin{equation}
x=\frac{\chi_{q}\left(\mathbf{j}\right)}{1-M_{q}\left(\mathbf{j}\right)}\overset{\textrm{by def}}{=}\frac{H_{\mathbf{j}}\left(0\right)}{1-H_{\mathbf{j}}^{\prime}\left(0\right)}
\end{equation}
which is our old friend (\ref{eq:Formula for X}). The trick to getting
something useful out of this is to realize that (\ref{eq:Ready for limit})
is not the \emph{only }way to express $H_{\mathbf{j}}^{\circ m}\left(x\right)$
in the form $ax+b$. Indeed:
\begin{equation}
H_{\mathbf{j}}^{\circ m}\left(x\right)=\underbrace{\left(H_{\mathbf{j}}\circ\cdots\circ H_{\mathbf{j}}\right)}_{m\textrm{ times}}\left(x\right)=H_{\mathbf{j}^{\wedge m}}\left(x\right)=M_{q}\left(\mathbf{j}^{\wedge m}\right)x+\chi_{q}\left(\mathbf{j}^{\wedge m}\right)
\end{equation}
where, recall, $\mathbf{j}^{\wedge m}$ is the concatenation of $m$
copies of $\mathbf{j}$.

Since $M_{q}$ transforms concatenation into multiplication (\textbf{Proposition
\ref{prop:M_H concatenation identity}}), we can write:
\begin{equation}
H_{\mathbf{j}}^{\circ m}\left(x\right)=M_{q}\left(\mathbf{j}^{\wedge m}\right)x+\chi_{q}\left(\mathbf{j}^{\wedge m}\right)=\left(M_{q}\left(\mathbf{j}\right)\right)^{m}x+\chi_{q}\left(\mathbf{j}^{\wedge m}\right)
\end{equation}
Combining this with (\ref{eq:Ready for limit}) yields:
\begin{equation}
\left(M_{q}\left(\mathbf{j}\right)\right)^{m}x+\chi_{q}\left(\mathbf{j}^{\wedge m}\right)=\left(M_{q}\left(\mathbf{j}\right)\right)^{m}x+\frac{1-\left(M_{q}\left(\mathbf{j}\right)\right)^{m}}{1-M_{q}\left(\mathbf{j}\right)}\chi_{q}\left(\mathbf{j}\right)
\end{equation}
and hence:
\begin{equation}
\chi_{q}\left(\mathbf{j}^{\wedge m}\right)=\frac{1-\left(M_{q}\left(\mathbf{j}\right)\right)^{m}}{1-M_{q}\left(\mathbf{j}\right)}\chi_{q}\left(\mathbf{j}\right)
\end{equation}
As we saw, the right-hand side converges $q$-adically as $m\rightarrow\infty$.
With this, we then have that if $\mathbf{j}$ is a string containing
at least one $1$:
\begin{equation}
\chi_{q}\left(\mathbf{j}^{\wedge\infty}\right)=\lim_{m\rightarrow\infty}\chi_{q}\left(\mathbf{j}^{\wedge m}\right)\overset{\mathbb{Z}_{q}}{=}\frac{\chi_{q}\left(\mathbf{j}\right)}{1-M_{q}\left(\mathbf{j}\right)}
\end{equation}
In other words, given a finite string $\mathbf{j}$, upon letting
$\mathbf{i}$ be the (necessarily infinite)\emph{ }string obtained
by concatenating infinitely many copies of $\mathbf{j}$, we have:
\begin{equation}
\chi_{q}\left(\mathbf{i}\right)=\frac{\chi_{q}\left(\mathbf{j}\right)}{1-M_{q}\left(\mathbf{j}\right)}
\end{equation}
Since infinite strings correspond to $2$-adic integers with infinitely
many digits (elements of $\mathbb{Z}_{2}^{\prime}$), the above identity
shows that $\chi_{q}$'s values for inputs in $\mathbb{Z}_{2}^{\prime}$
are related to the formula (\ref{eq:Formula for X}) for the fixed
points of $H_{\mathbf{j}}$s. This will give us the first half of
the CP.

The takeaway from this preliminary analysis tells us where to go first:
we need to study what happens when we take a given $\mathbf{j}\in\textrm{String}\left(2\right)$
and concatenate it with itself. Our first result is a simple observation
of how how concatenation of multiple copies of a single string $\mathbf{j}$
affects the integer $n$ that $\mathbf{j}$ represents.
\begin{prop}
\label{prop:Concatenation exponentiation}Let $n\in\mathbb{N}_{1}$,
and let $\mathbf{j}\in\textrm{String}\left(2\right)$ be the shortest
string representing $n$. Then, for all $m\in\mathbb{N}_{1}$: 
\begin{equation}
\mathbf{j}^{\wedge m}\sim n\frac{1-2^{m\lambda_{2}\left(n\right)}}{1-2^{\lambda_{2}\left(n\right)}}\label{eq:Proposition 1.2.10}
\end{equation}
\end{prop}
Proof: Since: 
\begin{equation}
n=j_{1}+j_{2}2+\cdots+j_{\lambda_{2}\left(n\right)-1}2^{\lambda_{2}\left(n\right)-1}
\end{equation}
we have that: 
\begin{align*}
\mathbf{j}^{\wedge m} & \sim j_{1}+j_{2}2+\cdots+j_{\lambda_{2}\left(n\right)-1}2^{\lambda_{2}\left(n\right)-1}\\
 & +2^{\lambda_{2}\left(n\right)}\left(j_{1}+j_{2}2+\cdots+j_{\lambda_{2}\left(n\right)-1}2^{\lambda_{2}\left(n\right)-1}\right)\\
 & +\cdots\\
 & +2^{\left(m-1\right)\lambda_{2}\left(n\right)}\left(j_{1}+j_{2}2+\cdots+j_{\lambda_{2}\left(n\right)-1}2^{\lambda_{2}\left(n\right)-1}\right)\\
 & =n+n2^{\lambda_{2}\left(n\right)}+n2^{2\lambda_{2}\left(n\right)}+\cdots+n2^{\left(m-1\right)\lambda_{2}\left(n\right)}\\
 & =n\frac{1-2^{m\lambda_{2}\left(n\right)}}{1-2^{\lambda_{2}\left(n\right)}}
\end{align*}
as desired.

Q.E.D.

\vphantom{}

As we saw in the preliminary analysis, the key step was to consider
the behavior of $\chi_{q}\left(\mathbf{j}^{\wedge m}\right)$ as $m\rightarrow\infty$.
Equivalently, we are studying what happens when we pre-compose $\chi_{q}$
with the map that accepts a string $\mathbf{j}\in\textrm{String}\left(2\right)$
and outputs the string $\mathbf{j}^{\wedge\infty}\in\textrm{String}_{\infty}\left(2\right)$.
We will call this map $B_{2}$.
\begin{defn}[$B_{2}$]
We define $B_{2}:\mathbb{N}_{0}\rightarrow\mathbb{Z}_{2}$ by: 
\begin{equation}
B_{2}\left(n\right)\overset{\textrm{def}}{=}\begin{cases}
0 & \textrm{if }n=0\\
\frac{n}{1-2^{\lambda_{2}\left(n\right)}} & \textrm{if }n\geq1
\end{cases}\label{eq:Definition of B projection function}
\end{equation}
\end{defn}
\begin{rem}
In terms of $2$-adic expansions, $B_{2}$ sends $n$ to the $2$-adic
integer $B_{2}\left(n\right)$ whose sequence of $2$-adic digits
consists of infinitely many concatenated copies of the sequence of
$2$-adic digits of $n$: 
\begin{equation}
B_{2}\left(n\right)\overset{\mathbb{Z}_{2}}{=}\lim_{m\rightarrow\infty}n\frac{1-2^{m\lambda_{2}\left(n\right)}}{1-2^{\lambda_{2}\left(n\right)}}\overset{\mathbb{Z}_{2}}{=}\frac{n}{1-2^{\lambda_{2}\left(n\right)}}
\end{equation}
In particular, since the sequence of $2$-adic digits of $B_{2}\left(n\right)$
is, by construction, periodic, the geometric series formula in $\mathbb{Z}_{2}$
guarantees that the $2$-adic integer $B_{2}\left(n\right)$ is an
element of $\mathbb{Q}$.

\end{rem}
\begin{rem}
\label{rem:Self-concatenation as a rising limit}Let $\mathbf{j}\in\textrm{String}\left(2\right)$
be non-zero and let $n$ be the (necessarily positive) integer represented
by $\mathbf{j}$. Because the sequence of the $2$-adic digits of
$B_{2}\left(n\right)$ is the concatenation of infinitely many copies
of $\mathbf{j}$, observe that for each integer $k\geq1$, the string
$\mathbf{j}^{\wedge k}$ represents $\left[B_{2}\left(n\right)\right]_{2^{k\left|\mathbf{j}\right|}}$
(the value of the $2$-adic integer $B_{2}\left(n\right)$ modulo
$2^{k\left|\mathbf{j}\right|}$). As such, for any $f:\mathbb{Z}_{2}\rightarrow\mathbb{Z}_{q}$
satisfying the rising-continuity condition (\ref{eq:unique p,q-adic rising continuity of Chi_H}),
we have that:
\begin{equation}
\lim_{k\rightarrow\infty}f\left(\mathbf{j}^{\wedge k}\right)\overset{\mathbb{Z}_{q}}{=}f\left(B_{2}\left(n\right)\right)
\end{equation}
where, by $f\left(\mathbf{j}^{\wedge k}\right)$, we mean the value
attained by $f$ at the integer ($\left[B_{2}\left(n\right)\right]_{2^{k\left|\mathbf{j}\right|}}$)
represented by $\mathbf{j}^{\wedge k}$.
\end{rem}
\begin{rem}
Note that we can extend $B_{2}$ to a function on $\mathbb{Z}_{2}$
by defining its restriction to $\mathbb{Z}_{2}^{\prime}$ as being
the identity map. In doing so, note that $B_{2}:\mathbb{Z}_{2}\rightarrow\mathbb{Z}_{2}$
becomes an idempotent map.
\end{rem}
With $B_{2}$ at our disposal, the behavior of $\chi_{q}\left(\mathbf{j}^{\wedge m}\right)$
as $m\rightarrow\infty$ can be realized as an instance of a functional
equation for $\chi_{q}\circ B_{2}$.
\begin{lem}[\textbf{Functional Equation for }$\chi_{q}\circ B_{2}$]
\label{lem:Chi_H o B_p functional equation}\ 
\begin{equation}
\chi_{q}\left(B_{2}\left(n\right)\right)\overset{\mathbb{Z}_{q}}{=}\frac{\chi_{q}\left(n\right)}{1-M_{q}\left(n\right)},\textrm{ }\forall n\in\mathbb{N}_{1}\label{eq:Chi_H B functional equation}
\end{equation}
\end{lem}
Proof: Let $n\in\mathbb{N}_{1}$, and let $\mathbf{j}\in\textrm{String}\left(2\right)$
be the shortest string representing $n$. Note that $\mathbf{j}$
contains at least one $1$. Now, for all $m\geq1$ and all $k\in\mathbb{N}_{1}$:
\begin{equation}
H_{\mathbf{j}^{\wedge m}}\left(k\right)=H_{\mathbf{j}}^{\circ m}\left(k\right)
\end{equation}
Given an affine linear map of the form $ax+b$, a simple application
of the geometric series shows that function obtained by iterating
$ax+b$ $m$ times is:
\begin{equation}
a^{m}x+\frac{1-a^{m}}{1-a}b
\end{equation}
Since $H_{\mathbf{j}}$ is an affine linear map:
\begin{equation}
H_{\mathbf{j}}\left(k\right)=M_{q}\left(\mathbf{j}\right)k+\chi_{q}\left(\mathbf{j}\right)
\end{equation}
this gives:
\begin{equation}
H_{\mathbf{j}^{\wedge m}}\left(k\right)=H_{\mathbf{j}}^{\circ m}\left(k\right)=\left(M_{q}\left(\mathbf{j}\right)\right)^{m}k+\frac{1-\left(M_{q}\left(\mathbf{j}\right)\right)^{m}}{1-M_{q}\left(\mathbf{j}\right)}\chi_{q}\left(\mathbf{j}\right)
\end{equation}
Noting also that: 
\begin{equation}
H_{\mathbf{j}^{\wedge m}}\left(k\right)=M_{q}\left(\mathbf{j}^{\wedge m}\right)k+\chi_{q}\left(\mathbf{j}^{\wedge m}\right)
\end{equation}
we equate our two formulae for $H_{\mathbf{j}^{\wedge m}}\left(k\right)$
to obtain:
\begin{equation}
M_{q}\left(\mathbf{j}^{\wedge m}\right)k+\chi_{q}\left(\mathbf{j}^{\wedge m}\right)=\left(M_{q}\left(\mathbf{j}\right)\right)^{m}k+\frac{1-\left(M_{q}\left(\mathbf{j}\right)\right)^{m}}{1-M_{q}\left(\mathbf{j}\right)}\chi_{q}\left(\mathbf{j}\right)
\end{equation}
By \textbf{Proposition \ref{prop:M_H concatenation identity}}, we
see that $M_{q}\left(\mathbf{j}^{\wedge m}\right)k=\left(M_{q}\left(\mathbf{j}\right)\right)^{m}k$.
Cancelling these terms from both sides leaves us with: 
\begin{equation}
\chi_{q}\left(\mathbf{j}^{\wedge m}\right)=\frac{1-\left(M_{q}\left(\mathbf{j}\right)\right)^{m}}{1-M_{q}\left(\mathbf{j}\right)}\chi_{q}\left(\mathbf{j}\right)=\frac{1-\left(M_{q}\left(\mathbf{j}\right)\right)^{m}}{1-M_{q}\left(n\right)}\chi_{q}\left(n\right)\label{eq:Chi_H B functional equation, ready to take limits}
\end{equation}

By \textbf{Proposition \ref{prop:Concatenation exponentiation}},
we have: 
\begin{equation}
\chi_{q}\left(\mathbf{j}^{\wedge m}\right)\overset{\mathbb{Q}}{=}\chi_{q}\left(n\frac{1-2^{m\lambda_{2}\left(n\right)}}{1-2^{\lambda_{2}\left(n\right)}}\right)
\end{equation}
where, as indicated, the equality occurs in $\mathbb{Q}$. Since $n\in\mathbb{N}_{1}$,
we have that: 
\begin{equation}
B_{2}\left(n\right)=\frac{n}{1-2^{\lambda_{2}\left(n\right)}}
\end{equation}
is a $2$-adic integer. Moreover, as is immediate from the proof of
\textbf{Proposition \ref{prop:Concatenation exponentiation}}, the
projection of this $2$-adic integer modulo $2^{m}$ is: 
\begin{equation}
\left[B_{2}\left(n\right)\right]_{2^{m}}=n\frac{1-2^{m\lambda_{2}\left(n\right)}}{1-2^{\lambda_{2}\left(n\right)}}
\end{equation}
which is, of course, exactly the rational integer represented by the
string $\mathbf{j}^{\wedge m}$. In other words: 
\begin{equation}
\frac{1-\left(M_{q}\left(\mathbf{j}\right)\right)^{m}}{1-M_{q}\left(n\right)}\chi_{q}\left(n\right)=\chi_{q}\left(\mathbf{j}^{\wedge m}\right)=\chi_{q}\left(n\frac{1-2^{m\lambda_{2}\left(n\right)}}{1-2^{\lambda_{2}\left(n\right)}}\right)=\chi_{q}\left(\left[B_{2}\left(n\right)\right]_{2^{m}}\right)
\end{equation}
Applying \textbf{Lemma \ref{lem:Unique rising continuation and p-adic functional equation of Chi_H}}
yields:
\begin{equation}
\chi_{q}\left(B_{2}\left(n\right)\right)\overset{\mathbb{Z}_{q}}{=}\lim_{m\rightarrow\infty}\chi_{q}\left(\left[B_{2}\left(n\right)\right]_{2^{m}}\right)\overset{\mathbb{Z}_{q}}{=}\lim_{m\rightarrow\infty}\frac{1-\left(M_{q}\left(\mathbf{j}\right)\right)^{m}}{1-M_{q}\left(n\right)}\chi_{q}\left(n\right)
\end{equation}
where, as indicated, the limits occur in $\mathbb{Z}_{q}$.

Finally, since $\mathbf{j}$ has at least one $1$, \textbf{Proposition
\ref{prop:Explicit Formulas for M_H}} shows that $\left|M_{q}\left(\mathbf{j}\right)\right|_{q}<1$,
and hence, that $\left(M_{q}\left(\mathbf{j}\right)\right)^{m}$ tends
to $0$ in $\mathbb{Z}_{q}$ as $m\rightarrow\infty$. Thus: 
\begin{equation}
\chi_{q}\left(B_{2}\left(n\right)\right)\overset{\mathbb{Z}_{q}}{=}\lim_{m\rightarrow\infty}\frac{1-\left(M_{q}\left(\mathbf{j}\right)\right)^{m}}{1-M_{q}\left(n\right)}\chi_{q}\left(n\right)\overset{\mathbb{Z}_{q}}{=}\frac{\chi_{q}\left(n\right)}{1-M_{q}\left(n\right)}
\end{equation}

Q.E.D.
\begin{rem}
The functional equation (\ref{eq:Chi_H B functional equation}) is
a vast generalization of the \textbf{Böhm-Sontacchi Criterion} (\ref{eq:The Bohm-Sontacchi Criterion}).
The fraction on the right-hand side of (\ref{eq:Chi_H B functional equation})
is the right-hand side of (\ref{eq:The Bohm-Sontacchi Criterion})
(albeit for $T_{q}$ instead of $T_{3}$) parameterized in terms of
the positive integer variable $n$.
\end{rem}
With that done, we now prepare for the second direction of the CPPP.

This time, instead of starting with a periodic point $x\in\mathbb{Z}\backslash\left\{ 0\right\} $,
we begin with a string $\mathbf{j}\in\textrm{String}\left(2\right)$
of length $L\geq1$. Since $\left|\mathbf{j}\right|\geq1$, $H_{\mathbf{j}}$
is an affine linear map of the form:
\begin{equation}
H_{\mathbf{j}}\left(x\right)=M_{H}\left(\mathbf{j}\right)x+\chi_{H}\left(\mathbf{j}\right)
\end{equation}
where $M_{H}\left(\mathbf{j}\right)\notin\left\{ 0,1\right\} $. As
a result, $H_{\mathbf{j}}$ has a unique fixed point $x_{0}\in\mathbb{Q}$,
given by:
\begin{equation}
x_{0}=\frac{H_{\mathbf{j}}\left(0\right)}{1-H_{\mathbf{j}}^{\prime}\left(0\right)}=\frac{\chi_{q}\left(\mathbf{j}\right)}{1-M_{q}\left(\mathbf{j}\right)}\label{eq:Thing}
\end{equation}
Obviously, if $x_{0}$ is not an integer, it cannot be a fixed point
of $H$. Thus, what we want is to conclude that $x_{0}$ being an
integer forces it to also be a fixed point of $H$ itself. To arrive
at this conclusion, we need to show that fixed points of $H_{\mathbf{j}}$
in $\mathbb{Z}$ are fixed points of $H$ in $\mathbb{Z}$. As we
shall see, this boils down to the notion of the ``correctness''
of $\mathbf{j}$.

For any given integer $x$, we once again observe that there is a
unique string $\mathbf{i}$ of length $L$ so that $H_{\mathbf{i}}\left(x\right)=H^{\circ L}\left(x\right)$.
This then leads to a notion of ``correctness'': a string \textbf{\emph{$\mathbf{i}$}}
will be ``correct'' for a given $x$ if $H_{\mathbf{i}}\left(x\right)=H^{\circ\left|\mathbf{i}\right|}\left(x\right)$,
and will be ``wrong'' for $x$ if $H_{\mathbf{i}}\left(x\right)\neq H^{\circ\left|\mathbf{i}\right|}\left(x\right)$.
In order to tackle the second half of the CPPP, we will show that
$\mathbf{i}$ is correct for $x$ precisely when $H_{\mathbf{i}}\left(x\right)\in\mathbb{Z}$;
likewise, $\mathbf{i}$ is wrong for $x$ if and only if $H_{\mathbf{i}}\left(x\right)\notin\mathbb{Z}$.
This will be taken care of by lemmata \textbf{\ref{lem:wrong values lemma}}
and \textbf{\ref{lem:properness lemma}}.

With these results in hand, we can then show that $x_{0}$ being an
element of $\mathbb{Z}$ forces $\mathbf{j}$ to be correct, and hence,
forces $x_{0}$ to satisfy $x_{0}=H_{\mathbf{j}}\left(x_{0}\right)=H^{\left|\mathbf{j}\right|}\left(x_{0}\right)$,
which proves that $x_{0}$ is indeed a periodic point.

As it turns out, the ``correctness'' issue described above is most
easily done by considering the extension of $H$ (that is, $T_{q}$)
to $\mathbb{Z}_{2}$. This is a well-established procedure in Collatz
studies \cite{Lagarias-Kontorovich Paper}, with the extension $H:\mathbb{Z}_{2}\rightarrow\mathbb{Z}_{2}$
being defined by:
\begin{equation}
H\left(\mathfrak{z}\right)\overset{\textrm{def}}{=}T_{q}\left(\mathfrak{z}\right)\overset{\textrm{def}}{=}\begin{cases}
\frac{\mathfrak{z}}{2} & \textrm{if }\mathfrak{z}\overset{2}{\equiv}0\\
\frac{q\mathfrak{z}+1}{2} & \textrm{if }\mathfrak{z}\overset{2}{\equiv}1
\end{cases}\label{eq:Definition of T_3 on Z_2}
\end{equation}
Note that $H$ is then continuous on $\mathbb{Z}_{2}$. Moreover,
the branches $H_{0}$ and $H_{1}$ are then continuous functions on
$\mathbb{Z}_{2}$, and in fact, on $\mathbb{Q}_{2}$ as well, and
the same applies to $H_{\mathbf{j}}$ for any $\mathbf{j}\in\textrm{String}\left(2\right)$.
With this in mind, we can approach the ``correctness'' question
by defining what it means for a given $\mathbf{j}$ to be \emph{``}wrong''.
\begin{defn}[\textbf{Wrong values \& seeds}]
We say a $2$-adic number $\mathfrak{y}\in\mathbb{Q}_{2}$ is a \textbf{wrong
value for $H$ }whenever there is a $\mathbf{j}\in\textrm{String}\left(2\right)$
and a $\mathfrak{z}\in\mathbb{Z}_{2}$ so that $\mathfrak{y}=H_{\mathbf{j}}\left(\mathfrak{z}\right)$
and $H_{\mathbf{j}}\left(\mathfrak{z}\right)\neq H^{\circ\left|\mathbf{j}\right|}\left(\mathfrak{z}\right)$.
We then call $\mathfrak{z}$ a \textbf{seed }of $\mathfrak{y}$.
\end{defn}
\begin{rem}
If $\mathfrak{y}$ is a wrong value of $H$\textemdash with seed $\mathfrak{z}$
and string $\mathbf{j}$, so that $\mathfrak{y}=H_{\mathbf{j}}\left(\mathfrak{z}\right)$\textemdash note
that for any $\mathbf{i}\in\textrm{String}\left(2\right)$, the $2$-adic
number:
\begin{equation}
H_{\mathbf{i}}\left(\mathfrak{y}\right)=H_{\mathbf{i}}\left(H_{\mathbf{j}}\left(\mathfrak{z}\right)\right)=H_{\mathbf{i}\wedge\mathbf{j}}\left(\mathfrak{z}\right)
\end{equation}
will \emph{also }be a wrong value of $H$ with seed $\mathfrak{z}$.
Thus, branches of $H$ \emph{always map wrong values to wrong values}.
Expressing this in terms of $H_{\mathbf{j}}$'s effect on $\mathbb{Q}_{2}/\mathbb{Z}_{2}$
(the set of all $2$-adic numbers with $2$-adic absolute values $>1$)
will then give us the results that we need.
\end{rem}
\begin{prop}
\label{prop:Q_p / Z_p prop}$H_{j}\left(\mathbb{Q}_{2}\backslash\mathbb{Z}_{2}\right)\subseteq\mathbb{Q}_{2}\backslash\mathbb{Z}_{2}$
for all $j\in\mathbb{Z}/2\mathbb{Z}$.
\end{prop}
Proof: Note that an element of $\mathbb{Q}_{2}$ is in $\mathbb{Q}_{2}\backslash\mathbb{Z}_{2}$
if and only if it has $2$-adic absolute value $>1$. So, let $\mathfrak{y}\in\mathbb{Q}_{2}\backslash\mathbb{Z}_{2}$.
Then:
\begin{equation}
\left|H_{0}\left(\mathfrak{y}\right)\right|_{2}=\left|\frac{\mathfrak{y}}{2}\right|_{2}=2\left|\mathfrak{y}\right|_{2}>2
\end{equation}
so, $H_{0}\left(\mathfrak{y}\right)\in\mathbb{Q}_{2}/\mathbb{Z}_{2}$.
On the other hand:
\begin{equation}
\left|H_{1}\left(\mathfrak{y}\right)\right|_{2}=\left|\frac{q\mathfrak{y}+1}{2}\right|_{2}=2\left|q\mathfrak{y}+1\right|_{2}
\end{equation}
Since $\left|q\mathfrak{y}\right|_{2}>1=\left|1\right|_{2}$, the
$2$-adic ultrametric inequality tells us that:
\begin{equation}
\left|q\mathfrak{y}+1\right|_{2}=\left|q\mathfrak{y}\right|_{2}=\left|\mathfrak{y}\right|_{2}
\end{equation}
Hence: 
\begin{equation}
\left|H_{1}\left(\mathfrak{y}\right)\right|_{2}=2\left|q\mathfrak{y}+1\right|_{2}=2\left|\mathfrak{y}\right|_{2}>2
\end{equation}
which shows that $H_{1}\left(\mathfrak{y}\right)\in\mathbb{Q}_{2}\backslash\mathbb{Z}_{2}$

Q.E.D.

\vphantom{}

The following two lemmata will take care of the ``correctness''
issue.
\begin{lem}
\label{lem:wrong values lemma}All wrong values of $H$ are elements
of $\mathbb{Q}_{2}\backslash\mathbb{Z}_{2}$.
\end{lem}
Proof: Let $\mathfrak{z}\in\mathbb{Z}_{2}$, and let $i\in\mathbb{Z}/2\mathbb{Z}$
be so that $\left[\mathfrak{z}\right]_{2}\neq i$. Then, by definition
of properness, the quantity: 
\begin{equation}
H_{i}\left(\mathfrak{z}\right)=\frac{a_{i}\mathfrak{z}+b_{i}}{2}
\end{equation}
has $2$-adic absolute value $>1$. By \textbf{Proposition \ref{prop:Q_p / Z_p prop}},
this forces $H_{\mathbf{j}}\left(H_{i}\left(\mathfrak{z}\right)\right)$
to be an element of $\mathbb{Q}_{2}\backslash\mathbb{Z}_{2}$ for
all $\mathbf{j}\in\textrm{String}\left(2\right)$. Since every wrong
value $\mathfrak{y}$ with seed $\mathfrak{z}\in\mathbb{Z}_{2}$ is
of the form $\mathfrak{y}=H_{\mathbf{j}}\left(H_{i}\left(\mathfrak{z}\right)\right)$
for some $\mathbf{j}\in\textrm{String}\left(2\right)$ and some $i\in\mathbb{Z}/2\mathbb{Z}$
for which $\left[\mathfrak{z}\right]_{2}\neq i$, this shows that
every wrong value of $H$ is in $\mathbb{Q}_{2}\backslash\mathbb{Z}_{2}$.

Q.E.D.

\vphantom{}
\begin{lem}
\label{lem:properness lemma}Let $\mathfrak{z}\in\mathbb{Z}_{2}$,
and let $\mathbf{j}\in\textrm{String}\left(2\right)$. If $H_{\mathbf{j}}\left(\mathfrak{z}\right)=\mathfrak{z}$,
then $H^{\circ\left|\mathbf{j}\right|}\left(\mathfrak{z}\right)=\mathfrak{z}$;
i.e. $\mathfrak{z}$ is a periodic point of $H$.
\end{lem}
Proof: Let $\mathfrak{z}$ and $\mathbf{j}$ be as given, with $\mathfrak{z}=H_{\mathbf{j}}\left(\mathfrak{z}\right)$.
By way of contradiction, suppose $H^{\circ\left|\mathbf{j}\right|}\left(\mathfrak{z}\right)\neq\mathfrak{z}$.
This shows us that $H^{\circ\left|\mathbf{j}\right|}\left(\mathfrak{z}\right)\neq H_{\mathbf{j}}\left(\mathfrak{z}\right)$.
Hence, $H_{\mathbf{j}}\left(\mathfrak{z}\right)$ is a wrong value
of $H$ with seed $\mathfrak{z}$. \textbf{Lemma} \textbf{\ref{lem:wrong values lemma}}
then forces $\mathfrak{z}=H_{\mathbf{j}}\left(\mathfrak{z}\right)\in\mathbb{Q}_{2}\backslash\mathbb{Z}_{2}$.
However, $H_{\mathbf{j}}\left(\mathfrak{z}\right)=\mathfrak{z}$,
and $\mathfrak{z}$ was given to be in $\mathbb{Z}_{2}$. This is
impossible!

Consequently, $H_{\mathbf{j}}\left(\mathfrak{z}\right)=\mathfrak{z}$
implies $H^{\circ\left|\mathbf{j}\right|}\left(\mathfrak{z}\right)=\mathfrak{z}$.

Q.E.D.

\vphantom{}

This completes the preparatory work.

\subsection{Proofs}

Though the summaries described in the previous subsection are, hopefully,
straight-forward, the actual business of proving the CP is a bit more
delicate. The proof of the CPPP is done in three stages, the last
of which (\textbf{Corollary \ref{cor:CP v4}}) is the version stated
in this paper's introduction.

The first version\textbf{ }we prove is \textbf{Theorem \ref{thm:CP v1}}.
The correspondence established by \textbf{Theorem \ref{thm:CP v1}
}is slightly weaker \textbf{Corollary \ref{cor:CP v4}}'s. The corollary
gives a correspondence between the values attained by $\chi_{q}$
and the periodic \emph{points }of $H$; the theorem, on the other
hand, gives a correspondence between the values attained by $\chi_{q}$
and the \emph{cycles }of $H$. \textbf{Corollary \ref{cor:CP v2}}
gives a Böhm-Sontacchi style expression for the periodic points of
$H$, while \textbf{Corollary \ref{cor:CP v4}} refines \textbf{Theorem
\ref{thm:CP v1}}'s correspondence from one of cycles to one of individual
periodic points of $\chi_{q}$.
\begin{thm}[\textbf{Correspondence Principle (Periodic Points), Ver. 1}]
\label{thm:CP v1}

\vphantom{}

I. Let $\Omega$ be any cycle of $T_{q}$ in $\mathbb{Z}$ with $\left|\Omega\right|\geq2$.
Then, there exist $x\in\Omega$ and $n\in\mathbb{N}_{1}$ so that:
\begin{equation}
\chi_{q}\left(B_{2}\left(n\right)\right)\overset{\mathbb{Z}_{q}}{=}x
\end{equation}
That is, there is an $n$ so that the infinite series defining $\chi_{q}\left(B_{2}\left(n\right)\right)$
converges\footnote{This series is obtained by evaluating $\chi_{q}\left(\left[B_{2}\left(n\right)\right]_{2^{m}}\right)$
using \textbf{Proposition \ref{prop:Explicit formula for Chi_H of bold j}},
and then taking the limit in $\mathbb{Z}_{q}$ as $m\rightarrow\infty$.} in $\mathbb{Z}_{q}$ to $x$.

\vphantom{}

II. Let $n\in\mathbb{N}_{1}$. If the rational number $\chi_{q}\left(B_{2}\left(n\right)\right)$
is in $\mathbb{Z}_{2}$, then $\chi_{q}\left(B_{2}\left(n\right)\right)$
is a periodic point of $T_{q}$ in $\mathbb{Z}_{2}$. In particular,
if $\chi_{q}\left(B_{2}\left(n\right)\right)$ is in $\mathbb{Z}$,
then $\chi_{q}\left(B_{2}\left(n\right)\right)$ is a periodic point
of $T_{q}$ in $\mathbb{Z}$.
\end{thm}
Proof: (Recall, we write $H$ to denote $T_{q}$.)

I. Let $\Omega$ be a cycle of $H$ in $\mathbb{Z}$ with $\left|\Omega\right|\geq2$.
Given any periodic point $x\in\Omega$ of $H$, there is going to
be a string $\mathbf{j}\in\textrm{String}\left(2\right)$ of length
$\left|\Omega\right|$ so that $H_{\mathbf{j}}\left(x\right)=x$.
In particular, since $\left|\Omega\right|\geq2$,\textbf{ }observe
that $\mathbf{j}$ contains \emph{at least one} non-zero entry; else,
$H_{\mathbf{j}}\left(x\right)=x/2^{\left|\mathbf{j}\right|}$\textemdash and
hence, $x$\textemdash must be $0$. A moment's thought shows that
for any $x^{\prime}\in\Omega$ other\emph{ }than $x$, the entries
of the string $\mathbf{i}$ for which $H_{\mathbf{i}}\left(x^{\prime}\right)=x^{\prime}$
must be a cyclic permutation of the entries of $\mathbf{j}$.
\begin{example}
For example, for the cycle $\left\{ 1,2\right\} $ of $T_{3}$, our
branches are: 
\begin{align*}
H_{0}\left(x\right) & =\frac{x}{2}\\
H_{1}\left(x\right) & =\frac{3x+1}{2}
\end{align*}
The composition sequence which sends $1$ back to itself first applies
$H_{1}$ to $1$, followed by $H_{0}$: 
\begin{equation}
H_{0}\left(H_{1}\left(1\right)\right)=H_{0}\left(2\right)=1
\end{equation}
On the other hand, the composition sequence which sends $2$ back
to itself first applies $H_{0}$ to $2$, followed by $H_{1}$: 
\begin{equation}
H_{1}\left(H_{0}\left(2\right)\right)=H_{1}\left(1\right)=2
\end{equation}
and the strings $\left(0,1\right)$ and $\left(1,0\right)$ are cyclic
permutations of one another.
\end{example}
In this way, note that there must exist an $x\in\Omega$ and a $\mathbf{j}$
(for which $H_{\mathbf{j}}\left(x\right)=x$) so that the\emph{ right-most
entry of $\mathbf{j}$ is non-zero}. From this point forward, we will
fix $x$ and $\mathbf{j}$ so that this condition is satisfied.

That being done, the next observation we make is that $x$ is a fixed
point of $H_{\mathbf{j}}$, and so:
\begin{equation}
x=H_{\mathbf{j}}^{\circ m}\left(x\right)=H_{\mathbf{j}^{\wedge m}}\left(x\right),\textrm{ }\forall m\geq0
\end{equation}
Writing out the affine linear map $H_{\mathbf{j}^{\wedge m}}$ in
$ax+b$ form gives us: 
\begin{equation}
x=H_{\mathbf{j}^{\wedge m}}\left(x\right)=M_{H}\left(\mathbf{j}^{\wedge m}\right)x+\chi_{H}\left(\mathbf{j}^{\wedge m}\right)=\left(M_{H}\left(\mathbf{j}\right)\right)^{m}x+\chi_{H}\left(\mathbf{j}^{\wedge m}\right)
\end{equation}
where the right-most equality follows from $M_{q}$'s concatenation
identity. Since $\mathbf{j}$ the right-most entry of $\mathbf{j}$
is non-zero, \textbf{Proposition \ref{prop:Explicit Formulas for M_H}}
implies that $\left|M_{q}\left(\mathbf{j}\right)\right|_{q}<1$. As
such, $\left|\left(M_{H}\left(\mathbf{j}\right)\right)^{m}\right|_{q}\rightarrow0$
as $m\rightarrow\infty$. Taking limits as $m\rightarrow\infty$ then
produces: 
\begin{equation}
x\overset{\mathbb{Z}_{q}}{=}\lim_{m\rightarrow\infty}\left(\left(M_{q}\left(\mathbf{j}\right)\right)^{m}x+\chi_{q}\left(\mathbf{j}^{\wedge m}\right)\right)\overset{\mathbb{Z}_{q}}{=}\lim_{m\rightarrow\infty}\chi_{q}\left(\mathbf{j}^{\wedge m}\right)
\end{equation}

Now, let $n$ be the integer represented by $\mathbf{j}$. Since $\mathbf{j}$'s
right-most entry is non-zero, $\mathbf{j}$ is necessarily the \emph{shortest}
string representing $n$, and $n$ is non-zero, with $\lambda_{2}\left(n\right)=\left|\mathbf{j}\right|>0$.
Using \textbf{Proposition \ref{prop:Concatenation exponentiation}},\textbf{
}we can then write: 
\begin{equation}
\mathbf{j}^{\wedge m}\sim n\frac{1-2^{m\lambda_{2}\left(n\right)}}{1-2^{\lambda_{2}\left(n\right)}},\textrm{ }\forall m\geq1
\end{equation}
As such, just as in the proof of \textbf{Lemma \ref{lem:Chi_H o B_p functional equation}},
we have: 
\begin{align}
x & \overset{\mathbb{Z}_{q}}{=}\lim_{m\rightarrow\infty}\chi_{q}\left(\mathbf{j}^{\wedge m}\right)\\
 & \overset{\mathbb{Z}_{q}}{=}\lim_{m\rightarrow\infty}\chi_{q}\left(n\frac{1-2^{m\lambda_{2}\left(n\right)}}{1-2^{\lambda_{2}\left(n\right)}}\right)\\
 & \overset{\mathbb{Z}_{q}}{=}\chi_{q}\left(\frac{n}{1-2^{\lambda_{2}\left(n\right)}}\right)\\
 & =\chi_{q}\left(B_{2}\left(n\right)\right)
\end{align}
This proves the existence of the desired $x$ and $n$.

\vphantom{}

II. \textbf{\emph{Before we assume the hypotheses of (II)}}, let us
do a brief computation. This will be the key step of the proof of
(II), and involves exploiting the fact that the branches of $H$ are
continuous self-maps of the \emph{$q$-adic integers}.
\begin{claim}
\label{claim:2.2}Let $n$ be any integer $\geq1$, and let $\mathbf{j}\in\textrm{String}\left(2\right)$
be the shortest string representing $n$. Then, $H_{\mathbf{j}}$
is a continuous map $\mathbb{Z}_{q}\rightarrow\mathbb{Z}_{q}$, and,
moreover, the $q$-adic integer $\chi_{q}\left(B_{2}\left(n\right)\right)$
is fixed by $H_{\mathbf{j}}$.

Proof of claim: Let $n$ and $\mathbf{j}$ be as given. Using $\chi_{q}$'s
concatenation identity, we can write: 
\begin{equation}
\chi_{q}\left(\mathbf{j}^{\wedge k}\right)=H_{\mathbf{j}}\left(\chi_{q}\left(\mathbf{j}^{\wedge\left(k-1\right)}\right)\right)=\cdots=H_{\mathbf{j}^{\wedge\left(k-1\right)}}\left(\chi_{q}\left(\mathbf{j}\right)\right)\label{eq:Concatenation identity for Chi_H}
\end{equation}
By \textbf{Proposition \ref{prop:Concatenation exponentiation}},
we know that $B_{2}\left(n\right)$ represents $\lim_{k\rightarrow\infty}\mathbf{j}^{\wedge k}$
(i.e. $\textrm{DigSum}_{2}\left(\lim_{k\rightarrow\infty}\mathbf{j}^{\wedge k}\right)=B_{2}\left(n\right)$).
In particular, we have that: 
\begin{equation}
\left[B_{2}\left(n\right)\right]_{2^{k}}=\mathbf{j}^{\wedge k}=n\frac{1-2^{k\lambda_{2}\left(n\right)}}{1-2^{\lambda_{2}\left(n\right)}}
\end{equation}
Letting $k\rightarrow\infty$, the limit (\ref{eq:Rising Continuity Formula for Chi_H})
tells us that: 
\begin{equation}
\lim_{k\rightarrow\infty}H_{\mathbf{j}^{\wedge\left(k-1\right)}}\left(\chi_{q}\left(n\right)\right)\overset{\mathbb{Z}_{q}}{=}\chi_{q}\left(\frac{n}{1-2^{\lambda_{2}\left(n\right)}}\right)=\chi_{q}\left(B_{2}\left(n\right)\right)\label{eq:Iterating H_bold-j on Chi_H}
\end{equation}
What we want to do is write:
\begin{equation}
H_{\mathbf{j}}\left(\lim_{k\rightarrow\infty}H_{\mathbf{j}^{\wedge\left(k-1\right)}}\left(\mathfrak{y}\right)\right)=\lim_{k\rightarrow\infty}H_{\mathbf{j}^{\wedge\left(k-1\right)}}\left(\mathfrak{y}\right),\textrm{ }\forall\mathfrak{y}\in\mathbb{Z}_{q}
\end{equation}
from which we can then conclude that that the $q$-adic integer $\chi_{q}\left(B_{2}\left(n\right)\right)$
is fixed by $H_{\mathbf{j}}$. To make this rigorous, note that $\gcd\left(a_{j},2\right)=1$
for all $j$. As such, for any element of $\textrm{String}\left(2\right)$\textemdash such
as our $\mathbf{j}$\textemdash the quantities $H_{\mathbf{j}}^{\prime}\left(0\right)$
and $H_{\mathbf{j}}\left(0\right)$ are rational numbers \emph{which
lie in} $\mathbb{Z}_{q}$. This guarantees that the function $\mathfrak{z}\mapsto H_{\mathbf{j}}^{\prime}\left(0\right)\mathfrak{z}+H_{\mathbf{j}}\left(0\right)$
(which is, of course $H_{\mathbf{j}}\left(\mathfrak{z}\right)$) is
then a well-defined \emph{continuous} map $\mathbb{Q}_{q}\rightarrow\mathbb{Q}_{q}$
whose restriction to $\mathbb{Z}_{q}$ homeomorphically maps $\mathbb{Z}_{q}$
onto $\mathbb{Z}_{q}$. Using \textbf{Remark \ref{rem:Self-concatenation as a rising limit}},\textbf{
}applying $H_{\mathbf{j}}$ to $\chi_{q}\left(B_{2}\left(n\right)\right)$,
we use the continuity of $H_{\mathbf{j}}$ on $\mathbb{Q}_{q}$ to
obtain: 
\begin{align}
H_{\mathbf{j}}\left(\chi_{q}\left(B_{2}\left(n\right)\right)\right) & \overset{\mathbb{Z}_{q}}{=}H_{\mathbf{j}}\left(\lim_{k\rightarrow\infty}\chi_{q}\left(\mathbf{j}^{\wedge\left(k-1\right)}\right)\right)\\
\left(H_{\mathbf{j}}:\mathbb{Q}_{q}\rightarrow\mathbb{Q}_{q}\textrm{ is cont.}\right); & \overset{\mathbb{Z}_{q}}{=}\lim_{k\rightarrow\infty}H_{\mathbf{j}}\left(\chi_{q}\left(\mathbf{j}^{\wedge\left(k-1\right)}\right)\right)\\
 & \overset{\mathbb{Z}_{q}}{=}\lim_{k\rightarrow\infty}H_{\mathbf{j}}\left(H_{\mathbf{j}^{\wedge\left(k-1\right)}}\left(0\right)\right)\\
 & \overset{\mathbb{Z}_{q}}{=}\lim_{k\rightarrow\infty}H_{\mathbf{j}^{\wedge k}}\left(0\right)\\
 & \overset{\mathbb{Z}_{q}}{=}\lim_{k\rightarrow\infty}\chi_{q}\left(\mathbf{j}^{\wedge k}\right)\\
\left(\textrm{by \textbf{Remark \ref{rem:Self-concatenation as a rising limit}}}\right); & \overset{\mathbb{Z}_{q}}{=}\chi_{q}\left(B_{2}\left(n\right)\right)
\end{align}
This proves the claim.
\end{claim}
\begin{rem}
The moral of \textbf{Claim \ref{claim:2.2}} is that we must\emph{
}consider both the $2$-adic \emph{and }$q$-adic aspects of the Shortened
$qx+1$ map in order to get something meaningful out of it.
\end{rem}
Now, let us actually \emph{assume} that $n$ satisfies the hypotheses
of (II):\textbf{\emph{ suppose $\chi_{q}\left(B_{2}\left(n\right)\right)\in\mathbb{Z}_{2}$}}.
By \textbf{Claim \ref{claim:2.2}}, letting $\mathbf{j}$ by the shortest
string representing $n$, we know that $H_{\mathbf{j}}\left(\chi_{q}\left(B_{2}\left(n\right)\right)\right)=\chi_{q}\left(B_{2}\left(n\right)\right)$.
Applying \textbf{Lemma \ref{lem:properness lemma}}, we conclude:
\begin{equation}
H^{\circ\left|\mathbf{j}\right|}\left(\chi_{q}\left(B_{2}\left(n\right)\right)\right)=\chi_{q}\left(B_{2}\left(n\right)\right)
\end{equation}
where $\left|\mathbf{j}\right|=\lambda_{2}\left(n\right)\geq1$. So,
$\chi_{q}\left(B_{2}\left(n\right)\right)$ is a periodic point of
$H$ in $\mathbb{Z}_{2}$. In particular, if $\chi_{q}\left(B_{2}\left(n\right)\right)$
is in $\mathbb{Z}$, then every element of the cycle generated by
applying the $2$-adic extension of $H$ to $\chi_{q}\left(B_{2}\left(n\right)\right)$
is an element of $\mathbb{Z}$, including $\chi_{q}\left(B_{2}\left(n\right)\right)$
itself when we return to $\chi_{q}\left(B_{2}\left(n\right)\right)$
after $\left|\mathbf{j}\right|$ applications of $H$. This proves
$\chi_{q}\left(B_{2}\left(n\right)\right)$ is a periodic point of
$H$ in $\mathbb{Z}$.

Q.E.D.

\vphantom{}

Using the functional equation in \textbf{Lemma \ref{lem:Chi_H o B_p functional equation}},
we can restate the Correspondence Principle in terms of $\chi_{q}$
and $M_{q}$ in an explicit formula which directly generalizes (and
simplifies) the \textbf{Böhm-Sontacchi Criterion} (\textbf{Theorem
\ref{thm:Bohm-Sontacchi}}).
\begin{cor}[\textbf{Correspondence Principle (Periodic Points), Ver. 2}]
\label{cor:CP v2}

\vphantom{}

I. Let $\Omega\subseteq\mathbb{Z}\backslash\left\{ 0\right\} $ be
any cycle of $T_{q}$. Then, viewing $\Omega$ as a subset of $\mathbb{Z}_{q}$,
the intersection $\chi_{q}\left(\mathbb{Z}_{2}\right)\cap\Omega$
is non-empty. Moreover, for every $x\in\chi_{q}\left(B_{2}\left(\mathbb{N}_{1}\right)\right)\cap\Omega$,
there is an $n\in\mathbb{N}_{1}$ so that: 
\begin{equation}
x=\frac{\chi_{q}\left(n\right)}{1-M_{q}\left(n\right)}
\end{equation}

\vphantom{}

II. Let $n\in\mathbb{N}_{1}$. If the quantity $x$ given by: 
\begin{equation}
x=\frac{\chi_{q}\left(n\right)}{1-M_{q}\left(n\right)}
\end{equation}
is a $2$-adic integer, then $x$ is a periodic point of $T_{q}$
in $\mathbb{Z}_{2}$; if $x$ is in $\mathbb{Z}$, then $x$ is a
periodic point of $T_{q}$ in $\mathbb{Z}$. Moreover, if $x\in\mathbb{Z}$,
then $x$ is positive if and only if $M_{q}\left(n\right)<1$, and
$x$ is negative if and only if $M_{q}\left(n\right)>1$.
\end{cor}
Proof: Re-write the results of \textbf{Theorem \ref{thm:CP v1}} using
\textbf{Lemma \ref{lem:Chi_H o B_p functional equation}}. The positivity/negativity
of $x$ stipulated in (II) follows by noting that $\chi_{q}\left(n\right)$
and $M_{q}\left(n\right)$ are positive rational numbers for all $n\in\mathbb{N}_{1}$.

Q.E.D.

\vphantom{}Finally, we can state and prove the most elegant version
of the CP. This version is also of independent interest because it
allows us to almost trivially establish a connection between $\chi_{q}$
and $T_{q}$'s divergent points.
\begin{cor}[\textbf{Correspondence Principle (Periodic Points), Ver. 4}\footnote{Version 3 can be found in the author's dissertation.}]
\label{cor:CP v4}The set of all non-zero periodic points of $T_{q}$
in $\mathbb{Z}$ is equal to $\mathbb{Z}\cap\chi_{q}\left(\mathbb{Q}\cap\mathbb{Z}_{2}^{\prime}\right)$,
where $\chi_{q}\left(\mathbb{Q}\cap\mathbb{Z}_{2}^{\prime}\right)$
is viewed as a subset of $\mathbb{Z}_{q}$.
\end{cor}
Proof: (Recall, we write $H$ to denote $T_{q}$.)

I. Let $x$ be a non-zero periodic point of $H$, and let $\Omega$
be the unique cycle of $H$ in $\mathbb{Z}$ containing $x$. By \textbf{Version
1 of the} \textbf{CP} (\textbf{Theorem \ref{thm:CP v1}}), there exists
a $y\in\Omega$ and a $\mathfrak{z}\in\mathbb{Z}_{2}$ (where $\mathfrak{z}=B_{2}\left(n\right)$
for some $n\in\mathbb{N}_{1}$) so that $\chi_{q}\left(\mathfrak{z}\right)=y$.
Since $y\in\Omega$, there is a $k\geq1$ so that $x=H^{\circ k}\left(y\right)$.
In particular, there is a \emph{unique} string $\mathbf{i}\in\textrm{String}\left(2\right)$
of length $k$ so that $H_{\mathbf{i}}\left(y\right)=H^{\circ k}\left(y\right)=x$.

Now, let $\mathbf{j}\in\textrm{String}_{\infty}\left(2\right)$ represent
$\mathfrak{z}$; note that $\mathbf{j}$ is infinite and that its
entries are periodic, and contain infinitely many $1$s. Using $\chi_{q}$'s
concatenation identity (\textbf{Lemma \ref{lem:Chi_H concatenation identity}}),
we can write: 
\begin{equation}
x=H_{\mathbf{i}}\left(y\right)=H_{\mathbf{i}}\left(\chi_{q}\left(\mathfrak{z}\right)\right)=\chi_{q}\left(\mathbf{i}\wedge\mathbf{j}\right)
\end{equation}

Next, let $\mathfrak{x}$ denote the $2$-adic integer represented
by $\mathbf{i}\wedge\mathbf{j}$; note that $\mathfrak{x}$ is then
\emph{not} an element of $\mathbb{N}_{0}$, because it has infinitely
many digits. By the above, we have that $\chi_{q}\left(\mathfrak{x}\right)=x$,
and hence that $x\in\mathbb{Z}\cap\chi_{q}\left(\mathbb{Z}_{2}\right)$.

Finally, since the $2$-adic digits of $\mathfrak{z}=B_{2}\left(n\right)$
form a periodic sequence, it must be that $\mathfrak{z}$ is a rational
number; in particular, an element of $\mathbb{Q}\cap\mathbb{Z}_{2}^{\prime}$.
Indeed, $\mathfrak{z}$ is not in $\mathbb{N}_{0}$ because $n\geq1$,
and $B_{2}\left(n\right)\in\mathbb{N}_{0}$ if and only if $n=0$.
So, letting $m$ be the rational integer represented by $\mathbf{i}$,
we have: 
\begin{equation}
\mathfrak{x}\sim\mathbf{i}\wedge\mathbf{j}\sim m+2^{\lambda_{2}\left(m\right)}\mathfrak{z}
\end{equation}
This shows that $\mathfrak{x}\in\mathbb{Q}\cap\mathbb{Z}_{2}^{\prime}$,
and hence, that $x=\chi_{q}\left(\mathfrak{x}\right)\in\mathbb{Z}\cap\chi_{q}\left(\mathbb{Q}\cap\mathbb{Z}_{2}^{\prime}\right)$.

\vphantom{}

II. Suppose $x\in\mathbb{Z}\cap\chi_{q}\left(\mathbb{Q}\cap\mathbb{Z}_{2}^{\prime}\right)$,
with $x=\chi_{q}\left(\mathfrak{z}\right)$ for some $\mathfrak{z}\in\mathbb{Q}\cap\mathbb{Z}_{2}^{\prime}$.
As a rational number which is both a $2$-adic integer \emph{and}
not an element of $\mathbb{N}_{0}$, the $2$-adic digits of $\mathfrak{z}$
are eventually periodic. As such, there are integers $m$ and $n$
(with $n\neq0$) so that: 
\begin{equation}
\mathfrak{z}=m+2^{\lambda_{2}\left(m\right)}B_{2}\left(n\right)
\end{equation}
Here, $n$'s $2$-adic digits generate the periodic part of $\mathfrak{z}$'s
digits, while $m$'s $2$-adic digits are the finite-length sequence
of $\mathfrak{z}$'s digits which occur before the periodicity sets
in.

Now, let $\mathbf{i}$ be the finite string representing $m$, and
let $\mathbf{j}$ be the infinite string representing $B_{2}\left(n\right)$.
Then, $\mathfrak{z}=\mathbf{i}\wedge\mathbf{j}$. So, by \textbf{Lemmata
\ref{lem:Chi_H concatenation identity}} and \textbf{\ref{lem:Chi_H o B_p functional equation}}
($\chi_{q}$'s concatenation identity and the $\chi_{q}\circ B_{2}$
functional equation, respectively): 
\begin{equation}
x=\chi_{q}\left(\mathbf{i}\wedge\mathbf{j}\right)=H_{\mathbf{i}}\left(\chi_{q}\left(\mathbf{j}\right)\right)=H_{\mathbf{i}}\left(\chi_{q}\left(B_{2}\left(n\right)\right)\right)=H_{\mathbf{i}}\left(\frac{\chi_{q}\left(n\right)}{1-M_{q}\left(n\right)}\right)
\end{equation}
where: 
\begin{equation}
y\overset{\textrm{def}}{=}\frac{\chi_{q}\left(n\right)}{1-M_{q}\left(n\right)}
\end{equation}
is a rational number.
\begin{claim}
\label{claim:2.2-1}$\left|y\right|_{2}\leq1$.

Proof of claim: Since $y$ is a rational number, it lies in $\mathbb{Q}_{2}$.
So, by way of contradiction, suppose $\left|y\right|_{2}>1$. By \textbf{Proposition
\ref{prop:Q_p / Z_p prop}}, this implies $\left|H_{\mathbf{i}}\left(y\right)\right|_{2}>1$.
However, $H_{\mathbf{i}}\left(y\right)=x$, and $x$ is a rational
integer; hence, $1<\left|H_{\mathbf{i}}\left(y\right)\right|_{2}=\left|x\right|_{2}\leq1$.
This is impossible! So, it must be that $\left|y\right|_{2}\leq1$.
This proves the claim.
\end{claim}
\begin{claim}
\label{claim:2.3}$x=H^{\circ\left|\mathbf{i}\right|}\left(y\right)$

Proof of claim: Suppose the equality failed. Then $H_{\mathbf{i}}\left(y\right)=x\neq H^{\circ\left|\mathbf{i}\right|}\left(y\right)$,
and so $x=H_{\mathbf{i}}\left(y\right)$ is a wrong value of $H$
with seed $y$. \textbf{Lemma \ref{lem:wrong values lemma}} then
forces $\left|x\right|_{2}=\left|H_{\mathbf{i}}\left(y\right)\right|_{2}>1$.
However, $\left|x\right|_{2}\leq1$. This is just as impossible as
it was in the previous paragraph. This proves the claim.
\end{claim}
Finally, let $\mathbf{v}$ be the shortest string representing $n$,
so that $\mathbf{j}$ (the string representing $B_{2}\left(n\right)$)
is obtained by concatenating infinitely many copies of $\mathbf{v}$.
Once again, we note that\emph{ $H_{\mathbf{v}}$ is a continuous self-map
of $\mathbb{Z}_{q}$} (\textbf{Claim \ref{claim:2.2}} from the proof
of \textbf{Version 1 of the} \textbf{CP}). As such:
\begin{equation}
\chi_{q}\left(B_{2}\left(n\right)\right)\overset{\mathbb{Z}_{q}}{=}\lim_{k\rightarrow\infty}\chi_{q}\left(\mathbf{v}^{\wedge k}\right)
\end{equation}
implies:
\begin{align}
H_{\mathbf{v}}\left(\chi_{q}\left(B_{2}\left(n\right)\right)\right) & \overset{\mathbb{Z}_{q}}{=}\lim_{k\rightarrow\infty}H_{\mathbf{v}}\left(\chi_{q}\left(\mathbf{v}^{\wedge k}\right)\right)\\
 & \overset{\mathbb{Z}_{q}}{=}\lim_{k\rightarrow\infty}\chi_{q}\left(\mathbf{v}^{\wedge\left(k+1\right)}\right)\\
 & \overset{\mathbb{Z}_{q}}{=}\chi_{q}\left(B_{2}\left(n\right)\right)
\end{align}
Hence, $H_{\mathbf{v}}\left(y\right)=y$. Since $\left|y\right|_{2}\leq1$
is fixed by $H_{\mathbf{v}}$,\textbf{ Lemma \ref{lem:properness lemma}}
applies, and we conclude that $y$ is a periodic point of $H:\mathbb{Z}_{2}\rightarrow\mathbb{Z}_{2}$.

By \textbf{Claim \ref{claim:2.3}}, $H$ iterates $y$ to $x$. Since
$y$ is a periodic point of $H$ in $\mathbb{Z}_{2}$, this forces
$x$ and $y$ to belong to the same cycle of $H$ in $\mathbb{Z}_{2}$,
and forces $x$ to be a periodic point of $H$ as well. As such, just
as $H$ iterates $y$ to $x$, so too does $H$ iterate $x$ to $y$.
Likewise, since $x$ is a rational integer, so too is $y$ (by \textbf{Proposition
\ref{prop:Q_p / Z_p prop}}). Thus, $x$ belongs to a cycle of $H$
in $\mathbb{Z}$, as desired.

Q.E.D.
\begin{example}
To illustrate the CPPP in action, observe that the cycle $1\rightarrow2\rightarrow1$
of $H=T_{3}$ applies the even branch ($H_{0}$) second and the odd
branch ($H_{1}$) first. Thus, the string $\mathbf{j}$ such that
$H_{\mathbf{j}}\left(1\right)=1$ is $\mathbf{j}=\left(0,1\right)$:
\begin{equation}
1=H_{0}\left(H_{1}\left(1\right)\right)=H_{0}\left(2\right)=1
\end{equation}
The integer $n$ represented by $\mathbf{j}$ is $n=0\cdot2^{0}+1\cdot2^{1}=2$.
Thus: 
\begin{equation}
\centerdot01010101\ldots\overset{\mathbb{Z}_{2}}{=}B_{2}\left(2\right)=\frac{2}{1-2^{\lambda_{2}\left(2\right)}}=\frac{2}{1-4}=-\frac{2}{3}
\end{equation}
and so: 
\begin{equation}
\chi_{3}\left(-\frac{2}{3}\right)=\chi_{3}\left(B_{2}\left(2\right)\right)=\frac{\chi_{3}\left(2\right)}{1-M_{3}\left(2\right)}=\frac{\frac{1}{4}}{1-\frac{3}{4}}=\frac{1}{4-3}=1
\end{equation}
where: 
\begin{equation}
\chi_{3}\left(2\right)=\frac{1}{2}\chi_{3}\left(1\right)=\frac{1}{2}\left(\frac{3\chi_{3}\left(0\right)+1}{2}\right)=\frac{1}{2}\cdot\frac{0+1}{2}=\frac{1}{4}
\end{equation}
and: 
\begin{equation}
M_{3}\left(2\right)=\frac{3^{\#_{1}\left(2\right)}}{2^{\lambda_{2}\left(2\right)}}=\frac{3^{1}}{2^{2}}=\frac{3}{4}
\end{equation}
\end{example}
\begin{rem}
\label{rem:The-general-procedure}Returning to \textbf{Corollary \ref{cor:CP v2}}
for a moment, what that result shows is that, while $\left(\chi_{q}\circ B_{2}\right)\left(n\right)$
is a periodic point of $T_{q}$ for any $n\in\mathbb{N}_{0}$ which
makes $\left(\chi_{q}\circ B_{2}\right)\left(n\right)\in\mathbb{Z}$,
not every periodic point of $T_{q}$ has the form $\left(\chi_{q}\circ B_{2}\right)\left(n\right)$
for some $n\in\mathbb{N}_{0}$. These exceptional periodic points
are simply those rational integer values attained by $\chi_{q}\left(\mathfrak{z}\right)$
for $\mathfrak{z}\in\mathbb{Z}_{2}\cap\mathbb{Q}$ which \emph{do
not }lie in $B_{2}\left(\mathbb{N}_{0}\right)$.

Given a cycle $\Omega\subset\mathbb{Z}\backslash\left\{ 0\right\} $
of $T_{q}$, there is a general procedure for finding the $\mathfrak{z}$s
for each element of $\Omega$, which we now describe. Suppose we have
found an $n\in\mathbb{N}_{1}$ so that $\mathfrak{z}_{0}\overset{\textrm{def}}{=}B_{2}\left(n\right)$
satisfies $\chi_{q}\left(\mathfrak{z}_{0}\right)\in\mathbb{Z}\backslash\left\{ 0\right\} $,
and hence, so that $\chi_{q}\left(\mathfrak{z}_{0}\right)$ is a periodic
point of $H$. Letting $x_{0}\overset{\textrm{def}}{=}\chi_{H}\left(\mathfrak{z}_{0}\right)$,
observe that:
\begin{equation}
\Omega=\left\{ H^{\circ k}\left(x_{0}\right):0\leq k\leq\left|\Omega\right|-1\right\} 
\end{equation}
Letting $x_{k}\overset{\textrm{def}}{=}H^{\circ k}\left(x_{0}\right)$,
we have $x_{k}=H^{\circ k}\left(\chi_{q}\left(\mathfrak{z}_{0}\right)\right)$.
If, for example, $\chi_{q}\left(\mathfrak{z}_{0}\right)$ is odd,
so that:
\begin{equation}
x_{1}=H\left(\chi_{q}\left(\mathfrak{z}_{0}\right)\right)=H_{1}\left(\chi_{q}\left(\mathfrak{z}_{0}\right)\right)
\end{equation}
the functional equation for $\chi_{q}$ tells us that:
\begin{equation}
x_{1}=H\left(\chi_{q}\left(\mathfrak{z}_{0}\right)\right)=H_{1}\left(\chi_{q}\left(\mathfrak{z}_{0}\right)\right)=\frac{q\chi_{q}\left(\mathfrak{z}_{0}\right)+1}{2}=\chi_{q}\left(2\mathfrak{z}_{0}+1\right)
\end{equation}
In particular, letting $\mathfrak{z}_{k}$ denote a $2$-adic integer
so that $x_{k}=\chi_{q}\left(\mathfrak{z}_{k}\right)$, we then have:
\begin{equation}
\mathfrak{z}_{k+1}=\begin{cases}
2\mathfrak{z}_{k} & \textrm{if }x_{k}\overset{2}{\equiv}0\\
2\mathfrak{z}_{k}+1 & \textrm{if }x_{k}\overset{2}{\equiv}1
\end{cases},\textrm{ }\forall k\in\left\{ 0,\ldots,\left|\Omega\right|-1\right\} 
\end{equation}

This can be understood in terms of cyclic permutations. In particular,
$\mathfrak{z}_{k}$ will be obtained from $\mathfrak{z}_{0}$ by adding
the $k$ right-most digits of $n$ to the left side of $\mathfrak{z}_{0}$.
So, for example, if $n$ has the $2$-adic digit sequence:
\begin{equation}
n=\centerdot_{2}1001=1+1\cdot2^{3}=9
\end{equation}
so that:
\begin{equation}
\mathfrak{z}_{0}=B_{2}\left(n\right)=\centerdot_{2}100110011001\ldots
\end{equation}
then we then have:
\begin{align}
\mathfrak{z}_{1} & =\centerdot_{2}1100110011001\ldots\\
\mathfrak{z}_{2} & =\centerdot_{2}01100110011001\ldots\\
\mathfrak{z}_{3} & =\centerdot_{2}001100110011001\ldots
\end{align}
and then $\mathfrak{z}_{4}=\mathfrak{z}_{0}$. Here, only $\mathfrak{z}_{0}$
and $\mathfrak{z}_{3}$ are the images of positive integers under
$B_{2}$. This is because they are generated by the sequences $1001$
and $0011$, respectively. The generating sequences of $\mathfrak{z}_{1}$
and $\mathfrak{z}_{2}$ are $1100$ and $0110$. Our identification
of strings with ($2$-adic) integers ignores all the $0$s to the
right of the final $1$, and so:
\begin{align}
B_{2}\left(\centerdot_{2}1100\right) & =B_{2}\left(\centerdot_{2}11\right)=\centerdot_{2}11111\ldots\\
B_{2}\left(\centerdot_{2}0110\right) & =B_{2}\left(\centerdot_{2}011\right)=\centerdot_{2}011011\ldots
\end{align}

Finally, noting that the rational number $\chi_{q}\left(B_{2}\left(n\right)\right)$
has a $2$-adic absolute value of $1$ for all $n\in\mathbb{N}_{1}$,
we can conclude that the elements of $\Omega$ of the form $\chi_{q}\left(B_{2}\left(n\right)\right)$
for some $n\in\mathbb{N}_{1}$ are precisely the odd integers in $\Omega$.
That is, an odd integer is a periodic point of $T_{q}$ if and only
if it is of the form $\chi_{q}\left(B_{2}\left(n\right)\right)$ for
some $n\in\mathbb{N}_{1}$.
\end{rem}
\vphantom{}

From the above example, we observe that the number of elements of
$\Omega\subseteq\mathbb{Z}\backslash\left\{ 0\right\} $ is equal
to the number of distinct cyclic permutations of the $2$-adic digits
of any $\mathfrak{z}\in\mathbb{Z}_{2}$ satisfying $\chi_{q}\left(\mathfrak{z}\right)\in\Omega$.
In particular, this implies that $\left|\Omega\right|=\lambda_{2}\left(n\right)$,
where $n$ is the smallest positive integer so that $\chi_{q}\left(B_{2}\left(n\right)\right)\in\Omega$. 

Having dealt with the periodic points, all that remains is to obtain
the partial CP for the divergent points of $T_{q}$ (the CPDP). To
do this, we need only make a simple observation: by \textbf{Corollary
\ref{cor:CP v4}}, we know that the periodic points of $T_{q}$ in
$\mathbb{Z}\backslash\left\{ 0\right\} $ are completely determined
by the $\mathbb{Z}$-values $\chi_{q}$ attains on $\mathbb{Z}_{2}^{\prime}\cap\mathbb{Q}$.
That is to say, we have a correspondence between periodic points of
$T_{q}$ and rational integer values attained by $\chi_{q}$ on \emph{rational
}$2$-adic integer inputs.

\vphantom{}

So, what happens if $\chi_{q}\left(\mathfrak{z}\right)\in\mathbb{Z}$
for some \emph{irrational }$2$-adic integer $\mathfrak{z}$ (i.e.,
$\mathfrak{z}\in\mathbb{Z}_{2}\backslash\mathbb{Q}$)?

\vphantom{}

Before we can answer this question, we need the a technical result:
\begin{prop}
\label{prop:non-vanishing of chi_3 off 0}$\chi_{q}\left(\mathfrak{z}\right)=0$
if and only if $\mathfrak{z}=0$.
\end{prop}
Proof: Since $\chi_{q}\left(0\right)$ is $0$ by definition, to complete
the proof, we need only show that $\chi_{q}\left(\mathfrak{z}\right)\neq0$
for all $\mathfrak{z}\in\mathbb{Z}_{2}\backslash\left\{ 0\right\} $.
So, suppose that $\mathfrak{z}\in\mathbb{Z}_{2}\backslash\left\{ 0\right\} $
makes $\chi_{q}\left(\mathfrak{z}\right)=0$. Then, since $\mathfrak{z}$
is non-zero, there are unique $n\in\mathbb{N}_{0}$ and $\mathfrak{u}\in\mathbb{Z}_{2}^{\times}$
so that $\mathfrak{z}=2^{n}\mathfrak{u}$. Since $\mathfrak{u}$ is
odd, the functional equations of $\chi_{q}$ yield:
\begin{equation}
0=\chi_{q}\left(\mathfrak{z}\right)=\chi_{q}\left(2^{n}\mathfrak{u}\right)=\frac{1}{2^{n}}\chi_{q}\left(\mathfrak{u}\right)=\frac{q\chi_{q}\left(\mathfrak{v}\right)+1}{2^{n+1}}
\end{equation}
where: 
\begin{equation}
\mathfrak{v}\overset{\textrm{def}}{=}\frac{\mathfrak{u}-1}{2}\in\mathbb{Z}_{2}
\end{equation}
Solving for $\chi_{q}\left(\mathfrak{v}\right)$ yields:
\begin{equation}
\chi_{q}\left(\mathfrak{v}\right)=-\frac{1}{q}
\end{equation}
However, $-1/q$ is not a $q$-adic integer, and $\chi_{q}$ only
outputs $q$-adic integers. This is a contradiction!

Thus, $\mathfrak{z}$ must be $0$. This proves $0$ is the only zero
of $\chi_{q}$ in $\mathbb{Z}_{2}$.

Q.E.D.

\vphantom{}

Answering our question about irrational $\mathfrak{z}$ immediately
yields the partial CP for divergent points.
\begin{thm}[\textbf{\textit{Correspondence Principle (Divergent Points)}}]
\label{thm:Divergent trajectories come from irrational z}Let $\mathfrak{z}\in\mathbb{Z}_{2}\backslash\mathbb{Q}$
be such that $\chi_{q}\left(\mathfrak{z}\right)\in\mathbb{Z}$. Then
$\chi_{q}\left(\mathfrak{z}\right)$ is a divergent point of $T_{q}$.
\end{thm}
\begin{rem}
\textbf{Theorem \ref{thm:Divergent trajectories come from irrational z}}
is false. See \textbf{Section \ref{sec:Corrections}} for further
discussion.
\end{rem}
Proof: (Recall we write $H$ to denote $T_{q}$) Let $\mathfrak{z}\in\mathbb{Z}_{2}\backslash\mathbb{Q}$
make $\chi_{q}\left(\mathfrak{z}\right)\in\mathbb{Z}$. Since $\mathfrak{z}$
is irrational, it cannot be $0$. So, \textbf{Proposition \ref{prop:non-vanishing of chi_3 off 0}
}tells us that $\chi_{q}\left(\mathfrak{z}\right)\neq0$. Now, observe
that $\chi_{q}\left(\mathfrak{z}\right)$ cannot be a \emph{periodic
point}\textbf{ }of $H$. If it were, then, as a non-zero periodic
point of $H$ in $\mathbb{Z}$,\textbf{ Corollary \ref{cor:CP v4}
}would guarantee that $\mathfrak{z}\in\mathbb{Q}\cap\mathbb{Z}_{2}^{\prime}$,
which is impossible: $\mathfrak{z}$ is irrational!

So, $\chi_{q}\left(\mathfrak{z}\right)$ is an integer which is not
a periodic point of $H$. By \textbf{Proposition \ref{prop:Classification of T_3's dynamics}},
every integer is either a a pre-periodic point of $H$ or a divergent
point of $H$. So, by way of contradiction, \emph{suppose $\chi_{q}\left(\mathfrak{z}\right)$
is a pre-periodic point of }$H$. Then, there is an $n\geq0$ so that
the integer $H^{\circ n}\left(\chi_{q}\left(\mathfrak{z}\right)\right)$
is a periodic point of $H$, and a necessarily non-zero one at that
(because the only integer $H$ sends to $0$ is $0$ itself). By \textbf{Corollary
\ref{cor:CP v4}}, there is then a $\mathfrak{y}\in\mathbb{Q}\cap\mathbb{Z}_{2}^{\prime}$
so that $\chi_{q}\left(\mathfrak{y}\right)=H^{\circ n}\left(\chi_{q}\left(\mathfrak{z}\right)\right)$
is a periodic point of $H$.

Letting $\mathbf{j}$ be the string representing $\mathfrak{z}$,
we have that:
\begin{equation}
\chi_{q}\left(\mathfrak{y}\right)=H^{\circ n}\left(\chi_{q}\left(\mathbf{j}\right)\right)
\end{equation}
Letting $\mathbf{i}\in\textrm{String}\left(2\right)$ be the length
$n$ string which gives the motions of $\chi_{q}\left(\mathbf{j}\right)$
under $H^{\circ n}$, we have
\begin{equation}
\chi_{q}\left(\mathfrak{y}\right)=H^{\circ n}\left(\chi_{q}\left(\mathbf{j}\right)\right)=H_{\mathbf{i}}\left(\chi_{q}\left(\mathbf{j}\right)\right)=\chi_{q}\left(\mathbf{i}\wedge\mathbf{j}\right)
\end{equation}
Thus, $\mathbf{i}\wedge\mathbf{j}$ represents $\mathfrak{y}$, while
$\mathbf{j}$ represents $\mathfrak{z}$. Since $\mathfrak{y}$ is
rational, its digits\textemdash the entires of $\mathbf{i}\wedge\mathbf{j}$\textemdash are
eventually periodic. Since $\mathbf{i}$ has finite length, this forces
the entires of $\mathbf{j}$ to be eventually periodic, as well. But
this is impossible: $\mathfrak{z}$ is an \emph{irrational }$2$-adic
integer, so its digits\textemdash the entries of $\mathbf{j}$\textemdash \emph{cannot
be }be eventually periodic. This is the desired contradiction: our
assumption that $\chi_{q}\left(\mathfrak{z}\right)$ was\emph{ }a
pre-periodic point of $H$ must have been in error.

So, by \textbf{Proposition \ref{prop:Classification of T_3's dynamics}}
(page \pageref{prop:Classification of T_3's dynamics}), since $\chi_{q}\left(\mathfrak{z}\right)$
is an element of $\mathbb{Z}$ which is not a pre-periodic point of
$H$, it must be a divergent point of $H$.

Q.E.D.
\begin{rem}
\newpage{}
\end{rem}

\section{Horizons\label{sec:3 Horizons}}

\subsection{\label{subsec:The-Tao-of}The Tao of Collatz}

The origins of the present work lie in ideas the author began exploring
in late 2019 and early 2020. Around the same time\textemdash late
2019\textemdash the esteemed Terence Tao released a paper on the Collatz
Conjecture (\textquotedblleft Almost all orbits of the Collatz map
attain almost bounded values\textquotedblright ) that applied probabilistic
techniques in a novel way, by use of what Tao called ``Syracuse Random
Variables'' to obtain the conclusion stated in his paper's title.
Despite the many substantial differences between our two approaches,
$\chi_{3}$ links them at a fundamental level.

Tao's approach involves constructing his Syracuse Random Variables
and then comparing them to a set-up involving tuples of geometric
random variables, the comparison in question being an estimate on
the 'distance' between the Syracuse Random Variables and the geometric
model, as measured by the total variation norm for discrete random
variables. The central challenge Tao overcomes to obtain these results
is to obtain explicit estimates for the decay of the characteristic
function of the Syracuse Random Variables.

Surprisingly, $\chi_{3}$ is equivalent to Tao's Syracuse Random Variables.
The Syracuse Random Variables are a family of random variables taking
values in $\mathbb{Z}/3^{n}\mathbb{Z}$, where $n\geq0$ is arbitrary.
Tao explicitly mentions early on in his paper that one can adopt the
point of view in which the Syracuse Random Variables are merely the
projections modulo $3^{n}$ of a \emph{single }$3$-adic valued random
variable; that random variable, of course, is $\chi_{3}$. Nevertheless,
Tao chose not to pursue this formalism.

The ``most difficult'' result Tao establishes in his paper\textemdash the
heart of his proof\textemdash is, in the terminology of \cite{my dissertation},
an estimate on the decay of complex-valued function $\varphi_{3}:\hat{\mathbb{Z}}_{3}\rightarrow\mathbb{C}$
defined by:
\begin{equation}
\varphi_{3}\left(t\right)\overset{\textrm{def}}{=}\int_{\mathbb{Z}_{2}}e^{-2\pi i\left\{ t\chi_{3}\left(\mathfrak{z}\right)\right\} _{3}}d\mathfrak{z},\textrm{ }\forall t\in\hat{\mathbb{Z}}_{3}\label{eq:Tao's Characteristic Function}
\end{equation}
Here, $\left\{ \cdot\right\} _{3}:\mathbb{Q}_{3}\rightarrow\hat{\mathbb{Z}}_{3}$
is the $3$-adic fractional part, and $d\mathfrak{z}$ is the real-valued
Haar probability measure on $\mathbb{Z}_{2}$. Tao's decay estimate
is given in \textbf{Proposition 1.17 }of his paper, where the above
integral appears in the form of an expected value \cite{Tao Probability paper}.
Using our terminology, Tao's result can be stated as:
\begin{quotation}
Let $n\geq1$. Then, for any real number $A>0$, there is a constant
$C$ so that:
\begin{equation}
\left|\varphi_{3}\left(t\right)\right|\leq\frac{C}{\left(\log_{3}\left|t\right|_{3}\right)^{A}},\textrm{ }\forall t\in\hat{\mathbb{Z}}_{3}\backslash\left\{ 0\right\} \label{eq:Tao's Result}
\end{equation}
\end{quotation}
Here:
\begin{equation}
\log_{3}x=\frac{\ln x}{\ln3}
\end{equation}
is the base-$3$ logarithm.

By exploiting the discrete Fourier transform, we can use $\varphi_{3}$
to write:
\begin{equation}
\textrm{P}\left(\chi_{3}\overset{3^{n}}{\equiv}m\right)=\frac{1}{3^{n}}\sum_{\left|t\right|_{3}\leq3^{n}}\varphi_{3}\left(t\right)e^{2\pi imt}\label{eq:Discrete Fourier transform of Phi_3 hat}
\end{equation}
Here, the sum is taken over all elements of $\hat{\mathbb{Z}}_{3}$
with denominators at most $3^{n}$. As indicated, $\textrm{P}\left(\chi_{3}\overset{3^{n}}{\equiv}m\right)$
is the probability that $\chi_{3}$ is congruent to $m$ mod $3^{n}$,
which is to say, it is measure (according to the real-valued $2$-adic
Haar probability measure) of the set of $\mathfrak{z}\in\mathbb{Z}_{2}$
for which $\chi_{3}\left(\mathfrak{z}\right)\overset{3^{n}}{\equiv}m$.
The $n$th Syracuse Random Variable is then precisely the real-valued
random variable $\mathfrak{z}\mapsto\left[\chi_{3}\left(\mathfrak{z}\right)\right]_{3^{n}}$,
which outputs the value of $\chi_{3}\left(\mathfrak{z}\right)$ modulo
$3^{n}$. The sample space\footnote{In this respect, our investigations of $\chi_{3}$ as a function of
$\mathfrak{z}\in\mathbb{Z}_{2}$ are extremely unusual from the vantage
point of Probability theory, where the fact that random variables
are functions of points in the sample space is, in practice, almost
always swept under the rug.} is $\mathbb{Z}_{2}$, with the usual Borel $\sigma$-algebra. The
reader can verify the equivalence of $\chi_{3}$ and the Syracuse
Random Variables for themselves by using the functional equation satisfied
by $\varphi_{3}$:
\begin{equation}
\varphi_{3}\left(t\right)=\frac{1}{2}\varphi_{3}\left(\frac{t}{2}\right)+\frac{1}{2}e^{-2\pi i\left\{ \frac{t}{2}\right\} _{3}}\varphi\left(\frac{3t}{2}\right),\textrm{ }\forall t\in\hat{\mathbb{Z}}_{3}\label{eq:Phi_3 hat functional equation}
\end{equation}
This is obtained by splitting the integral (\ref{eq:Tao's Characteristic Function})
into the sum of an integral over $2\mathbb{Z}_{2}$ and an integral
over $2\mathbb{Z}_{2}+1$, using the substitutions $\mathfrak{y}=\mathfrak{z}/2$
and $\mathfrak{y}=\left(\mathfrak{z}-1\right)/2$ to get:
\begin{equation}
\varphi_{3}\left(t\right)=\frac{1}{2}\int_{\mathbb{Z}_{2}}e^{-2\pi i\left\{ t\chi_{3}\left(2\mathfrak{y}\right)\right\} _{3}}d\mathfrak{y}+\frac{1}{2}\int_{\mathbb{Z}_{2}}e^{-2\pi i\left\{ t\chi_{3}\left(2\mathfrak{y}+1\right)\right\} _{3}}d\mathfrak{y}
\end{equation}
using the functional equations for $\chi_{3}$ to express the exponents
in terms of $\chi_{3}\left(\mathfrak{y}\right)$ and then expressing
the right-hand side in terms of $\varphi_{3}\left(t\right)$. (\ref{eq:Phi_3 hat functional equation})
can be used to obtain a recursive formula for $\varphi_{3}\left(t\right)$.
Using (\ref{eq:Discrete Fourier transform of Phi_3 hat}) converts
this into a recursive formula for $\textrm{P}\left(\chi_{3}\overset{3^{n}}{\equiv}m\right)$.
Tao gives the formula for $\textrm{P}\left(\chi_{3}\overset{3^{n}}{\equiv}m\right)$
in \textbf{Lemma 1.12 }of his paper. The conversion between Tao's
notation and our own is:
\begin{equation}
\mathbb{P}\left(\mathbf{Syrac}\left(\mathbb{Z}/3^{n}\mathbb{Z}\right)=x\right)=\textrm{P}\left(\chi_{3}\overset{3^{n}}{\equiv}x\right),\textrm{ }\forall x\in\mathbb{Z}/3^{n}\mathbb{Z}
\end{equation}

\subsection{\label{subsec:The-Converse-of-CPDP}An Observation About Backward
Orbits}

We begin an important function.
\begin{defn}
The \textbf{$2$-adic (left) shift map }$\theta_{2}:\mathbb{Z}_{2}\rightarrow\mathbb{Z}_{2}$\textbf{
}is defined by:
\begin{equation}
\theta_{2}\left(\mathfrak{z}\right)\overset{\textrm{def}}{=}\frac{\mathfrak{z}-\left[\mathfrak{z}\right]_{2}}{2},\textrm{ }\forall\mathfrak{z}\in\mathbb{Z}_{2}\label{eq:Definition of the 2-adic left shift map-1}
\end{equation}
\end{defn}
With the shift map, we can give the following characterization of
$\chi_{q}$'s interaction with $\theta_{2}$:
\begin{prop}
\label{prop:Chi_q and the shift map}
\begin{equation}
\chi_{q}\left(\mathfrak{z}\right)=H_{\left[\mathfrak{z}\right]_{2}}\left(\chi_{q}\left(\theta_{2}\left(\mathfrak{z}\right)\right)\right),\textrm{ }\forall\mathfrak{z}\in\mathbb{Z}_{2}\label{eq:Chi in terms of the shift map}
\end{equation}
\end{prop}
Proof: Let $\mathfrak{z}\in\mathbb{Z}_{2}$, and let $c_{n}$ denote
the $n$th $2$-adic digit of $\mathfrak{z}$, with $\mathfrak{z}=\sum_{n=0}^{\infty}c_{n}2^{n}$.
Noting that:
\begin{equation}
\mathfrak{z}=c_{0}+2\sum_{n=1}^{\infty}c_{n}2^{n-1}=\left[\mathfrak{z}\right]_{2}+2\left(\frac{\mathfrak{z}-\left[\mathfrak{z}\right]_{2}}{2}\right)
\end{equation}
Observe that the right-hand side will be of the form $2\mathfrak{y}$
for some $\mathfrak{y}\in\mathbb{Z}_{2}$ if and only if $\mathfrak{z}\overset{2}{\equiv}0$
and will be of the form $1+2\mathfrak{y}$ for some $\mathfrak{y}\in\mathbb{Z}_{2}$
if and only if $\mathfrak{z}\overset{2}{\equiv}1$. If $\mathfrak{z}\overset{2}{\equiv}0$,
we have:
\begin{equation}
\chi_{q}\left(\mathfrak{z}\right)=\chi_{q}\left(2\left(\frac{\mathfrak{z}}{2}\right)\right)=\frac{1}{2}\chi_{q}\left(\frac{\mathfrak{z}}{2}\right)
\end{equation}
If $\mathfrak{z}\overset{2}{\equiv}1$, on the other hand:
\begin{equation}
\chi_{q}\left(\mathfrak{z}\right)=\chi_{q}\left(2\left(\frac{\mathfrak{z}-1}{2}\right)+1\right)=\frac{q\chi_{q}\left(\frac{\mathfrak{z}-1}{2}\right)+1}{2}
\end{equation}
Both cases can be expressed in a single formula by writing:
\begin{equation}
\chi_{q}\left(\mathfrak{z}\right)=H_{\left[\mathfrak{z}\right]_{2}}\left(\chi_{q}\left(\frac{\mathfrak{z}-\left[\mathfrak{z}\right]_{2}}{2}\right)\right)\label{eq:Pre left shift formula}
\end{equation}
where: 
\begin{align}
H_{0}\left(x\right) & \overset{\textrm{def}}{=}\frac{x}{2}\\
H_{1}\left(x\right) & \overset{\textrm{def}}{=}\frac{qx+1}{2}
\end{align}
with $H_{\left[\mathfrak{z}\right]_{2}}$ being the map where the
subscript is the value of $\mathfrak{z}$ mod $2$. Rewriting $\frac{\mathfrak{z}-\left[\mathfrak{z}\right]_{2}}{2}$
as $\theta_{2}\left(\mathfrak{z}\right)$ then gives the desired result.

Q.E.D.

\vphantom{}Suppose $\mathfrak{z}_{0}\in\mathbb{Z}_{2}\backslash\mathbb{Q}$
satisfies $\chi_{q}\left(\mathfrak{z}_{0}\right)\in\mathbb{Z}$. Then,
let $x_{0}\overset{\textrm{def}}{=}\chi_{q}\left(\mathfrak{z}_{0}\right)$.
If $\mathfrak{z}_{0}\overset{2}{\equiv}0$, then:
\begin{equation}
x_{0}=\chi_{q}\left(\mathfrak{z}_{0}\right)=H_{0}\left(\chi_{q}\left(\theta_{2}\left(\mathfrak{z}_{0}\right)\right)\right)=\frac{\chi_{q}\left(\theta_{2}\left(\mathfrak{z}_{0}\right)\right)}{2}
\end{equation}
and so, $\mathfrak{z}_{1}\overset{\textrm{def}}{=}\theta_{2}\left(\mathfrak{z}_{0}\right)$
is an element of $\mathbb{Z}_{2}\backslash\mathbb{Q}$ so that $\chi_{q}\left(\mathfrak{z}_{1}\right)=2x_{0}\in\mathbb{Z}$,
moreover, note that $T_{q}\left(2x_{0}\right)=x_{0}$. Alternatively,
if $\mathfrak{z}_{0}\overset{2}{\equiv}1$, we get:
\begin{equation}
x_{0}=\chi_{q}\left(\mathfrak{z}_{0}\right)=H_{1}\left(\chi_{q}\left(\mathfrak{z}_{1}\right)\right)=\frac{q\chi_{q}\left(\mathfrak{z}_{1}\right)+1}{2}
\end{equation}
and so:
\begin{equation}
\chi_{q}\left(\mathfrak{z}_{1}\right)=\frac{2x_{0}-1}{q}
\end{equation}
Since $\chi_{q}$ is a map from $\mathbb{Z}_{2}$ to $\mathbb{Z}_{q}$,
$\left(2x_{0}-1\right)/q$ must be a $q$-adic integer. This forces
$2x_{0}-1\overset{q}{\equiv}0$ (else, $\chi_{q}\left(\mathfrak{z}_{1}\right)\in\mathbb{Q}_{q}\backslash\mathbb{Z}_{q}$). 

Now, let $\mathfrak{z}_{-1}\in\mathbb{Z}_{2}$ satisfy $\theta_{2}\left(\mathfrak{z}_{-1}\right)=\mathfrak{z}_{0}$.
Then, we can write $\mathfrak{z}_{-1}=\left[\mathfrak{z}_{-1}\right]_{2}+2\mathfrak{z}_{0}$.
As such, define $x_{-1}\overset{\textrm{def}}{=}\chi_{q}\left(\mathfrak{z}_{-1}\right)$.
Then:
\begin{align}
x_{-1} & =\chi_{q}\left(\mathfrak{z}_{-1}\right)\\
 & =\chi_{q}\left(2\mathfrak{z}_{0}+\left[\mathfrak{z}_{-1}\right]_{2}\right)\\
 & =H_{\left[\mathfrak{z}_{-1}\right]_{2}}\left(\chi_{q}\left(\mathfrak{z}_{0}\right)\right)\\
 & =H_{\left[\mathfrak{z}_{-1}\right]_{2}}\left(x_{0}\right)\\
 & =\begin{cases}
\frac{x_{0}}{2} & \textrm{if }\mathfrak{z}_{-1}\overset{2}{\equiv}0\\
\frac{qx_{0}+1}{2} & \textrm{if }\mathfrak{z}_{-1}\overset{2}{\equiv}1
\end{cases}
\end{align}
If $x_{0}\in\mathbb{Z}$, one and only one of $x_{0}/2$ and $\left(qx_{0}+1\right)/2$
will be an element of $\mathbb{Z}$. This occurs precisely when $\left[\mathfrak{z}_{-1}\right]_{2}=\left[x_{0}\right]_{2}$.
That is to say:
\begin{equation}
\mathfrak{z}_{-1}=\left[x_{0}\right]_{2}+2\mathfrak{z}_{0}
\end{equation}
In this manner, we can define $\mathfrak{z}_{-n}$ recursively for
all $n\geq1$ by declaring:
\begin{align}
\mathfrak{z}_{-n} & \overset{\textrm{def}}{=}\left[x_{-n+1}\right]_{2}+2\mathfrak{z}_{-n+1}\\
x_{-n} & \overset{\textrm{def}}{=}\chi_{q}\left(\mathfrak{z}_{-n}\right)
\end{align}
for all $n\geq1$.

In this way, we have shown the following:
\begin{prop}
\label{prop:Insights for the divergent case}Let $\mathfrak{z}_{0}\in\mathbb{Z}_{2}$
and let $x_{0}\overset{\textrm{def}}{=}\chi_{q}\left(\mathfrak{z}_{0}\right)$,
and, for each $n\geq1$, define $\mathfrak{z}_{n}\overset{\textrm{def}}{=}\theta_{2}^{\circ n}\left(\mathfrak{z}_{0}\right)$
and $x_{n}\overset{\textrm{def}}{=}\chi_{q}\left(\mathfrak{z}_{n}\right)$.
If $\mathfrak{z}_{0}\in\mathbb{Z}_{2}\backslash\mathbb{Q}$ and $x_{0}\in\mathbb{Z}$,
define:
\begin{align}
\mathfrak{z}_{-n} & \overset{\textrm{def}}{=}\left[x_{-n+1}\right]_{2}+2\mathfrak{z}_{-n+1}\\
x_{-n} & \overset{\textrm{def}}{=}\chi_{q}\left(\mathfrak{z}_{-n}\right)
\end{align}
for all $n\geq1$. Then, for all $n\in\mathbb{Z}$:

\vphantom{}

I. $x_{n}\in\mathbb{Z}$.

\vphantom{}

II. $T_{q}\left(x_{n+1}\right)=x_{n}$.

\vphantom{}

III. $x_{n}=H_{\left[\mathfrak{z}_{n}\right]_{2}}\left(x_{n+1}\right)$.

\vphantom{}

IV $\theta_{2}\left(\mathfrak{z}_{n}\right)=\mathfrak{z}_{n+1}$.

\vphantom{}

V. $\mathfrak{z}_{n}=\sum_{k=0}^{\infty}\left[x_{n+k+1}\right]_{2}2^{k}$
\end{prop}
Letting $n\geq0$, we have that the $x_{-n}$s is the $n$th term
of the forward orbit of $x_{0}$, while $x_{n}$ is a specific element
of the $n$th pre-image of $x_{0}$ under $T_{q}$. In this respect,
we have a sort of conveyer belt, where shifting the digits of $\mathfrak{z}_{0}$
to the left corresponds to applying an inverse image of $T_{q}$ to
$x_{0}$. On the other hand, if we shift the digits of $\mathfrak{z}_{0}$
to the right, we would need to pick the new left-most digit of the
shifted $\mathfrak{z}_{0}$, however, as the above shows, our choice
is completely determined by the parity of $x_{-1}=T_{q}\left(x_{0}\right)$.
(V) is of particular interest, because it shows that if $\mathfrak{z}_{0}$
is irrational and makes $\chi_{q}\left(\mathfrak{z}_{0}\right)$ into
the integer $x_{0}$, the $x_{n}$s for $n\geq1$ are completely determined
by $\mathfrak{z}_{0}$. These $x_{n}$s correspond to a particular
choice of a backward orbit of $x_{0}$. Indeed, the $n$th $2$-adic
digit of $\mathfrak{z}_{0}$ is the index ($0$ or $1$) of the $n$th
branch of $T_{q}$ (either $H_{0}$ or $H_{1}$) which we use to go
from $x_{n}$ (an $n$th pre-image of $x_{0}$) to $x_{n+1}$ (an
$\left(n+1\right)$th pre-image of $x_{0}$).

\subsection{\label{subsec:Implications-for-Transcendental}Implications for Transcendental
Number Theory}

Any proof of the \textbf{Weak Collatz Conjecture} (\textbf{Conjecture
\ref{conj:Weak Collatz Conjecture}} on page \pageref{conj:Weak Collatz Conjecture})
will necessarily entail a significant advancement in transcendental
number theory. A proof of the Weak Collatz Conjecture would necessarily
yield a proof of \textbf{Baker's Theorem} far simpler than any currently
known method \cite{Tao Blog}. Baker's Theorem, recall, concerns lower
bounds on the absolute values of linear forms of logarithms, which
are expressions of the form: 
\begin{equation}
\beta_{1}\ln\alpha_{1}+\cdots+\beta_{N}\ln\alpha_{N}\label{eq:linear form in logarithm}
\end{equation}
where the $\beta_{n}$s and $\alpha_{n}$s are complex algebraic numbers
(with all of the $\alpha_{n}$s being non-zero). Most proofs of Baker's
Theorem employ a variant of what is known as Baker's Method, in which
one constructs of an analytic function (called an \textbf{auxiliary
function})\textbf{ }with a large number of zeroes of specified degree
so as to obtain contradictions on the assumption that (\ref{eq:linear form in logarithm})
is small in absolute value. The import of Baker's Theorem is that
it allows one to obtain lower bounds on expressions such as $\left|2^{m}-3^{n}\right|$,
where $m$ and $n$ are positive integers, and it was for applications
such as these in conjunction with the study of Diophantine equations
that Baker earned the Fields Medal. See \cite{Baker's Transcendental Number Theory}
for a comprehensive account of the subject; also, parts of \cite{Cohen Number Theory}.

With $\chi_{q}$ and the \textbf{CP }at our disposal, we can get a
glimpse at the kinds of advancements in transcendental number theory
that might be needed in order to resolve the Weak Collatz Conjecture.
To begin, it is instructive to consider the following table of values
for $\chi_{q}$ and related functions in the case of the Shortened
$qx+1$ maps.
\begin{center}
\begin{table}
\begin{centering}
\begin{tabular}{|c|c|c|c|c|c|}
\hline 
$n$ & $\chi_{q}\left(n\right)$ & $\chi_{q}\left(B_{2}\left(n\right)\right)$ & $\chi_{3}\left(B_{2}\left(n\right)\right)$ & $\chi_{5}\left(B_{2}\left(n\right)\right)$ & $\chi_{7}\left(B_{2}\left(n\right)\right)$\tabularnewline
\hline 
\hline 
$0$ & $0$ & $0$ & $0$ & $0$ & $0$\tabularnewline
\hline 
$1$ & $\frac{1}{2}$ & $\frac{1}{2-q}$ & $-1$ & $-\frac{1}{3}$ & $-\frac{1}{5}$\tabularnewline
\hline 
$2$ & $\frac{1}{4}$ & $\frac{1}{4-q}$ & $1$ & $-1$ & $-\frac{1}{4}$\tabularnewline
\hline 
$3$ & $\frac{2+q}{4}$ & $\frac{2+q}{4-q^{2}}$ & $-1$ & $-\frac{1}{3}$ & $-\frac{1}{5}$\tabularnewline
\hline 
$4$ & $\frac{1}{8}$ & $\frac{1}{8-q}$ & $\frac{1}{5}$ & $\frac{1}{3}$ & $1$\tabularnewline
\hline 
$5$ & $\frac{4+q}{8}$ & $\frac{4+q}{8-q^{2}}$ & $-7$ & $-\frac{9}{17}$ & $-\frac{11}{41}$\tabularnewline
\hline 
$6$ & $\frac{2+q}{8}$ & $\frac{2+q}{8-q^{2}}$ & $-5$ & $-\frac{7}{17}$ & $-\frac{9}{41}$\tabularnewline
\hline 
$7$ & $\frac{4+2q+q^{2}}{8}$ & $\frac{4+2q+q^{2}}{8-q^{3}}$ & $-1$ & $-\frac{1}{3}$ & $-\frac{67}{335}$\tabularnewline
\hline 
$8$ & $\frac{1}{16}$ & $\frac{1}{16-q}$ & $\frac{1}{13}$ & $\frac{1}{11}$ & $\frac{1}{9}$\tabularnewline
\hline 
$9$ & $\frac{8+q}{16}$ & $\frac{8+q}{16-q^{2}}$ & $\frac{11}{7}$ & $-\frac{13}{9}$ & $-\frac{5}{11}$\tabularnewline
\hline 
$10$ & $\frac{4+q}{16}$ & $\frac{4+q}{16-q^{2}}$ & $1$ & $-1$ & $-\frac{1}{3}$\tabularnewline
\hline 
$11$ & $\frac{8+4q+q^{2}}{16}$ & $\frac{8+4q+q^{2}}{16-q^{3}}$ & $-\frac{29}{11}$ & $-\frac{53}{109}$ & $-\frac{85}{327}$\tabularnewline
\hline 
$12$ & $\frac{2+q}{16}$ & $\frac{2+q}{16-q^{2}}$ & $\frac{5}{7}$ & $-\frac{7}{9}$ & $-\frac{3}{11}$\tabularnewline
\hline 
$13$ & $\frac{8+2q+q^{2}}{16}$ & $\frac{8+2q+q^{2}}{16-q^{3}}$ & $-\frac{23}{11}$ & $-\frac{43}{109}$ & $-\frac{71}{327}$\tabularnewline
\hline 
$14$ & $\frac{4+2q+q^{2}}{16}$ & $\frac{4+2q+q^{2}}{16-q^{3}}$ & $\frac{19}{7}$ & $-\frac{39}{109}$ & $-\frac{67}{327}$\tabularnewline
\hline 
$15$ & $\frac{8+4q+2q^{2}+q^{3}}{16}$ & $\frac{8+4q+2q^{2}+q^{3}}{16-q^{4}}$ & $-1$ & $-\frac{1}{3}$ & $-\frac{149}{786}$\tabularnewline
\hline 
\end{tabular}
\par\end{centering}
\caption{Values of $\chi_{q}\left(n\right)$ and related functions}
\end{table}
\par\end{center}
\begin{example}
As per the Correspondence Principle,\textbf{ }note that the integer
values attained by $\chi_{3}\left(B_{2}\left(n\right)\right)$ and
$\chi_{5}\left(B_{2}\left(n\right)\right)$ are all periodic points
of the maps $T_{3}$ and $T_{5}$, respectively; this includes fixed
points at negative integers, as well. Examining $\chi_{q}\left(B_{2}\left(n\right)\right)$,
we see certain patterns, such as:
\begin{equation}
\chi_{q}\left(B_{2}\left(2^{n}-1\right)\right)=\frac{1}{2-q},\textrm{ }\forall n\in\mathbb{N}_{1}
\end{equation}
More significantly, $\chi_{3}\left(B_{2}\left(n\right)\right)$ appears
to be more likely to be positive than negative, whereas the opposite
appears to hold true for $q=5$ (and, heuristically, for all $q\geq5$).
Of special interest is: 
\begin{equation}
\chi_{q}\left(B_{2}\left(10\right)\right)=\frac{4+q}{16-q^{2}}=\frac{1}{4-q}
\end{equation}
\end{example}
By Version 1 of the CP (\textbf{Theorem \ref{thm:CP v1}}), every
cycle $\Omega\subseteq\mathbb{Z}$ of $T_{q}$ with $\left|\Omega\right|\geq2$
contains an integer $x$ of the form: 
\begin{equation}
x=\chi_{q}\left(B_{2}\left(n\right)\right)=\frac{\chi_{q}\left(n\right)}{1-\frac{q^{\#_{1}\left(n\right)}}{2^{\lambda_{2}\left(n\right)}}}=\frac{2^{\lambda_{2}\left(n\right)}\chi_{q}\left(n\right)}{2^{\lambda_{2}\left(n\right)}-q^{\#_{1}\left(n\right)}}\label{eq:Rational expression of odd integer periodic points of ax+1}
\end{equation}
for some $n\in\mathbb{N}_{1}$. In fact, every periodic point $x$
of $T_{q}$ in the odd integers can be written in this form. As such
$\left|2^{\lambda_{2}\left(n\right)}-q^{\#_{1}\left(n\right)}\right|=1$
is a \emph{sufficient condition }for $\chi_{q}\left(B_{2}\left(n\right)\right)$
to be a periodic point of $T_{q}$. However, as the $n=10$ case shows,
this is not\emph{ }a \emph{necessary }condition: there can be values
of $n$ where $2^{\lambda_{2}\left(n\right)}-q^{\#_{1}\left(n\right)}$
is large in archimedean absolute value, and yet nevertheless divides
the numerator on the right-hand side of (\ref{eq:Rational expression of odd integer periodic points of ax+1}),
thereby reducing $\chi_{q}\left(B_{2}\left(n\right)\right)$ to an
integer. In fact, thanks to P. Mih\u{a}ilescu's resolution of \textbf{Catalan's
Conjecture}, it would seem that Baker's Method-style estimates on
the archimedean\emph{ }size of $2^{\lambda_{2}\left(n\right)}-q^{\#_{1}\left(n\right)}$
will be of little use in understanding $\chi_{q}\left(B_{2}\left(n\right)\right)$:
\begin{thm}[\textbf{Mih\u{a}ilescu's Theorem}\footnote{Presented in \cite{Cohen Number Theory}.}]
The only choice of $x,y\in\mathbb{N}_{1}$ and $m,n\in\mathbb{N}_{2}$
for which: 
\begin{equation}
x^{m}-y^{n}=1\label{eq:Mihailescu's Theorem}
\end{equation}
are $x=3$, $m=2$, $y=2$, $n=3$ (that is, $3^{2}-2^{3}=1$). 
\end{thm}
With Mih\u{a}ilescu's Theorem, it is easy to see that for any odd
integer $q\geq3$, $\left|q^{\#_{1}\left(n\right)}-2^{\lambda_{2}\left(n\right)}\right|$
will never be equal to $1$ for any $n\geq8$, because the exponent
of $2$ (that is, $\lambda_{2}\left(n\right)$) will be $\geq4$ for
all $n\geq8$ (any such $n$ has $4$ or more binary digits). Consequently,
for any odd prime $q$, if $n\geq8$ makes $\chi_{q}\left(B_{2}\left(n\right)\right)$
into a rational integer (and hence, a periodic point of $T_{q}$),
it \emph{must }be that the numerator $2^{\lambda_{2}\left(n\right)}\chi_{q}\left(n\right)$
in (\ref{eq:Rational expression of odd integer periodic points of ax+1})
is a multiple of $2^{\lambda_{2}\left(n\right)}-q^{\#_{1}\left(n\right)}$.

In light of this, rather than the Archimedean/Euclidean \emph{size
}of $2^{\lambda_{2}\left(n\right)}-q^{\#_{1}\left(n\right)}$, it
appears we must study its \emph{multiplicative }structure, as well
as that of $2^{\lambda_{2}\left(n\right)}\chi_{q}\left(n\right)$\textemdash which
is to say, the set of these integers' prime divisors, and how this
set depends on $n$. In other words, we ought to study the $p$-adic
absolute values of $\left|\chi_{q}\left(n\right)\right|_{p}$ and
$\left|1-M_{q}\left(n\right)\right|_{p}$ for various values of $n\geq1$
and primes $p$. It may also be of interest to use $p$-adic methods
(Fourier series, Mellin transform) to study $\left|\chi_{q}\left(\mathfrak{z}\right)\right|_{p}$
as a real-valued function of a $2$-adic integer variable for varying
values of $p$.

\section{\label{sec:Corrections}Corrections}

As stated repeatedly throughout the paper, the Correspondence Principle
for Divergent Points is false: the existence of irrational $2$-adic
integers $\mathfrak{z}\in\mathbb{Z}_{2}\backslash\mathbb{Q}$ for
which $\chi_{q}\left(\mathfrak{z}\right)\in\mathbb{Z}$ \emph{does
not }guarantee that the integer $\chi_{q}\left(\mathfrak{z}\right)$
is a divergent point of $T_{q}$. The author owes credit to Zhongshu
Wu and Wu's use of ChatGPT for the following counterexample.
\begin{example}
Consider the following reverse orbit of $1$ under $T_{3}$:
\begin{equation}
1\leftarrow2\leftarrow4\leftarrow8\leftarrow5\leftarrow10\leftarrow20\leftarrow40\leftarrow80\leftarrow53\leftarrow\cdots
\end{equation}
More generally, starting with any integer $x_{0}$ which is not a
multiple of $3$ (for the above example, $x_{0}=1$), set:
\begin{equation}
x_{k+1}=\frac{2^{a_{k}+1}x_{k}-1}{3},\textrm{ }\forall k\geq0\label{eq:x_k plus 1}
\end{equation}
where the $a_{k}$s are any sequence of positive integers with $a_{k+1}>a_{k}$
for all $k\geq0$ so that $2^{a_{k}}x_{k}$ is congruent to $8$ mod
$9$. (Note that such a choice of $a_{k}$s is possible only if the
$x_{k}$s are not multiples of $3$.) Moreover, since $2$ is a primitive
root mod $9$, whenever $x_{k}$ is not a multiple of $3$, there
will always be arbitrarily large $a_{k}$ so that $2^{a_{k}}x_{k}$
is congruent to $8$ mod $9$. Moreover, as constructed, we have that
(\ref{eq:x_k plus 1}) guarantees that $x_{k+1}$ is congruent to
$2$ mod $3$

Next, set $s_{0}=0$ and define:
\[
s_{k+1}=s_{k}+a_{k}+1,\textrm{ }\forall k\geq0
\]
Then, for each $k\geq0$ and each integer $r$ with $0\leq r\leq a_{k}$,
let: 
\begin{equation}
d_{s_{k}+r}=\begin{cases}
0 & \textrm{if }0\leq r<a_{k}\\
1 & \textrm{if }r=a_{k}
\end{cases}
\end{equation}
and define a $2$-adic integer $\mathfrak{z}_{0}\in\mathbb{Z}_{2}$
by:
\begin{equation}
\mathfrak{z}_{0}=\sum_{n=0}^{\infty}d_{n}2^{n}
\end{equation}
Equivalently, we have:
\begin{equation}
\mathfrak{z}_{0}=2^{a_{0}}+2^{a_{0}+a_{1}+1}+2^{a_{0}+a_{1}+a_{2}+2}+\cdots
\end{equation}
Since the $a_{k}$s are strictly increasing, the $d_{k}$s consist
of $1$s separated by increasingly large numbers of $0$s, which guarantees
that $\mathfrak{z}_{0}$ is irrational. 

Now, set $R_{0}=x_{0}$, and let: 
\begin{equation}
R_{n+1}=\begin{cases}
2R_{n} & \textrm{if }d_{n}=0\\
\frac{2R_{n}-1}{3} & \textrm{if }d_{n}=1
\end{cases}
\end{equation}
Then:
\begin{equation}
\chi_{3}\left(\left[\mathfrak{z}_{0}\right]_{2^{N}}\right)=\sum_{n=0}^{N-1}d_{n}\frac{3^{\#_{1}\left(\left[\mathfrak{z}_{0}\right]_{2^{n}}\right)}}{2^{n+1}}
\end{equation}
where $\#_{1}\left(\left[\mathfrak{z}_{0}\right]_{2^{n}}\right)$
is the sum of $d_{0}+d_{1}+\cdots+d_{n-1}$. Because of this, we have:
\begin{equation}
x_{0}-\chi_{3}\left(\left[\mathfrak{z}_{0}\right]_{2^{N}}\right)=\frac{3^{\#_{1}\left(\left[\mathfrak{z}_{0}\right]_{2^{N}}\right)}}{2^{N}}R_{N}
\end{equation}
As constructed, the $R_{N}$s are all integers. Hence, upon taking
$3$-adic absolute values:
\[
\left|x_{0}-\chi_{3}\left(\left[\mathfrak{z}_{0}\right]_{2^{N}}\right)\right|_{3}=\left|\frac{3^{\#_{1}\left(\left[\mathfrak{z}_{0}\right]_{2^{N}}\right)}}{2^{N}}R_{N}\right|_{3}\leq3^{-\#_{1}\left(\left[\mathfrak{z}_{0}\right]_{2^{N}}\right)}
\]
Since $\mathfrak{z}_{0}$ has infinitely many $1$s digits, the upper
bound tends to $0$ as $N\rightarrow\infty$, the rising-continuity
of $\chi_{3}$ proves that $\chi_{3}\left(\mathfrak{z}_{0}\right)$
converges $3$-adically to $x_{0}$. This shows that for each integer
$x_{0}$ co-prime to $3$, there is an irrational $\mathfrak{z}_{0}$
at which $\lim_{n\rightarrow\infty}\chi_{3}\left(\left[\mathfrak{z}_{0}\right]_{2^{n}}\right)$
converges $3$-adically to $x_{0}$. This contradicts the Correspondence
Principle for Divergent Points (\textbf{Theorem \ref{thm:Divergent trajectories come from irrational z})},
which asserts that $\chi_{3}\left(\mathfrak{z}_{0}\right)\in\mathbb{Z}$
for $\mathfrak{z}_{0}\in\mathbb{Z}_{2}\backslash\mathbb{Q}$ implies
$\chi_{3}\left(\mathfrak{z}_{0}\right)$ is a divergent point of $T_{3}$.
Indeed, setting $x_{0}=1$ shows that even a periodic point can be
a value attained by $\chi_{3}$ at an irrational $2$-adic input. 
\end{example}
While this invalidates the CPDP, the rest of this paper's findings
still hold. That being said, the lack of the CPDP also has implications
for the Wiener Tauberian approach to the study of the value distribution
of $\chi_{q}$. This will be discussed in Part III of this series
of papers.

\end{document}